\documentclass[reqno]{amsart}
\usepackage{epsfig,amsmath,amsfonts,latexsym}
\usepackage{amsthm}
\usepackage[abs]{overpic}
\usepackage[usenames,dvipsnames]{xcolor}
\usepackage{palatino}
\usepackage[hyperfootnotes=false]{hyperref}
\usepackage{caption}
\usepackage{dsfont}
\usepackage{mathtools}
\usepackage{cite}
\usepackage{tikz-cd}
\usepackage{amssymb}
\usepackage{comment}
\hypersetup{
  colorlinks,
  citecolor=Blue,
  linkcolor=Black,
  urlcolor=Green}
  
\theoremstyle{plain}

\newtheorem{theorem}{Theorem}
\newtheorem{theoremM}{Theorem}

\newtheorem{dummy}{anything}[section]

\newtheorem{lemma}[dummy]{Lemma}

\newtheorem{proposition}[dummy]{Proposition}

\newtheorem{corollary}[dummy]{Corollary}

\theoremstyle{definition}
\newtheorem{definition}[dummy]{Definition}

\newtheorem{example}[dummy]{Example}

\newtheorem{remark}[dummy]{Remark}

\newtheorem*{assumptions}{Assumptions}
\newtheorem*{assumption}{Assumption}
\theoremstyle{remark}
\newcommand{\R}{\mathbb{R}}
\newcommand{\C}{\mathbb{C}}

\DeclareMathOperator{\re}{re}

\DeclareMathOperator{\im}{im}

\DeclareMathOperator{\atan}{atan}

\def\R{\mathbb{R}}

\def\C{\mathbb{C}}

\title{Global hypersurfaces of section in the spatial restricted three-body problem}
\author{Agustin Moreno, Otto van Koert}

\address[A.\ Moreno]{\emph{Current:} School of Mathematics\\ Institute for Advanced Study\\ Princeton NJ\\ USA \newline \emph{Previous:} Department of Mathematics \\ Uppsala University \\  Uppsala  \\ Sweden}

\email{agustin.moreno2191@gmail.com}

\address[O.\ van Koert]{Department of Mathematical Sciences and Research Institute of Mathematics, Seoul National University\\
Building 27, room 402\\
San 56-1, Sillim-dong, Gwanak-gu, Seoul, South Korea\\
Postal code 08826 }

\email{okoert@snu.ac.kr}

\date{}

\begin{document}

\begin{abstract}
We propose a contact-topological approach to the \emph{spatial} circular restricted three-body problem, for energies below and slightly above the first critical energy value. We prove the existence of a circle family of global \emph{hypersurfaces} of section for the regularized dynamics.
Below the first critical value, these hypersurfaces are diffeomorphic to the unit disk cotangent bundle of the $2$-sphere, and they carry symplectic forms on their interior, which are each deformation equivalent to the standard symplectic form.
The boundary of the global hypersurface of section is an invariant set for the regularized dynamics that is equal to a level set of the Hamiltonian describing the regularized planar problem.
The first return map is Hamiltonian, and restricts to the boundary as the time-$1$ map of a positive reparametrization of the Reeb flow in the planar problem. This construction holds for any choice of mass ratio, and is therefore \emph{non-perturbative}. We illustrate the technique in the completely integrable case of the rotating Kepler problem, where the return map can be studied explicitly.
\end{abstract}
\maketitle

{\centering\footnotesize  \textit{C'est avec l'intuition qu'on trouve, c'est avec la logique qu'on prouve.} H. Poincar\'e.\par}

\tableofcontents

\section{Introduction}

In this article, we will discuss global hypersurfaces of section for the well-known circular restricted spatial three-body problem.
This problem concerns the motion of a massless particle in $\R^3$ under the influence of two heavy primaries with mass $\mu$ and $1-\mu$, where $\mu\in(0,1)$ is a mass parameter. By a time-dependent rotation, we can fix these primaries at $\vec m=(\mu-1,0,0)$, which we will call Moon, and $\vec e=(\mu,0,0)$, which we will call Earth.
In the setting of symplectic geometry, the restricted three-body problem is then most easily described as the Hamiltonian dynamics of the following Hamiltonian on $(T^*\R^3\setminus \{ \vec m,\vec e\},d\vec p \wedge d \vec q)$:
$$
H(\vec q,\vec p)=\frac{1}{2}\Vert \vec p\Vert^2 - \frac{\mu}{\Vert \vec q-\vec m\Vert } - \frac{1-\mu}{\Vert \vec q-\vec e\Vert } +p_1q_2-p_2q_1. 
$$
The planar case of this problem is obtained by setting $q_3=p_3=0$. 

The Hamiltonian flow of this dynamical system has singularities caused by two-body collisions, namely collisions of the massless particle with $\vec m$ and with $\vec e$.
The resulting flow can be extended across these singularities using various schemes. We will use Moser-regularization, \cite{moser}, to do so.

\subsection*{Poincar\'e--Birkhoff theorem, and the planar three-body problem.} The problem of finding closed orbits in the planar case goes back to ground-breaking work in celestial mechanics of Poincar\'e \cite{P87,P12}, building on work of G.W.\ Hill on the lunar problem \cite{H78}. The basic scheme for his approach may be reduced to:
\begin{itemize}
    \item[(1)] Finding a \emph{global surface of section} for the dynamics of the regularized problem;
    \item[(2)] Proving a fixed point theorem for the resulting first return map.
\end{itemize}
This is the setting for the celebrated Poincar\'e-Birkhoff theorem.
%, proposed and confirmed in special cases by Poincar\'e and later proved in full generality by Birkhoff in \cite{Bi13}. 
%The statement can be summarized as: if $f: A\rightarrow A$ is an area-preserving homeomorphism of the annulus $A=[-1,1]\times S^1$ that satisfies a \emph{twist} condition at the boundary, then it admits infinitely many periodic points of arbitrary large period. 
%In the case where the annulus $A$ arises as a global surface of section for some dynamics in a $3$-manifold $M$ (i.e.\ $\partial A$ consists of closed orbits, and the orbit of every point in $M\backslash \partial A$ meets int$(A)$ in the future and past), and $f$ is the associated return map mapping a point to its first return point, the periodic points of $f$ correspond to closed orbits for the dynamics. 
In this paper, we address the above first step in the spatial case; the second will be addressed in \cite{MvK}, where we prove a generalized version of the Poincar\'e-Birkhoff theorem.

\subsection*{Moser regularization.}
We denote by $L_1$ the critical point of the Jacobi Hamiltonian $H$ with the smallest critical value, and by $p_{\R^3}:T^*\R^3\to \R^3$ the projection to the $\vec q$-coordinates.
Let us fix an energy level $c<0$, and consider the component $S_c\subset H^{-1}(c)$ with the property that $\overline{ p_{\R^3}(S_c)}$ contains $\vec m$: this is the component of the level set containing the Moon.
As explained in Section~\ref{sec:Moserregularization}, Moser regularization applied to this setting gives us a smooth Hamiltonian $Q_{\mu,c}$ on 
$$
T^*S^3=\{ (\xi,\eta) \in T^*\R^4 ~|~\Vert \xi \Vert^2=1,~
\langle \xi , \eta \rangle =0 \}
$$ 
with the following property.
The level set $Q_{\mu,c}^{-1}(\mu^2/2)$ contains a component $\Sigma_c$ that projects to $S_c$ under stereographic projection, and the Hamiltonian dynamics of $Q_{\mu,c}$ are reparametrization of the Hamiltonian dynamics of $H$.
This procedure extends the dynamics across collisions with the Moon. We may also do the same for the Earth.

The topology is as follows.
For $c <H(L_1)$, $\Sigma_c$ is diffeomorphic to the unit cotangent bundle of $S^3$.
We refer to this energy range as the \emph{low-energy range}.
If $c$ is in the interval $(H(L_1),H(L_2))$, where $L_2$ is the critical point with second lowest critical value, then $\overline{ p_{\R^3}(S_c)}$ contains both the Earth and the Moon, and we need to perform regularization in both points. We also denote the doubly regularized level set by $\Sigma_c$.
This doubly regularized hypersurface $\Sigma_c$ is diffeomorphic to the connected sum of two unit cotangent bundles of $S^3$.

\subsection*{Statement of results} 
We first recall the concept of a global hypersurface of section. On an oriented smooth manifold $\Sigma$, we consider a flow $\varphi_t$ of an autonomous vector field, which we assume to be nowhere vanishing.

\begin{definition}
  A \emph{global hypersurface of section} for $\varphi_t$ consists of an embedded, compact, oriented hypersurface $P$ with the following properties:
  \begin{enumerate}
      \item The boundary $B=\partial P$, if non-empty, is an invariant set for $\varphi_t$;
      \item $\varphi_t$ is positively transverse to $P\backslash B$;
      \item For all $x\in \Sigma \setminus B$, there exist $\tau_+>0$ and $\tau_-<0$ such that $\varphi_{\tau_+}(x) \in \mbox{int}(P)$ and $\varphi_{\tau_-}(x) \in \mbox{int}(P)$.
  \end{enumerate}
\end{definition}

The main purpose of this object to reduce the dynamics of flows to the (discrete) dynamics of diffeomorphisms.
This concept has been used very fruitfully for 3-manifolds, where the boundary is necessarily empty or a collection of periodic orbits. Higher-dimensional invariant sets are hard to find in general, and usually don't have good stability properties, making this notion less ubiquitous for higher-dimensional dynamical systems. 
However, the spatial restricted three-body problem has special symmetries, and these allow us to prove the following result.

\begin{theorem}
\label{thm:gss_openbooks}
Fix a mass parameter $\mu \in [0,1]$, and let $\Sigma_c$ as above denote a connected component of the regularized, spatial, circular, restricted three-body problem for energy level $c$ that contains $\vec m$ or $\vec e$ in its projection.
We have the following:
\begin{itemize}
\item for $c<H(L_1)$, the set $\{ \xi_3=0 \}$ is a global hypersurface of section;
\item for $c \in (H(L_1),H(L_2))$, the manifold $\Sigma_c$ also admits a global hypersurface of section.\footnote{if $\mu=0,1$, then $H(L_1)=H(L_2)$, so the statement is empty in this case. For $c>H(L_2)$, the energy levels are necessarily non-compact, and satellites can escape in the unbounded component. }
\end{itemize} 
Moreover, in each of these cases there is an $S^1$-family of global hypersurfaces of section.  
\end{theorem}

\begin{remark}
The main advantage of the first description is that it is simple and explicit, and so lends itself well to numerical work.
Below we will investigate the underlying contact topology, but we point out that the above theorem does not require any contact topology; it also does not imply the contact condition.
\end{remark}

The above result can be rephrased and clarified using the language of contact topology.
We know from \cite{AFvKP,ChoKim} that the bounded components of the regularized energy hypersurfaces of the circular (planar and spatial) restricted three-body problem are of contact type for energy levels $c$ below and slightly above the first critical energy value, say up to $H(L_1)+\epsilon$.
We will now impose this extra assumption, namely $c<H(L_1)+\epsilon$.

The $S^1$-family of global hypersurfaces of section form the pages of a so-called \emph{open book decomposition} for the regularized $5$-dimensional energy level sets. 
Such a decomposition, for a closed manifold $\Sigma$, consists of a fiber bundle $\pi: \Sigma \backslash B \rightarrow S^1$, where $B\subset \Sigma$ is a codimension $2$ submanifold, known as \emph{binding}, which has a trivial normal bundle, such that $\pi$ coincides with the angular coordinate along some choice of neighbourhood $B\times \mathbb{D}^2$ of $B$. 
Topologically, this decomposition is also determined by the data of the \emph{page} $P$ (the closure of the typical fiber of $\pi$) with $\partial P=B$, and the \emph{monodromy} $\phi:P \rightarrow P$ of the fiber bundle $\pi$, satisfying $\phi=id$ near $B$. 
We denote $\Sigma=\mathbf{OB}(P,\phi)$ whenever $\Sigma$ admits an open book decomposition with data $(P,\phi)$: this is sometimes called an abstract open book.
Smoothly, we have
$$
\mathbf{OB}(P,\phi)\cong B \times \mathbb{D}^2 \bigcup_\partial \mbox{Map}(P,\phi), $$
where $\bigcup_\partial$ denotes the boundary union, and $\mbox{Map}(P,\phi)=P \times \R/(x,t)\sim (\phi(x),t+1)$ is the associated mapping torus. 
This manifold inherits a projection $\pi$ as above. 
We will denote the $\theta$-page by $P_\theta:=\overline{\pi^{-1}(\theta)}$.

A relationship between open books and hypersurfaces of section for a given flow is contained in the following notion. 
If $\varphi_t:\Sigma\rightarrow\Sigma$ is the time-$t$ flow of some autonomous vector field $X$ on $\Sigma$, we say that the open book decomposition is \emph{adapted to the dynamics} if each page $P_\theta$ of the open book decomposition is a global hypersurface of section for $\varphi_t$. 
This definition implies that $B$ is invariant under the flow. 

We now rephrase our theorem in contact-topological terms.

\begin{theoremM}\label{thm:openbooks}
Fix a mass ratio $\mu\in (0,1]$.
As above, denote a connected, bounded component of the regularized, spatial, circular restricted three-body problem for energy level $c$ by $\Sigma_c$. 
Then $\Sigma_c$ is of contact-type and admits a supporting open book decomposition for energies $c<H(L_1)$ that is adapted to the Hamiltonian dynamics of $Q_{\mu,c}$. 
Furthermore, if $\mu<1$, then there is $\epsilon>0$ such that the same holds for $c\in (H(L_1),H(L_1)+\epsilon)$.
The open books have the following abstract form:
\[
\Sigma_c
\cong
\begin{cases}
\textbf{OB}(\mathbb{D}^*S^2,\tau^2), & \mbox{if }c<H(L_1) \\
\textbf{OB}(\mathbb{D}^*S^2 \natural \mathbb{D}^*S^2,\tau_1^2 \circ \tau_2^2), & \mbox{if }c\in (H(L_1),H(L_1)+\epsilon) \text{ and }\mu<1.
\end{cases}
\]
Here, $\mathbb{D}^*S^2$ is the unit cotangent bundle of the $2$-sphere, $\tau$ is the positive Dehn-Seidel twist along the Lagrangian zero section $S^2\subset \mathbb{D}^*S^2$, and $\mathbb{D}^*S^2 \natural \mathbb{D}^*S^2$ denotes the boundary connected sum of two copies of $\mathbb{D}^*S^2$. 
The monodromy of the second open book is the composition of the square of the positive Dehn-Seidel twists along both zero sections (they commute). The binding is the planar problem: for  $c<H(L_1)$ this is $\Sigma_c^P\cong \mathbb{R}P^3$, and for $c\in (H(L_1), H(L_2)\,)$ we have $\Sigma_c^P\cong \mathbb{R}P^3 \# \mathbb{R}P^3$. 
\end{theoremM}

\begin{figure}
    \centering
    \includegraphics[width=0.4 \linewidth]{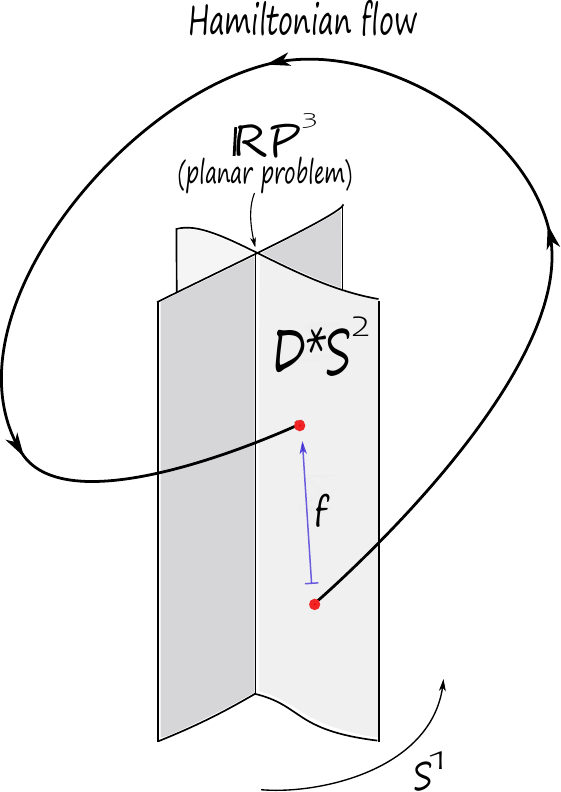}
    \caption{The open book for $\Sigma_c$, with $c<H(L_1)$, and the first return map $f$.}
    \label{fig:openbook}
\end{figure}

See Figure \ref{fig:openbook} for an abstract representation.

\smallskip

Let us draw analogies to the planar situation, discuss some history and speculate a bit: 
for $c<H(L_1)$, we have, smoothly, the following: $B=\Sigma_c^P\cong \mathbb{R}P^3\cong \mathbf{OB}(\mathbb{D}^*S^1,\tau_P^2)$, where $\tau_P$ is the positive Dehn twist along $S^1\subset \mathbb{D}^*S^1$, and one would hope that this open book is adapted to the planar dynamics, and that the return map is a Birkhoff twist map. 
Let us briefly review what is known about annular global surfaces of section in the planar restricted three-body problem.
For $c<H(L_1)$ and $\mu \sim 1$, one can interpret from this perspective that Poincar\'e \cite{P12} proved this by perturbing the rotating Kepler problem (when $\mu=1$, an integrable system; note that in other conventions this is $\mu=0$). 
In the case where $c\ll H(L_1)$ is very negative and $\mu \in (0,1]$, the existence of an annular global surface of section and the twist map condition was established by Conley \cite{C63} (also perturbatively). 
The most recent result concerning annular global surfaces of section is the result due to Hryniewicz, Salom\~ao and Wysocki in \cite[Theorem 1.18]{HSW}:
this result asserts the existence of an adapted open book in the case where $(\mu,c)$ lies in the \emph{convexity range}. This is a subset of the low-energy range, consisting of pairs $(\mu,c)$ for which the so-called Levi-Civita regularization is dynamically convex, \cite{AFFHvK}.

Disk-like global surfaces of section were found by McGehee, in \cite{M69}, for the rotating Kepler problem and small perturbations thereof, so $\mu \sim 1$. He computed the return map for $\mu=1$, and used KAM theory to establish the existence of invariant tori.
More recently, \emph{non-perturbative} holomorphic curve methods due to Hofer-Wysocki-Zehnder \cite{HWZ98} have been used to establish disk-like global surfaces of section in the convexity range, \cite{AFFHvK,AFFvK}. 

\begin{figure}
    \centering
    \includegraphics{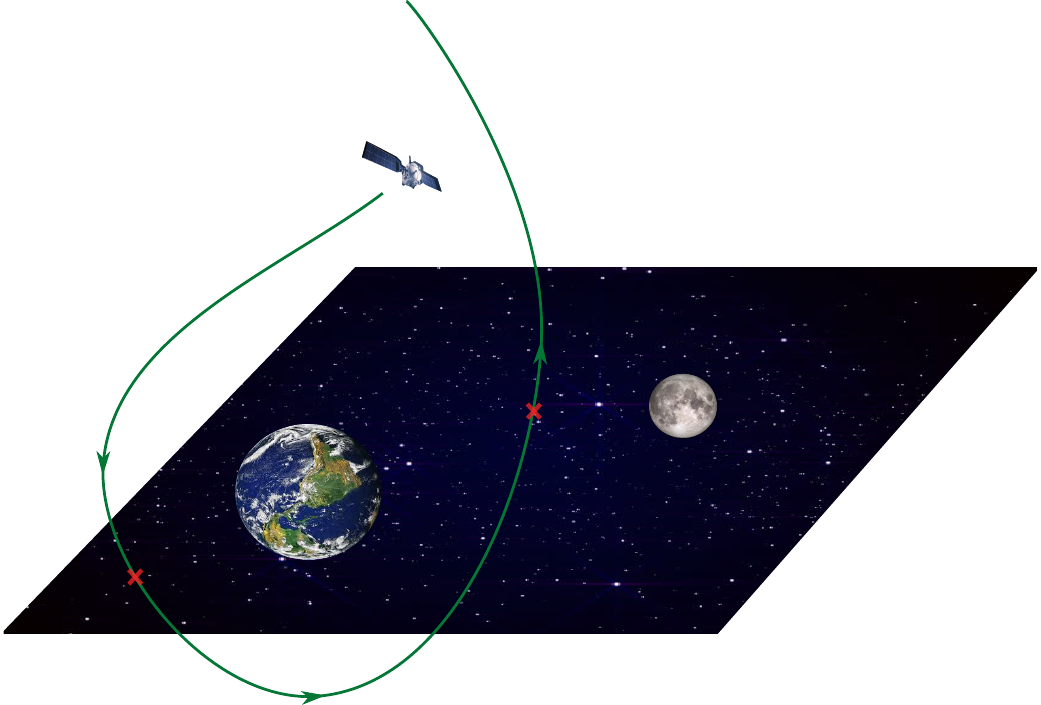}
    \caption{Theorem~\ref{thm:openbooks} admits a physical interpretation: away from collisions, the orbits of the negligible mass point intersect the plane containing the primaries transversely. This is intuitively clear from a physical perspective, and translates (after regularization) to the fact that the pages $\{q_3=0,p_3>0\}$, $\{q_3=0,p_3<0\}$ of the ``physical'' open book are global hypersurfaces of section outside of the collision locus. Unfortunately this does not extend continuously to the latter, since for instance there exist (regularized) collision orbits which are periodic and ``bounce'' vertically over a primary, always staying on the region $q_3>0$ (or $q_3<0$). This will be addressed by interpolating with the ``geodesic'' open book near the collision locus.}
    \label{fig:transversalitypic}
\end{figure}

We wish to emphasize that the results in this paper hold for $c$ in the whole low-energy range, independently of mass ratio, and even extend to higher energies. One partial reason is the following: while in the planar case finding a suitable invariant subset is non-trivial, one usable invariant subset in the spatial case is immediately obvious; it is the planar problem.

\begin{figure}
    \centering
    \includegraphics[width=0.7 \linewidth]{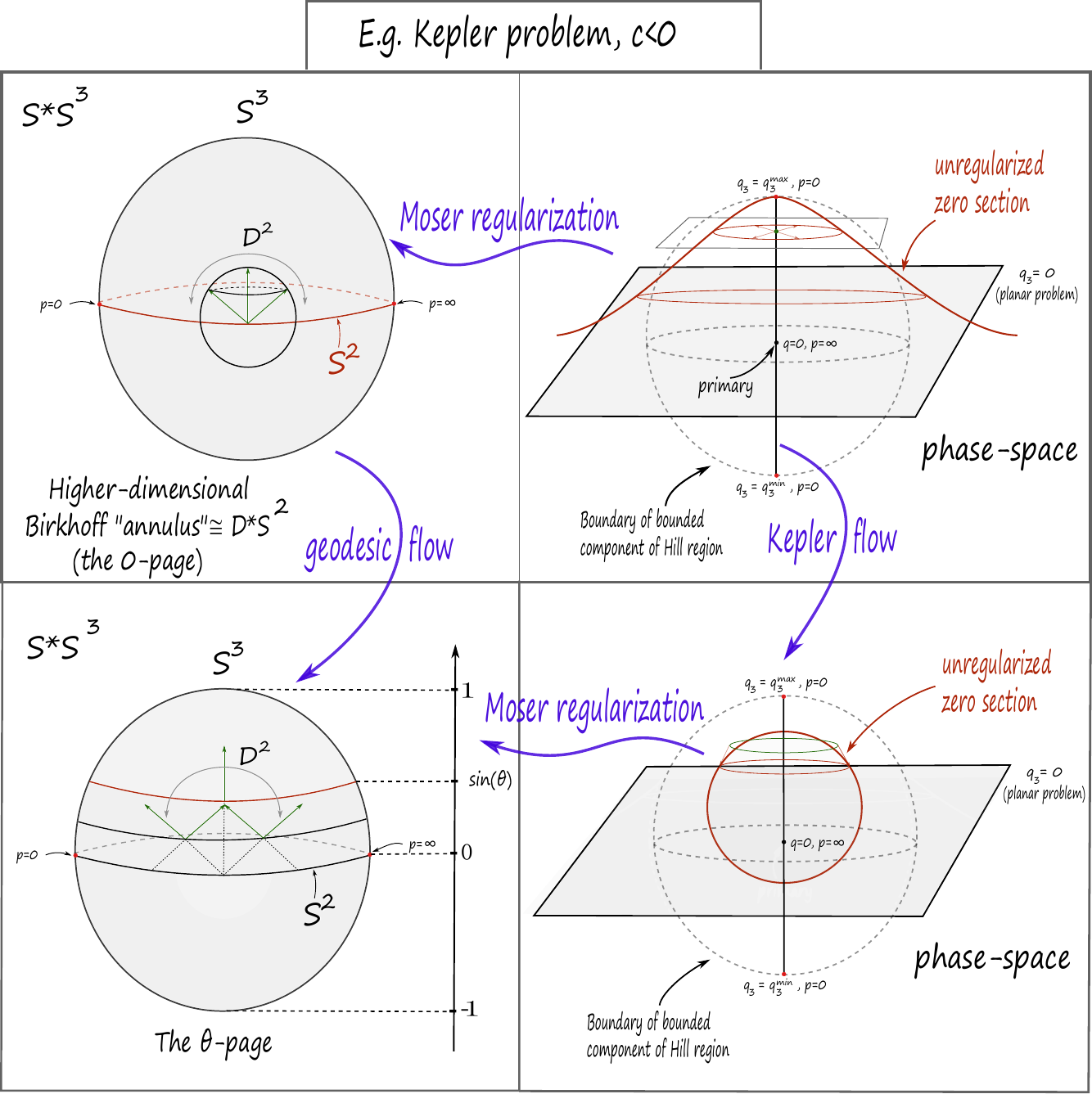}
    \caption{The page $\mathbb{D}^*S^2$ of the open book for $\mathbb S^*S^3$ with $q_3>0, p_3=0$ corresponds physically to orbits achieving their maximal height $q_3>0$ for fixed energy. In the Kepler problem, see Example~\ref{ex:Kepler}, the zero section of this page sits in phase space as a plane, emanating from $q_1=q_2=0, q_3=q_3^{max}$ (the highest allowed position for given energy), consisting of horizontal momenta which become infinite as $q_3 \rightarrow 0^+$ along the $q_3$-axis. This compactifies to a $2$-sphere after regularization so that $p=\infty$ becomes the north pole of $S^2\subset S^3$. The fiber over a horizontal momentum $p \in S^2$ sits as a $2$-disk of positions $q$ with $q_3\geq 0$ so that $(q,p)$ has given energy; for instance, the fiber over $p=0$ is the intersection of the boundary of the bounded Hill region with the upper half-plane $q_3\geq 0$, while the others are $2$-disks with boundary in $\{q_3=0\}$ emanating from the $q_3$-axis and which shrink down to a point at $q=0$. The other pages sit in phase-space as the image of the above one under the Kepler flow. Note that the $\theta+\pi$-page is obtained from the $\theta$-page by reflection along the equator in the regularized picture, and by reflection along the plane $\{q_3=0\}$ in the unregularized one. This is a general fact; see Proposition \ref{prop:symmetries}.}
    \label{fig:Birkhoff_annulus}
\end{figure}

Topologically and abstractly, the situation may be understood as follows: the Stein manifold $\mathbb{D}^*S^2$ carries a Lefschetz fibration structure, whose smooth fibers are the annuli $\mathbb{D}^*S^1$, and its monodromy is precisely $\tau^2_P$ along the vanishing cycle $S^1 \subset \mathbb{D}^*S^1$. 
We write $\mathbb{D}^*S^2=\mathbf{LF}(\mathbb{D}^*S^1,\tau_P^2)$. 
By restricting this Lefschetz fibration to the boundary, we obtain the above open book for $\mathbb{R}P^3$. 
See Figure \ref{fig:LF} in Appendix \ref{app:rotKepler}. 
The Lefschetz fibration on the pages $\mathbb{D}^*S^2$ gives $(\mathbb S^*S^3,\xi_{std})$ the structure of an \emph{iterated planar} contact $5$-manifold, which has been studied in \cite{Acu, AEO, AM18}. 
Similar remarks apply for $c \in (H(L_1),H(L_1)+\epsilon).$ 

%where 
%\begin{equation*}\begin{split}\mathbb{R}P^3\# \mathbb{R}P^3&=\mathbf{OB}(\mathbb{D}^*S^1 \natural \mathbb{D}^*S^1,\tau_{P,1}^2 \circ \tau_{P,2}^2)\\&=\partial \mathbf{LF}(\mathbb{D}^*S^1 \natural \mathbb{D}^*S^1,\tau_{P,1}^2 \circ \tau_{P,2}^2)\\&=\partial \left(\mathbb{D}^*S^2 \natural \mathbb{D}^*S^2\right).\end{split}\end{equation*}

Figure \ref{fig:Birkhoff_annulus} describes one of the pages of the ``geodesic'' open book (which is well-behaved with respect to collisions, as opposed to the more ``physical'' version discussed in Figure \ref{fig:transversalitypic}), in the simplest case of the Kepler problem, see Example~\ref{ex:Kepler}. 
This is a higher-dimensional version of the well-known Birkhoff annulus, and topologically consists of unit directions in $\mathbb S^*S^3$ which are positively transverse to the equator $S^2\subset S^3$. 
Its boundary is precisely the unit directions tangent to $S^2$ (i.e.\ $\mathbb{R}P^3$). 
The other pages are obtained by flowing with the geodesic flow for the round metric on $S^3$. By construction, this open book is adapted to the dynamics of the geodesic flow, which is the regularized version of the Kepler problem for $c<0$. 
In Figure \ref{fig:Birkhoff_annulus} we also sketch how the pages of this open book sit in phase space, for which one should recall that in Moser regularization momenta lie in the base, and positions lie in the fiber; the fiber over $p=\infty$ as a point in $S^2\subset S^3$ is a Legendrian $2$-sphere -the collision locus- which is collapsed to the point $q=0$ in unregularized coordinates (we will review this regularization method in Section \ref{sec:Moserregularization}). For the general case, the open book in Theorem~\ref{thm:openbooks} coincides with the physical open book away from collisions, and with the geodesic one on the collision locus (but where the unit cotangent bundle for the round metric is replaced with a low-energy level set of the appropriate Hamiltonian); see the proof of Lemma \ref{lemma:GSS_Stark_Zeeman}. One may think of a page as a Liouville filling of the planar problem; see Figure~\ref{fig:pageasfilling}. Here,  there is a slight subtlety due to the fact that the symplectic form on a page degenerates at the boundary, but the symplectic form may be modified by a continuous conjugation to make it a filling; see Appendix \ref{monodromyvsreturnmap}. 

The open books of Theorem~\ref{thm:openbooks} support the corresponding contact structure on $\Sigma_c$, in the sense of Giroux (see Definition \ref{def:Giroux}). One can interpret this as the smooth topology being adapted to the given geometry. However, unlike in our current situation, the contact form (and hence the dynamics) in the setting of Giroux's notion is never fixed; only the contact structure is. One then adapts the dynamics to the open book via a Giroux form, whose Reeb dynamics is normally taken as simple as possible. The content of Theorem~\ref{thm:openbooks} is stronger: the smooth topology is actually adapted to the given dynamics, which is a posteriori given by a Giroux form.

\medskip

\textbf{Symmetries.} Consider the symplectic involution of $(\mathbb{R}^6,dp\wedge dq)$ given by
$$
r: (q_1,q_2,q_3,p_1,p_2,p_3)\mapsto (q_1,q_2,-q_3,p_1,p_2,-p_3).
$$
We also have the anti-symplectic involutions
$$
\rho_1: (q_1,q_2,q_3,p_1,p_2,p_3) \mapsto (q_1,-q_2,-q_3,-p_1,p_2,p_3)
$$
$$
\rho_2: (q_1,q_2,q_3,p_1,p_2,p_3) \mapsto (q_1,-q_2,q_3,-p_1,p_2,-p_3),
$$
satisfying the relations $\rho_1\circ \rho_2=\rho_2 \circ \rho_1=r$, and so generating the abelian group $\{1,r,\rho_1,\rho_2\}\cong \mathbb{Z}_2 \oplus \mathbb{Z}_2$, which is the natural symmetry group of the spatial circular restricted three-body problem.

After regularization, the symplectic involution admits the following intrinsic description. Consider the smooth reflection $R:S^3 \rightarrow S^3$ along the equatorial sphere $S^2\subset S^3$. Then $r$ is the physical transformation it induces on $T^*S^3$, given by
$$
r: T^*S^3\rightarrow T^*S^3
$$
$$
r(q,p)=(R(q),[(d_qR)^*]^{-1}(p)).
$$
This map preserves the unit cotangent bundle $\mathbb S^*S^3$. The maps $\rho_1,\rho_2$ also have regularized versions. The following emphasizes the symmetries present in our setup:

\begin{proposition}\label{prop:symmetries}
Let $c<H(L_1)$, and consider the symplectic involution $r:\mathbb S^*S^3\rightarrow \mathbb S^*S^3$. The open book decomposition $\Sigma_c=\mathbf{OB}(\mathbb{D}^*S^2,\tau^2)$ is symmetric with respect to $r$, in the sense that
$$
r(P_\theta)=P_{\theta+\pi},\;\;\mbox{Fix}(r)=B=\Sigma_c^P.
$$
Moreover, the anti-symplectic involutions preserve $B$ and satisfy
$$
\rho_1(P_\theta)=P_{-\theta},\;\rho_2(P_\theta)=P_{-\theta + \pi}.
$$
In particular, $\rho_1$ preserves $P_0$ and $P_\pi$, whereas $\rho_2$ preserves $P_{\pi/2}$ and $P_{-\pi/2}$.  
\end{proposition}

In other words, the open book is compatible with all the symmetry group $\mathbb{Z}_2 \oplus \mathbb{Z}_2$.

\subsection*{The return map.} First, we recall a standard definition. We say that a symplectomorphism $f:(M,\omega)\rightarrow (M,\omega)$ is \emph{Hamiltonian} if $f=\phi_K^1$, where $K: \mathbb{R}\times M\rightarrow \mathbb{R}$ is a smooth (time-dependent) Hamiltonian, and $\phi_K^t$ is the Hamiltonian isotopy it generates. This is defined by $\phi_K^0=id$, $\frac{d}{dt}\phi_K^t=X_{K_t}\circ \phi_K^t$, and $X_{H_t}$ is the Hamiltonian vector field of $H_t$ defined via $i_{X_{H_t}}\omega=-dH_t$. Here we write $K_t=K(t,\cdot)$. If $\omega=d\alpha$ is exact, an exact symplectomorphism is a self-diffeomorphism $f$ such that $f^*\alpha=\alpha+d\tau$ for some smooth function $\tau$. All Hamiltonian maps are exact symplectomorphisms. 
 
In our setup, for $c<H(L_1)$, and after fixing a page $P=\pi^{-1}(1)$ of the corresponding open book, Theorem \ref{thm:openbooks} implies the existence of a Poincar\'e return map $f:\mbox{int}(P)\rightarrow \mbox{int}(P)$, defined via $f(p)=\varphi_{\tau(p)}(p),$ where $\varphi_t:\Sigma_c \rightarrow \Sigma_c$ is the Hamiltonian flow, $p\in \mbox{int}(P)$, and $\tau(p)$ is the smallest positive $t$ for which $\varphi_t(p)$ lies in $\mbox{int}(P)$. 

Moreover, we consider the $2$-form $\omega$ obtained by restriction to $P$ of $d\alpha$, where $\alpha$ is the contact form on $\Sigma_c$ for the spatial problem, whose restriction to the binding $\alpha_P=\alpha\vert_B$ is the contact form for the planar problem. Then $\omega$ is symplectic only along the \emph{interior} of $P$.
%(which may be thought of as an ideal Liouville domain, \cite{giroux2}): actually not. 
Moreover, one may find a diffeomorphism $G:\mbox{int}(P)\rightarrow \mbox{int}(\mathbb{D}^*S^2)$ on the interior which extends smoothly to the boundary $B$, but its inverse $G^{-1}$, although continuous at $B$, is \emph{not} differentiable along $B$, in such a way that $\widetilde\omega=G_*\omega$ is now non-degenerate also at $B$ (cf.\ the proof of Lemma \ref{lemma:monodromyvsreturnmap}). After conjugating $f$ with $G$, we obtain a symplectomorphism $\widetilde{f}:=G \circ f \circ G^{-1}:(\mbox{int}(\mathbb{D}^*S^2),\widetilde\omega)\rightarrow (\mbox{int}(\mathbb{D}^*S^2),\widetilde \omega)$, where $\widetilde\omega$ is a Liouville filling of $(B,\alpha_P)$. 

\begin{theorem}\label{thm:returnmap}
For every $\mu \in (0,1]$, $c<H(L_1)$, the associated Poincar\'e return map $f$ extends smoothly to the boundary $B=\partial P$, and in the interior it is an exact symplectomorphism 
$$
f=f_{c,\mu}:(\mbox{int}(P),\omega)\rightarrow (\mbox{int}(P),\omega),
$$
where $\omega=d\alpha$ (depending on $c,\mu$). We have $f(\partial P)=\partial P$, and $f\vert_{\partial P}$ is the time-$1$ map of a positively reparametrized Reeb flow giving the planar three body problem for energy $c$. Moreover, $f$ is Hamiltonian in the interior, generated by a Hamiltonian isotopy which extends smoothly to the boundary. 

After conjugating with $G$, $\widetilde{f}$ extends continuously to the boundary, is Hamiltonian in the interior, generated by a Hamiltonian isotopy which extends continuously to the boundary, and $\widetilde\omega$ has Liouville completion symplectomorphic to the standard symplectic form $\omega_{std}$ on $T^*S^2$. 
\end{theorem}

\begin{figure}
    \centering
    \includegraphics{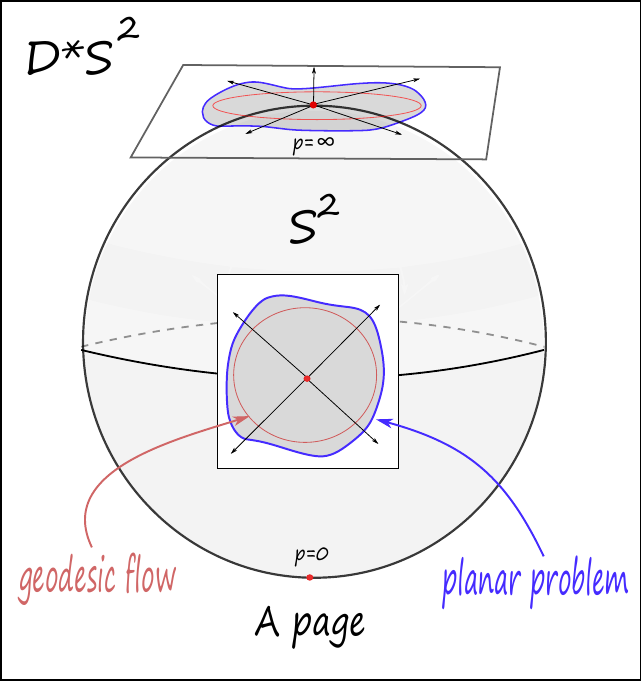}
    \caption{A page of the open book as a filling of the planar problem. After conjugation with a homeomorphism $G$ which makes the page an honest Liouville domain by removing the degeneracy at the boundary, it may be viewed as a fiber-wise star-shaped domain in $T^*S^2$ with the standard symplectic form. The geodesic flow corresponds to the unit cotangent bundle.}
    \label{fig:pageasfilling}
\end{figure}

The form $\widetilde\omega$ can be symplectically deformed, in the class of Liouville fillings of the fixed contact \emph{structure} on $B$, to the standard symplectic form by deforming to the Kepler problem (this can be seen as the limit $c\rightarrow -\infty$, for which $\widetilde f$ converges to the identity). Equivalently, we can think of this Liouville filling as having the standard symplectic form, but non-standard contact boundary (as in Figure~\ref{fig:pageasfilling}).

The fact that $f$ is an exact symplectomorphism follows from the fact that the ambient dynamics comes from a Reeb flow (see Lemma \ref{lemma:symplecto}); this of course implies that $f$ preserves the symplectic volume. The fact that $f$ extends to the boundary is non-trivial, and relies on second order estimates near the binding: it suffices to show that the Hamiltonian giving the spatial problem is positive definite on the symplectic normal bundle to the binding (see Section \ref{sec:secondorderestimates}). This non-degeneracy condition can be interpreted as a convexity condition that plays the role, in this setup, of the notion of \emph{dynamical convexity} due to Hofer-Wysocki-Zehnder (see Definition 3.6 and the proof of Theorem 3.4 in \cite{HWZ98}). Note that if a continuous extension exists, then by continuity it is unique.

The fact that $f$ is Hamiltonian in the interior follows from the following:
\begin{enumerate}
    \item The monodromy of the open book is Hamiltonian (as an isotopy class; here, the Hamiltonian is allowed to move the boundary);
    \item The general fact that the return map $f$ is always symplectically isotopic to a representative of the monodromy, via a boundary-preserving isotopy (Lemma \ref{lemma:monodromyvsreturnmap});
    \item $H^1(P;\mathbb{R})=0$, so that every symplectic isotopy is Hamiltonian.
\end{enumerate}

\begin{remark}[Boundary behaviour] The statement that $f$ is a reparametrized Reeb flow is purely \emph{localized} at the boundary, and follows directly from the construction of a return map; this is not to be confused with the \emph{Hamiltonian twist condition} introduced in \cite{MvK}, which is of a global nature. We remark that positively reparametrized Reeb vector fields are the same as \emph{rotational (or nonsingular) Beltrami fields}, which are relevant in hydrodynamics, and hence both classes are equivalent from a dynamical perspective; see \cite{EG}.
\end{remark}

\begin{remark}[Fixed points] One can extract the following consequence from classical topology: since $f$ is homotopic to the identity, its Lefschetz number is $L(f)=\chi(\mathbb{D}^*S^2)=2$. Therefore it has a least one fixed point $x$, which could be degenerate; if not, there are at least two. Unfortunately, a priori it could lie in the boundary (in \cite{MvK}, to deal with this kind of problem, we will impose suitable convexity assumptions at the boundary). If it does not, we have two possibilities. Indeed, applying the symplectic involution $r$, for which the Hamiltonian vector field is equivariant, we obtain another fixed point $r(x)$ (in the opposite page of the open book, also equivariant by Proposition \ref{prop:symmetries}). If $r(x)$ lies in the same orbit $\gamma$ as $x$, then by uniqueness of solutions this (spatial) orbit is symmetric, i.e.\ $r(\gamma(t))=\gamma(t+1/2)$ for a suitable parametrization. If not, we obtained two spatial orbits. 

\end{remark}

\textbf{Rotating Kepler problem.} In Appendix \ref{app:rotKepler}, we discuss the completely integrable limit case of the rotating Kepler problem, where $\mu=1$ and so there is only one primary. The return map can be studied completely explicitly. Geometrically, this map may be understood via the following proposition (see also Figure \ref{fig:LF}):

\begin{theorem}[Integrable case]\label{prop:integrablecase} In the rotating Kepler problem, the return map $f$ preserves the annuli fibers of a concrete symplectic Lefschetz fibration of abstract type $\mathbb{D}^*S^2=\mathbf{LF}(\mathbb{D}^*S^1,\tau_P^2)$, where it acts as a classical integrable twist map on regular fibers, and fixes the two (unique) nodal singularities on the singular fibers. The boundary of each of the symplectic fibers coincides with the direct/retrograde planar circular orbits (a Hopf link in $\mathbb{R}P^3$).
\end{theorem}

The two fixed points are the north and south poles of the zero section $S^2$, and correspond to the two periodic collision orbits bouncing on the primary (one for each of the half-planes $q_3>0$, $q_3<0$), which we call the \emph{polar} orbits. See Appendix~\ref{app:rotKepler} for an extended discussion, where we derive an explicit formula for the return map, and also describe the Liouville tori. In \cite{M20}, the first author also proves that the structure of a symplectic Lefschetz fibration always exists whenever the planar problem admits adapted open books (which holds e.g.\ when the planar problem is dynamically convex \cite[Theorem 1.18]{HSW}). The symplectic form that makes its fibers symplectic annuli is $d\alpha$, with $\alpha$ the contact form giving the dynamics. We remark that this Lefschetz fibration might not in general be invariant under the return map (but the boundary of its fibers is).

\subsection*{Outlook.} We expect the framework discussed above to provide means for studying more general Hamiltonian systems than the three-body problem. We shall illustrate this expectation by discussing how this works for the more general class of \emph{Stark-Zeeman} systems (see Section \ref{sec:starkzeeman} below), under suitable conditions. An example of such a system describes the dynamics of an electron in an external electric and magnetic field, as well as many other systems in classical and celestial mechanics (see \cite{CFZ,CFvK}). 

The general framework is then the following. Assume that a given Hamiltonian system admits a contact-type and closed energy hypersurface $(\Sigma,\xi=\ker \alpha)=\mathbf{OB}(P,\phi)$, where the Reeb dynamics of the contact form $\alpha$ is the Hamiltonian dynamics, and is adapted to an open book decomposition with data $(P,\phi)$. Here, $P$ is a Liouville domain and $\phi$ is the symplectic monodromy. Then, by considering the return map, one expects to extract dynamical information from the Floer theory of the page $P$. If the return map is Hamiltonian, then one might extract information from Floer theory (e.g.\ symplectic homology). This is the direction pursued in \cite{MvK}.

\subsection*{Proofs}
The main technical ingredients come in the form of various estimates that are scattered over the paper. For the convenience of the reader we have therefore added a small section summarizing the ingredients of the proofs in Section~\ref{sec:summary_proofs}.
We also work out some details concerning the monodromies of the open books in that section.

\subsection*{Acknowledgements.} The authors thank Urs Frauenfelder, for suggesting this problem to the first author, for his generosity with his ideas and for insightful conversations throughout the project; Lei Zhao, Murat Saglam, Alberto Abbondandolo, and Richard Siefring, for further helpful inputs, interest in the project, and discussions. The first author has also significantly benefited from several conversations with Kai Cieliebak in Germany and Sweden; as well as with Alejandro Passeggi in Montevideo, Uruguay. We also thank Gabriele Benedetti, for pointing out a mathematical flaw in an earlier version. This research was started while the first author was affiliated to Augsburg Universität, Germany. The first author is also indebted to a Research Fellowship funded by the Mittag-Leffler Institute in Djursholm, Sweden, where this manuscript was finalized. The second author was supported by NRF grant NRF-2019R1A2C4070302, which was funded by the Korean Government.

\section{The circular restricted three-body problem}

The setup of the classical three-body problem consists of three bodies in $\mathbb{R}^3$, subject to the gravitational interactions between them, which are governed by Newton's laws of motion. We consider three bodies: Earth (E), Moon (M) and Satellite (S), with masses $m_E, m_M, m_S$. We have the following special cases:
\begin{itemize}
    \item(restricted) $m_S=0$ (S is negligible with respect to the \emph{primaries} E and M);
    \item(circular) Each primary moves in a circle, centered around the common center of mass of the two (as opposed to general ellipses); 
    \item(planar) S moves in the plane containing the primaries;
    \item(spatial) The planar assumption is dropped, and S is allowed to move in three-space.
\end{itemize}

The problem then consists in understanding the dynamics of the trajectories of the Satellite, whose motion is affected by the primaries, but not vice-versa. We denote the \emph{mass ratio} by $\mu = \frac{m_M}{m_E+ m _M} \in (0,1],$ and we normalize so that $m_E+ m _M=1$, and so $\mu=m_M$. 

Choose rotating coordinates, in which both primaries are at rest. While the Hamiltonian for the inertial coordinate system is time-dependent, it is autonomous for the rotating ones; the price to pay is the appearance of the angular momentum term. Assuming that the positions of Earth and Moon are $\vec e=(\mu,0,0),$ $\vec m=(-1+\mu,0,0)$\footnote{We will use the symbol $\vec{\phantom{~}}\;$ for vectors in $\R^3$ to make our formulas for Moser regularization simpler. We will use the convention that $\xi\in \R^4$ has the form $(\xi_0,\vec \xi)$.}, the so-called Jacobi Hamiltonian is 
$$
H: \mathbb{R}^3\backslash\{\vec e,\vec m\}\times \mathbb{R}^3 \rightarrow \mathbb{R}
$$
\begin{equation}
\label{eq:Jacobi_Hamiltonian}
H(\vec q,\vec p)=\frac{1}{2}\Vert \vec p\Vert^2 - \frac{\mu}{\Vert \vec q-\vec m\Vert } - \frac{1-\mu}{\Vert \vec q-\vec e\Vert } +p_1q_2-p_2q_1. 
\end{equation}
There are precisely five critical points of $H$, called the \emph{Lagrangian points} $L_i,$ $i=1,\dots,5$, ordered so that $H(L_1)<H(L_2)< H(L_3)<H(L_4)=H(L_5)$ (in the case $\mu <1/2)$. For $c \in \mathbb{R}$, consider the energy hypersurface $\Sigma_c=H^{-1}(c)$. If 
$$
\pi: \mathbb{R}^3\backslash\{\vec e,\vec m\}\times \mathbb{R}^3 \rightarrow  \mathbb{R}^3\backslash\{\vec e,\vec m\},\;\pi(\vec q,\vec p)=\vec q, 
$$
is the projection onto the position coordinate, we define the \emph{Hill's region} of energy $c$ as 
$$
\mathcal{K}_c=\pi(\Sigma_c) \in \mathbb{R}^3\backslash\{\vec e,\vec m\}.
$$
This is the region in space where the satellite with Jacobi energy $c$ is allowed to move. If $c < H(L_1)$ lies below the first critical energy value, then $\mathcal{K}_c$ has three connected components: a bounded one around the Earth, another bounded one around the Moon, and an unbounded one. Denote the first two components by $\mathcal{K}_c^E$ and $\mathcal{K}_c^M$, as well as $\Sigma_c^E=\pi^{-1}(\mathcal{K}_c^E)\cap \Sigma_c$, $\Sigma_c^M=\pi^{-1}(\mathcal{K}_c^M)\cap \Sigma_c$, the bounded components of the corresponding energy hypersurface. As $c$ crosses the first critical energy value, the two connected components $\mathcal{K}_c^E$ and $\mathcal{K}_c^M$ get glued to each other into a new connected component $\mathcal{K}_c^{E,M}$, which topologically is their connected sum. Denote $\Sigma_c^{E,M}=\pi^{-1}(\mathcal{K}_c^{E,M})$.

While the $5$-dimensional energy hypersurfaces are non-compact, due to collisions of the massless body with one of the primaries, two body collisions can be regularized via Moser's recipe (see Section \ref{sec:Moserregularization}). 
This consists in interchanging position and momenta, and compactifying by adding a $2$-sphere at infinity (one point on $S^2$ for each direction), corresponding to collisions where the momentum explodes. 
The bounded components $\Sigma_c^E$ and $\Sigma_c^M$ (for $c < H(L_1)$), as well as $\Sigma_c^{E,M}$ (for $c \in (H(L_1),H(L_1)+\epsilon$), are thus compactified to compact manifolds $\overline{\Sigma}_c^E$, $\overline{\Sigma}_c^M$, and $\overline{\Sigma}_c^{E,M}$. 
The first two are diffeomorphic to $\mathbb S^*S^3=S^3 \times S^2$, and the latter to the connected sum $\mathbb S^*S^3 \# \mathbb S^*S^3$. 

\section{Stark-Zeeman systems}\label{sec:starkzeeman}
Let us assume that $c<H(L_1)$. By restricting the Jacobi Hamiltonian to the Earth or Moon component, we can view it as a Stark-Zeeman system. To define such systems in general, consider a twisted symplectic form
$$
\omega=d\vec p\wedge d\vec q +\pi^*\sigma_B,
$$
with $\sigma_B=\frac{1}{2}\sum B_{ij}dq_i\wedge dq_j$, the magnetic field, where the $B_{ij}$ are smooth functions of $\vec q$. 
A {\em Stark-Zeeman system} is then a Hamiltonian dynamical system for such a symplectic form with a Hamiltonian of the form
$$
H_{SZ}(\vec q, \vec p)=\frac{1}{2} \Vert \vec p \Vert^2+V_0(\vec q)+V_1(\vec q),
$$
where $V_0(\vec q)=-\frac{g}{\Vert \vec q\Vert}$ for some positive coupling constant $g$, and $V_1$ is an extra potential. 

We will make two further assumptions.
\begin{assumptions}
\begin{enumerate}
%$\;$
\item[(A1)] We assume that the magnetic field $\sigma_B$ is exact with primitive $1$-form $\vec A$.
Then with respect to $\omega_0=d\vec p\wedge d\vec q$, we obtain the following Hamiltonian
$$
H(\vec q, \vec p)=\frac{1}{2}\Vert \vec p+\vec A(\vec q) \Vert^2 +V_0(\vec q) +V_1(\vec q),
$$
which has the same dynamics as the above Stark-Zeeman system for the twisted form.
\item[(A2)] We assume that $\vec A(\vec q)=(A_1(q_1,q_2),A_2(q_1,q_2),0)$, and that the potential satisfies the symmetry $V_1(q_1,q_2,-q_3)=V_1(q_1,q_2,q_3)$.
\end{enumerate}
\end{assumptions}
Observe that these assumptions imply that the planar problem, defined as the subset $\{(\vec q,\vec p):q_3=p_3=0\}$, is an invariant set of the Hamiltonian flow. Indeed, we have
\begin{equation}\label{Hamvf}
\dot q_3 = \frac{\partial H}{\partial p_3}= p_3,
\text{ and }
\dot p_3 = -\frac{\partial H}{\partial q_3}
= -\frac{g q_3}{\Vert \vec q\Vert^3} -\frac{\partial V_1}{\partial q_3}.
\end{equation}
Both these terms vanish on the subset $q_3=p_3=0$ by noting that the symmetry implies that $\frac{\partial V_1}{\partial q_3}|_{q_3=0}=0$.

\begin{example}
\label{ex:Kepler}
The Kepler problem is an important example of a Stark-Zeeman system without a magnetic term.
Its dynamics can be described as the Hamiltonian flow of the following Hamiltonian defined on $(T^*\R^3 \setminus \{ 0 \}, d\vec p \wedge d \vec q)$:
$$
K=\frac{1}{2}\Vert \vec p \Vert^2 -\frac{1}{\Vert \vec q \Vert}.
$$
\end{example}

\begin{remark}
\label{rem:SZ_to_std}
The assumption in (A1) allows us to transform a Stark-Zeeman system with Hamiltonian 
$$
H_{SZ}(\vec q, \vec p)=\frac{1}{2} \Vert \vec p \Vert^2+V_0(\vec q)+V_1(\vec q),
$$
for a twisted symplectic form $\omega=d\vec p\wedge d\vec q +\pi^*\sigma_B$ to a Hamiltonian
$$
H(\vec q, \vec p)=\frac{1}{2}\Vert \vec p+\vec A(\vec q) \Vert^2 +V_0(\vec q) +V_1(\vec q),
$$
with the same dynamics for the untwisted form $\omega_0 = d\vec p\wedge d\vec q$.
In the remainder of the paper, we will always transform Stark-Zeeman systems to Hamiltonian systems with untwisted symplectic forms by changing the kinetic term as above.
\end{remark}

\section{Moser regularization}\label{sec:Moserregularization}
For non-vanishing $g$, Stark-Zeeman systems have a singularity corresponding to two-body collisions, which we will regularize by Moser regularization.
To do so, we will define a new Hamiltonian $Q$ on $T^*S^3$ whose dynamics correspond to a reparametrization of the dynamics of $H$ with the symplectic form $\omega_0$ (this Hamiltonian $H$ is the one obtained in Remark~\ref{rem:SZ_to_std}).
We will describe the scheme for energy levels $H=c$ with $c<0$.
Define the intermediate Hamiltonian
$$
K(\vec q,\vec p):=(H(\vec q, \vec p)-c)\Vert\vec q\Vert.
$$
For $\vec q\neq 0$, this function is smooth, and its Hamiltonian vector field equals
$$
X_K= \Vert\vec q\Vert\cdot X_H +(H-c) X_{\Vert\vec q\Vert}. 
$$
We observe that $X_K$ is a multiple of $X_H$ on the level set $K=0$.
Writing out $K$ gives
$$
K=
\left( \frac{1}{2}(\Vert \vec p\Vert^2+1)-(c+1/2) +\langle\vec p,\vec A \rangle +\frac{1}{2}\Vert\vec A\Vert^2+V_1(\vec q)
\right) \Vert\vec q\Vert
-g
. 
$$
\begin{remark}
At this point, it is worth pointing out that the round metric for a sphere $S^n_R$ with radius $R$ has the following form in stereographic coordinates
$$
g_x(v,v)= \frac{4R^4}{(|x|^2+R^2)^2}|v|^2.
$$
The above expression is hence a deformation of the norm of $\vec q$, interpreted as a covector; the point $\vec p$ plays the role of the coordinate on the base.
We will explain these coordinates below.
\end{remark}
\textbf{Stereographic projection.} We set 
$$
\vec x=- \vec p, \;\;\vec y=\vec q.
$$ We view $T^*S^3$ as a symplectic submanifold of $T^*\R^{4}$, via
$$
T^*S^3
=
\{
(\xi,\eta)\in T^*\mathbb{R}^{4} \vert\; \Vert \xi\Vert^2=1,\;
\langle \xi, \eta \rangle =0
\}.
$$
Let $N=(1,0,0,0)\in S^3$ be the north pole. To go from $T^*S^3\backslash T^*_N S^3$ to $T^*\R^3$ we use the stereographic projection. Recalling the notation $\xi=(\xi_0,\vec \xi)\in \mathbb{R}\times \mathbb R^4$, $\eta=(\eta_0,\vec \eta)\in \mathbb{R}\times \mathbb R^4$, this is given by
\begin{equation}
\label{eq:regularized_to_unregularized}
\begin{split}
\vec x &= \frac{\vec \xi}{1-\xi_0} \\
\vec y &= \eta_0 \vec \xi +(1-\xi_0) \vec \eta .
\end{split}
\end{equation}

To go from $T^*\R^3$ to $T^*S^3\backslash T^*_NS^3$, we use the inverse given by
\begin{equation}
\label{eq:unregularized_to_regularized}
\begin{split}
\xi_0 &= \frac{\Vert \vec x\Vert ^2-1}{\Vert \vec x\Vert ^2+1} \\
\vec \xi &= \frac{2 \vec x}{\Vert \vec x\Vert ^2+1} \\
\eta_0 &= \langle \vec x,\vec y \rangle \\
\vec \eta &= \frac{\Vert \vec x\Vert ^2+1}{2}\vec y 
- \langle \vec x, \vec y \rangle \vec x.
\end{split}
\end{equation}

These formulas imply the following identities
$$
\frac{2}{\Vert \vec x\Vert ^2+1}= 1-\xi_0, \;\;
\Vert \vec y\Vert =\frac{2\Vert\eta \Vert}{\Vert \vec x\Vert^2+1}
=(1-\xi_0) \Vert \eta\Vert, 
$$
which allows us to simplify the expression for $K$. We obtain a Hamiltonian $\tilde K$ defined on $T^*S^3$, given by
\[
\begin{split}
\tilde K
&= 
\left( 
\frac{1}{1-\xi_0}
-(c+1/2) 
-\frac{1}{1-\xi_0} \langle \vec \xi , \vec A(\xi,\eta) \rangle
+
\frac{1}{2}\Vert \vec A(\xi,\eta)\Vert ^2
+
V_1(\xi,\eta)
\right)
(1-\xi_0) \Vert \eta \Vert -g 
\\
& =\Vert \eta \Vert
\left(
1-(1-\xi_0)(c+1/2)
-\langle \vec \xi, \vec A(\xi,\eta) \rangle
+
(1-\xi_0)\left(\frac{1}{2}\Vert \vec A(\xi,\eta)\Vert ^2
+
V_1(\xi,\eta)
\right)
\right)
-g.
\end{split}
\]
Put
\begin{equation}
\label{eq:Stark_Zeeman_f}
\begin{split}
f(\xi,\eta)
&=
1+(1-\xi_0)\left(-(c+1/2)
+
\frac{1}{2}\Vert \vec A(\xi,\eta)\Vert ^2
+
V_1(\xi,\eta)
\right)
-\langle \vec \xi, \vec A(\xi,\eta) \rangle
\\
&=1+(1-\xi_0) b(\xi,\eta) +M(\xi,\eta),
\end{split}
\end{equation}
where 
$$
b(\xi,\eta)=-(c+1/2)
+
\frac{1}{2}\Vert \vec A(\xi,\eta)\Vert ^2
+
V_1(\xi,\eta)
$$ 
$$
M(\xi,\eta)=-\langle \vec \xi, \vec A(\xi,\eta) \rangle.
$$
Note that the collision locus corresponds to $\xi_0=1$, i.e.\ the cotangent fiber over $N$. We then have that 
$$
\tilde K=\Vert \eta \Vert f(\xi,\eta)-g.$$ To obtain a smooth Hamiltonian, we define the Hamiltonian
$$
Q(\xi,\eta):= \frac{1}{2} f(\xi,\eta)^2 \Vert \eta \Vert^2.
$$
The dynamics on the level set $Q=\frac{1}{2} g^2$ are a reparametrization of the dynamics of $\tilde K=0$, which in turn correspond to the dynamics of $H=c$.
\begin{remark}
We have chosen this form to stress that $Q$ is a deformation of the Hamiltonian $Q_{round} =\frac{1}{2} \Vert \eta \Vert^2$ describing the geodesic flow on the round sphere. This is the the regularized Kepler problem, corresponding to the Reeb dynamics of the standard contact form, a Giroux form (in the sense of definition \ref{def:Giroux} below) for the open book $\mathbb S^*S^3=\mathbf{OB}(\mathbb{D}^*S^2,\tau^2)$, supporting the standard contact structure on $\mathbb S^*S^3$. 
\end{remark}

\subsection{Formula for the restricted three-body problem.}\label{sec:form3BP} By completing the squares, we obtain
$$
H(\vec q, \vec p)=\frac{1}{2}\left((p_1+q_2)^2+(p_2-q_1)^2+p_3^2\right)-\frac{\mu}{\Vert \vec q - \vec m\Vert}-\frac{1-\mu}{\Vert \vec q -\vec e\Vert} -\frac{1}{2}(q_1^2+q_2^2).
$$
Keeping Remark~\ref{rem:SZ_to_std} in mind, we can view the Hamiltonian dynamics of this Hamiltonian for the untwisted symplectic form $\omega_0$ as a Stark-Zeeman system with primitive 
$$\vec A=(q_2,-q_1,0)$$
and potential 
\begin{equation}\label{potential}
V_1(\vec q)=-\frac{1-\mu}{\Vert \vec q -\vec e\Vert} -\frac{1}{2}(q_1^2+q_2^2),
\end{equation}
both of which satisfy Assumptions (A1) and (A2).

After a computation, we obtain
\begin{equation}\label{eq:f3bp}
f(\xi,\eta)=
1+ \left( 1-\xi_{{0}} \right)  \left( -(c+1/2)+ \xi_{{2}}\eta_{{1}}-\xi_{{1}}
\eta_{{2}} \right) -\xi_{{2}} \left( 1-\mu \right)-\frac {(1-\mu) (1-\xi_{{0}}) }{\Vert \vec \eta (1-\xi_0)+\vec \xi \eta_0 + \vec m - \vec e \Vert}
\end{equation}
\begin{equation}\label{eq:b}
b(\xi,\eta)= -(c+1/2) -\frac {(1-\mu) }{\Vert \vec \eta (1-\xi_0)+\vec \xi \eta_0 + \vec m - \vec e \Vert}
\end{equation}
\begin{equation}\label{eq:M1}
M(\xi,\eta)=(1-\xi_0) (\xi_{{2}}\eta_{{1}}-\xi_{{1}}
\eta_{{2}}) - \xi_2(1-\mu).
\end{equation}

\subsection{Hamiltonian vector field} Consider the Hamiltonian
$$
Q^0=\frac{1}{2} \Vert \eta \Vert^2 f^2(\xi,\eta)
$$
on $(T^*\R^n,\omega_0=d\eta \wedge d\xi)$, for some smooth function $f$, and let $Q:=Q^0|_{T^*S^{n-1}}$. A computation using the submanifold constraints $\Vert \xi \Vert^2=1$ and $\langle \xi,\eta \rangle=0$ gives:

\begin{lemma}
We have
\begin{equation}\label{HamQvf}
\begin{split}
X_Q=&f\left(f \eta +\Vert \eta \Vert^2 \left(f_\eta-\xi f_\eta \cdot \xi\right)\right)\partial_\xi\\
&+\Vert \eta \Vert^2 f \left(\eta f_\eta\cdot \xi-f_\xi-\xi \left(f+ 
f_\eta \cdot \eta
-
f_\xi \cdot \xi
\right)\right)\partial_\eta. 
\end{split}
\end{equation}

\end{lemma}

\begin{example}
In the case of the round sphere ($f\equiv 1$), the above reduces to
\begin{equation}\label{roundsphere}
X_Q=\eta \cdot \partial_\xi - \Vert \eta \Vert^2\xi \cdot \partial_\eta.
\end{equation}
\end{example}

\section{Contact topology and dynamics}

In this section, we recollect important notions, which were also discussed in the Introduction.

\subsection{Open book decompositions and global hypersurfaces of section}\label{OBDS} A {\em concrete open book} $(B,\pi)$ on a manifold $Y$ consists of
\begin{enumerate}
\item a codimension $2$ closed submanifold $B$ with trivial normal bundle, and
\item a fiber bundle $\pi: Y\setminus B\to S^1$ such that for some collar neighbourhood $B\times \mathbb{D}^2$ of $B$ we have
$$
\pi: B\times \mathbb{D}^2 \longrightarrow S^1,\quad
(b;r,\theta)  \longmapsto \theta,
$$
where $(r,\theta)$ are polar coordinates on $\mathbb{D}^2$.
\end{enumerate}
The submanifold $B$ will be called the {\em binding}. The closure of the fibers of $\pi$ are called {\em pages}. Abstractly, given a manifold $P$ with boundary $B$, and a diffeomorphism $\phi: P \rightarrow P$ with $\phi\vert_B=id$, we construct an \emph{abstract} open book 
$$
\mathbf{OB}(P,\phi)=B\times \mathbb{D}^2\bigcup_\partial \mbox{Map}(P,\phi),
$$
where $\bigcup_\partial$ denotes boundary union, and $\mbox{Map}(P,\phi)=P \times [0,1]/(x,0)\sim (\phi(x),1)$ is the associated mapping torus. 
An abstract open book induces an obvious concrete open book, and vice-versa (uniquely up to isotopy). 

\medskip

The above is so far a notion of smooth topology, and its relationship to contact topology is encoded in the following definition. Recall that a (positive) \emph{contact form} on an (oriented) odd-dimensional manifold $Y^{2n+1}$ is a $1$-form $\alpha \in \Omega^1(Y)$ satisfying the contact condition $\alpha \wedge d\alpha^n>0$. The induced \emph{contact structure} is the hyperplane distribution $\xi=\ker \alpha \subset TY$, and $(Y,\xi)$ is a \emph{contact manifold}.
We have the following definition due to Giroux, \cite{Gir}.

\begin{definition}[Giroux]\label{def:Giroux} Suppose that $Y$ is equipped with a concrete open book $(B,\pi)$ and a contact form $\alpha$ satisfying the following:
  \begin{itemize}
  \item{} $\alpha$ induces a positive contact structure on $B$, and
  \item{} $d\alpha$ induces a positive symplectic structure on the interior of each page of $\pi$.
  \end{itemize}
Then we say that $\alpha$ is \emph{adapted} to $(B, \pi)$, or that $\alpha$ is a \emph{Giroux form} for the open book, and that $\xi=\ker \alpha$ is \emph{supported} by the open book. 
\end{definition}

 In the above situation, $(B,\xi_B=(\ker \alpha\vert_B))$ is a \emph{contact submanifold} of $Y$ (i.e.\ $\xi \vert_{TB}=\xi_B$). We usually write $(Y,\xi)=\mathbf{OB}(P,\phi)$ whenever $\xi$ is supported by the abstract open book with data $(P,\phi)$.

\begin{example}
We have $(\mathbb S^*S^3,\xi_{std})=\mathbf{OB}(\mathbb{D}^*S^2, \tau^2)$, where $\tau$ is the Dehn-Seidel twist along the Lagrangian zero section $S^2\subset \mathbb{D}^*S^2$. See Section \ref{sec:geodesic_OB}.
\end{example}

Any contact form $\alpha$ on $Y$ gives rise to a dynamical system, given by the flow of the \emph{Reeb vector field} $R_\alpha$, which is defined implicitly via the equations
$$
\alpha(R_\alpha)=1,\;\;d\alpha(R_\alpha,\cdot)=0
$$

There is a relation between open books and global hypersurfaces of section for the Reeb flow.
One part of this relation is clearly expressed in the following lemma.
\begin{lemma}
\label{lemma:Reeb}
Suppose that $B$ is a connected contact submanifold of a contact manifold $(Y,\xi)$.  A contact form $\alpha$ for $(Y,\xi)$ is adapted to an open book $(B,\pi)$ if and only if 
\begin{itemize}
    \item $B$ is invariant under the flow of $R_\alpha$;
    \item $R_\alpha$ is positively transverse to the fibers of $\pi$,
  i.e.~$d\pi(R_\alpha) > 0$.
\end{itemize}
Assume now furthermore that we have a bound on the return time.
Then it follows that every page is a global hypersurface of section for the Reeb dynamics.
If the contact condition does not hold, transversality to all pages is of course still enough for a single page to be a global hypersurface of section.
\end{lemma}
We note that, in the situation of the above lemma, we have the following general fact: 
\begin{lemma}\label{lemma:symplecto}
The associated return map $f$ on each page is automatically an exact symplectomorphism with respect to the symplectic form induced by the restriction of $d\alpha$.
\end{lemma}
\begin{proof}
Let $\omega=d\alpha\vert_P=d\lambda$ where $P$ is a fixed page and $\lambda=\alpha\vert_P$, denote the time-$t$ Reeb flow by $\varphi_t$, and let $\tau:\mbox{int}(P)\rightarrow \mathbb{R}^+$, $\tau(x)=\min\{t>0:\varphi_t(x)\in \mbox{int}(P)\}.$ Then $f(x)=\varphi_{\tau(x)}(x)$, and so, for $x \in \mbox{int}(P)$, $v \in T_xP$, we have
$$
d_xf(v)=d_x\tau(v)R_{\alpha}(f(x))+d_x\varphi_{\tau(x)}(v).
$$
Using that $\varphi_t$ satisfies $\varphi_t^*\alpha=\alpha$, we obtain
\begin{equation}
    \begin{split}
        (f^*\lambda)_x(v)&=\alpha_{f(x)}(d_xf(v))\\
        &=d_x\tau(v)+(\varphi_{\tau(x)}^*\alpha)_x(v)\\
        &=d_x\tau(v)+\lambda_x(v).\\
    \end{split}
\end{equation}
Therefore $f^*\lambda=d\tau + \lambda$, which shows the claim.
\end{proof}

It is not always the case that the return map extends continuously to the boundary, see for instance Remark 7.1.11 in \cite{FvK18} for a return map that twists ``infinitely fast'' along the boundary.
This is a rather subtle point, to which we will come back to in Section \ref{sec:secondorderestimates} in order to obtain Theorem \ref{thm:returnmap}; see Proposition~\ref{proposition:convexity}.

\begin{comment}
One may also consider hypersurfaces of section in situations where the contact property does not necessarily hold.
\begin{lemma}[First order estimates: general case]\label{firstordergeneral}
Suppose that $Y^n$ is a compact smooth manifold equipped with a non-vanishing autonomous vector field $X$.
Assume that 
\begin{itemize}
    \item $B^{n-2}$ is a closed $n-2$-dimensional, smooth submanifold of $Y$ that is invariant under the flow of $X$.
    \item  $\pi: Y\setminus B \to S^1$ is a smooth fiber bundle such that $d \pi(X) > 0$. 
\end{itemize}
Then each fiber of $\pi$ is a global hypersurface of section for the flow of $X$.
\end{lemma}

In the above situation, there is no a priori reason why the arising return map preserves some volume form, so the most interesting case is indeed the contact-type one.
\end{comment}

\section{First order estimates}\label{sec:firstorderestimates}

In this section, we will setup the geometric situation, and conclude the proof of the first part of Theorem~\ref{thm:gss_openbooks}. 
We will assume nothing on the energy $c$, and prove directly that the Hamiltonian vector field $X_Q$ is transverse to the planar problem, for Stark-Zeeman systems satisfying Assumptions (A1)--(A4), yet to be fully determined. 
This gives global hypersurfaces of section even if the contact-type condition fails. 
In the case where $c<H(L_1)$ or $c \in (H(L_1),H(L_1)+\epsilon),$ so that $\Sigma_c$ is contact-type, we will obtain our result in Theorem~\ref{thm:openbooks} via Lemma \ref{lemma:Reeb}. 

\subsection{The physical open book.} For Stark-Zeeman systems satisfying Assumptions (A1) and (A2) we will define a natural candidate for an open book. As noted before, the planar problem defines an invariant subset. In unregularized coordinates, we put 
$$
B_u:=\{(\vec q,\vec p)\in H^{-1}(c)~|~q_3=p_3=0 \}.
$$
Its normal bundle is trivial, and we have the following map to $S^1$, 
\begin{equation}\label{pinonreg}
\pi_p: H^{-1}(c) \setminus B_u \longrightarrow S^1,
(\vec q,\vec p) \longmapsto \frac{q_3+ip_3}{\Vert q_3+ip_3\Vert}.
\end{equation}
We will refer to this map as the \emph{physical} open book. We consider the angular $1$-form
$$
\omega_p:=\frac{\Omega^u_p}{p_3^2+q_3^2},
$$
where
\begin{equation}\label{omegapunreg}
\Omega^u_p=p_3 dq_3 -q_3 dp_3,
\end{equation}
is the unregularized numerator. In view of Lemma~\ref{lemma:Reeb}, we need to see whether $\omega_p(X_H)$ is non-negative, and vanishes only along the planar problem.

From Equation~\eqref{Hamvf}, we have
\begin{equation}\label{omegaXH}
\omega_p(X_H)=\frac{p_3^2+q_3^2\left( \frac{g}{\Vert \vec q\Vert^3}+\frac{1}{q_3} \frac{\partial V_1}{\partial q_3}(\vec q)\right)}{p_3^2+q_3^2}.
\end{equation}

Assumption (A2) implies that $\frac{\partial V_1}{\partial q_3}(\vec q)=a q_3+o(q_3^2)$ near $q_3=0$, and so $\frac{1}{q_3} \frac{\partial V_1}{\partial q_3}(\vec q)$ is well-defined at $q_3=0$. In order for the above expression to satisfy the required non-negativity condition, we impose the following:

\begin{assumption}

(A3) We assume that the function 
    $$
    F(\vec q)=\frac{g}{\Vert \vec q\Vert^3}+\frac{1}{q_3} \frac{\partial V_1}{\partial q_3}(\vec q)
    $$ is everywhere positive.
\end{assumption}
Note that it suffices that the second summand be non-negative.

\begin{remark}\label{rem:A3}
In the restricted three-body problem, from Equation~\eqref{potential}, we obtain
$$
\frac{\partial V_1}{\partial q_3}(\vec q)=q_3\frac{1-\mu}{\Vert \vec q - \vec e\Vert^3},
$$
and therefore Assumption (A3) is satisfied.
\end{remark}

This observation only applies to the unregularized problem, and we will want to look at the \emph{compact} hypersurfaces of section. We hence need the following expression for the above map $\pi_p$ in $(\xi,\eta)$ coordinates. Let $S: T^* (S^3 \setminus \{ N\}) \longrightarrow T^*\R^3$ be the stereographic projection map.
\begin{lemma}\label{thetap}
With 
$$
\Theta_p(\xi,\eta) = \xi_3+i\left( (1-\xi_0)\eta_0 \xi_3 +(1-\xi_0)^2 \eta_3 \right), 
$$
we have
$$
i \pi_p\circ S(\xi,\eta) =\frac{\Theta_p(\xi,\eta)}{|\Theta_p(\xi,\eta)|}.
$$
\end{lemma}

\begin{proof}
We consider the denominator of $i \pi_p(q,p)$ and use the formulas for Moser regularization
$$
-p_3+iq_3
=
\frac{\xi_3}{1-\xi_0}+i\left( \eta_0 \xi_3 +(1-\xi_0)\eta_3 \right).
$$
We rescale by $(1-\xi_0)$, which doesn't change the map to $S^1$, to obtain the claim.
\end{proof}

The associated angular 1-form in $(\xi,\eta)$-coordinates is
$$
\omega_p = 
\frac{\Omega_p}
{\re (\Theta_p(\xi,\eta))^2 +\im (\Theta_p(\xi,\eta))^2}
$$
where
\begin{equation}\label{Omegap}
\Omega_p
=\re(\Theta_p(\xi,\eta))d \im(\Theta_p(\xi,\eta))
-
\im(\Theta_p(\xi,\eta))
d\re(\Theta_p(\xi,\eta))
\end{equation}
is the regularized numerator. By construction, we have the following relationship between the unregularized and regularized numerators:
\begin{equation}\label{regvsunreg}
  \Omega_p=(1-\xi_0)^2\Omega_p^u. 
\end{equation}

\subsection{The geodesic open book: a simple example.}
\label{sec:geodesic_OB}
Before doing the general case, let us examine a higher-dimensional analogue of the famous Birkhoff open book in the simple case of the round sphere. The Hamiltonian is $Q= \frac{1}{2} \Vert \eta \Vert^2 |_{T^*S^{n}}$ with Hamiltonian vector field
$$
X_Q=\eta \cdot \partial_{\xi} -\xi \cdot \partial_{\eta}.
$$
This is the Reeb vector field of the standard Liouville form $\lambda_{std}$ on the energy hypersurface $\Sigma=Q^{-1}(\frac{1}{2})=\mathbb S^*S^n$. We have the invariant set
$$
B:=\{ (\xi_0,\ldots,\xi_n;\eta_0,\ldots, \eta_n) \in \Sigma~|~\xi_n=\eta_n=0\}=\mathbb S^*S^{n-1}.
$$
Define the circle-valued map
$$
\pi_g: \Sigma \setminus B \longrightarrow S^1,
\quad
(\xi_0,\ldots,\xi_n;\eta_0,\ldots, \eta_n)
\longmapsto \frac{\eta_n +i \xi_n}{\Vert \eta_n+i\xi_n\Vert}
.
$$
This is the projection of the concrete open book which was discussed in the introduction, which we shall refer to as the \emph{geodesic} open book; note that the page $\xi_n=0$ and $\eta_n>0$ corresponds to the higher-dimensional Birkhoff ``annulus'' $\mathbb{D}^*S^{n-1}$. The angular form is then
$$
\omega_g =\frac{\eta_n d\xi_n -\xi_n d \eta_n}{\xi_n^2 +\eta_n^2}.
$$
We see that $\omega_g(X_Q)=1>0$, so Lemma~\ref{lemma:Reeb} tells us that $(B,\pi_g)$ is a supporting open book for $\Sigma$ and the pages of $\pi_g$ are global hypersurfaces of section for $X_Q$. Abstractly, this gives $(\mathbb S^*S^{n},\xi_{std})=\mathbf{OB}(\mathbb{D}^*S^{n-1},\tau^2)$.

\subsection{First order estimates for spatial Stark-Zeeman systems} We now do the general case. We consider a connected component of the energy hypersurface of a regularized Stark-Zeeman system, which we denote by $\Sigma \subset Q^{-1}(g^2/2)$.
From Formula~\eqref{eq:Stark_Zeeman_f}, we know that $f$ is of the form
$$
f=1+(1-\xi_0) b(\xi,\eta)+M(\xi,\eta),
$$
for some smooth functions $b,M$. We want to consider the analogue of the geodesic open book we considered in Section \ref{sec:geodesic_OB} (for $n=3)$. With the same formulas for $\pi_g$ we have again the angular form
$$
\omega_g = \frac{\Omega_g}{\xi_3^2+\eta_3^2},
$$
where
$$
\Omega_g := \eta_3 d\xi_3-\xi_3 d\eta_3.
$$
Note that $\pi_g$ does \emph{not} agree with $\pi_p$ in regularized coordinates (see Lemma~\ref{thetap}). 

From Equation~\eqref{HamQvf}, it follows from direct computation that:

\begin{lemma}\label{OmegaXQ}
We have the expression
$$
\Omega_g(X_Q)=f^2\eta_3^2+\Vert \eta \Vert^2 f(f+f_\eta \cdot \eta - f_\xi \cdot \xi)\xi_3^2+\Vert \eta \Vert^2 f(f_{\eta_3}\eta_3-f_{\xi_3}\xi_3-2\xi_3\eta_3f_\eta \cdot \xi)
$$
\medskip
If we write $f=1+(1-\xi_0)b+M,$ then this further gives
\begin{equation}\label{OmegaXQwithb}
\begin{split}
\Omega_g(X_Q)=&f^2\eta_3^2+\Vert \eta \Vert^2f(1+\xi_0 b+M+M_\eta \cdot \eta-M_\xi\cdot \xi+(1-\xi_0)(b+b_\eta \cdot \eta - b_\xi\cdot \xi))\xi_3^2\\
&+\Vert \eta \Vert^2f((1-\xi_0)(\eta_3b_{\eta_3}+\xi_3b_{\xi_3}-2\xi_3\eta_3b_\eta\cdot \xi)+\eta_3 M_{\eta_3}+\xi_3M_{\xi_3}-2\eta_3\xi_3 M_\eta\cdot \xi)
\end{split}
\end{equation}
\end{lemma}
$\square$

\medskip

Note that setting $f=1=\Vert \eta \Vert$ in the above expression recovers the example from Section \ref{sec:geodesic_OB}. We will not verify on the whole set $\Sigma$ whether $\Omega_g(X_Q)$ is non-negative, but instead only whether this holds near the collision locus, and combine this with Expression~\eqref{omegaXH} from our earlier computation in unregularized coordinates. The basic observation is that $\Omega_p$ works away from the collision locus, whereas $\Omega_g$ works near the collision locus. 
We therefore interpolate between the two. This creates an interpolation region where we need finer estimates in order to obtain global hypersurfaces of section. This is the content of what follows.

Assume that $Q=\frac{1}{2} \Vert \eta \Vert^2 f(\xi,\eta)^2$ is the regularized Hamiltonian of Stark-Zeeman system satisfying (A1), (A2) and (A3), and write $f=1+(1-\xi_0) b+M$. Let $\Sigma$ be a connected component of the regularized energy hypersurface $Q^{-1}(g^2/2)$, which we assume to be closed. Further assume that:
\begin{assumption}
\begin{enumerate}
    \item[(A4)] $1+b(\xi,\eta)+M(\xi,\eta) -M_\xi(\xi,\eta) \cdot \xi>0$ for all $(\xi,\eta)\in \Sigma$ with $\xi=(1,0,0,0)$.
\end{enumerate}
\end{assumption}

\begin{lemma}
\label{lemma:GSS_Stark_Zeeman}
Under Assumptions (A1)-(A4) as above, there exists an open book decomposition on $\Sigma$, with binding the planar problem, so that each page is a global hypersurface of section for the Hamiltonian flow $X_Q$.
\end{lemma}

\begin{proof}
Recall that the collision locus corresponds to $\xi_0=1$, is diffeomorphic to $S^2$ and its points have the form $(\xi,\eta)=(1,0,0,0;\eta)$. Define $B:= \{ (\xi,\eta)\in \Sigma~|~\xi_3=\eta_3=0 \}$.
By Assumption (A2) and Moser regularization we see that $B$ is an invariant set for the flow of $X_Q$ (it is the regularization of $B_u$). Choose a smooth, non-negative function $\rho=\rho(\xi_0)$, which is positive near $\xi_0=1$, and define the map
\begin{equation}\label{thetamap}
\Theta:
\Sigma \longrightarrow \C,
\quad
(\xi,\eta) \longmapsto
\Theta_p(\xi,\eta)+i \rho(\xi_0) \eta_3,
\end{equation}
where $\Theta_p$ is as in Lemma \ref{thetap}. We note that $\Theta(\xi,\eta)=0$ if and only if $(\xi,\eta) \in B$, since, by Lemma \ref{thetap}, $\re (\Theta)=\xi_3$ and $\im (\Theta)=(1-\xi_0) \eta_0 \xi_3+( (1-\xi_0)^2+\rho(\xi_0))\eta_3$.
Hence we obtain the circled-valued map
$$
\pi: \Sigma \setminus B \longrightarrow S^1,\quad
(\xi,\eta) \longmapsto \Theta(\xi,\eta)/\Vert\Theta(\xi,\eta)\Vert.
$$
Note that $\pi=\pi_p$ away from the collision locus (where $\rho=0$), and $\pi=\pi_g$ at the collision locus.
The associated angular $1$-form is
$$
\omega = \frac{ \Omega }
{\re(\Theta(\xi,\eta))^2 +\im(\Theta(\xi,\eta))^2},
$$
where
$$
\Omega = \re(\Theta(\xi,\eta))d\im(\Theta(\xi,\eta))
-
\im(\Theta(\xi,\eta))
d\re(\Theta(\xi,\eta))
=
\Omega_p +\rho(\xi_0) \Omega_g +\xi_3 \eta_3 d\rho.
$$
We will apply Lemma~\ref{lemma:Reeb} to verify the open book condition, so we need to check whether $\Omega(X_Q)>0$. To achieve this, we will impose further conditions on $\rho$ as we go.

\medskip

{\bf Claim 1: } \emph{For any small $\epsilon>0$, there exists $C_\epsilon >0$ such that if $\xi_0\leq 1-\epsilon$, then $\Omega_p(X_Q)\geq C_\epsilon (\xi_3^2+ \eta_3^2)$. Furthermore, we have $\Omega_p(X_Q)\geq 0$ everywhere, with equality only along $B$.}

\medskip

To see that this inequality holds, we recall Equation~\eqref{omegaXH}.
This equation tells us, under Assumption~(A3), that $\Omega_p(X_H)_{(q,p)}\geq 0$ with equality only along $B$. Since $X_Q=hX_H$ for some positive function $h$ away from the collision locus, we also have $\Omega_p(X_Q)_{(\xi,\eta)}\geq 0$ for all $(\xi,\eta)\in \Sigma$ with $\xi_0\neq 1$, with equality only along $B\cap \{\xi_0\neq 1\}$. To see that $\Omega_p(X_Q)\geq 0$ on the collision locus (and hence everywhere), note from Equation~\eqref{regvsunreg} that $\Omega_p$ is a smooth $1$-form that vanishes at $\xi_0=1$. We conclude that $\Omega_p(X_Q)\geq 0$ for all $(\xi,\eta) \in \Sigma $, with equality only along $B$. To obtain the lower quadratic bound for $\xi_0<1$, we use Equations \eqref{omegaXH} and \eqref{regvsunreg} to get
$$
\Omega_p(X_Q)=h\cdot\Omega_p(X_H)=h\cdot (1-\xi_0)^2\Omega_p^u(X_H)=h\cdot(1-\xi_0)^2(p_3^2+q_3^2F(\vec q)).
$$
Given $\epsilon>0$ small, compactness of $\Sigma$ implies uniform lower bounds $h\geq K_\epsilon$ and $F\geq K_\epsilon^\prime$ along $\Sigma$ (using (A3)), and so we may choose $C^\prime_\epsilon=K_\epsilon \min\{1,K_\epsilon^\prime\}$ so that
$$
\Omega_p(X_Q)\geq C_\epsilon^\prime (1-\xi_0)^2(p_3^2+q_3^2) \geq  C_\epsilon^\prime \epsilon^2 (p_3^2+q_3^2).
$$
From Moser regularization, we get 
$$
p_3^2 = \frac{\xi_3^2}{(1-\xi_0)^2},\quad
q_3^2 = \eta_0^2 \xi_3^2+(1-\xi_0)^2 \eta_3^2 +2 \eta_0 (1-\xi_0) \xi_3 \eta_3,
$$
so 
\[
\begin{split}
p_3^2+q_3^2 &=  \frac{\xi_3^2}{(1-\xi_0)^2}+\eta_0^2 \xi_3^2+(1-\xi_0)^2 \eta_3^2 +2 \eta_0 (1-\xi_0) \xi_3 \eta_3 \\
&=
\left(
\begin{array}{c}
\xi_3 \\
\eta_3
\end{array}
\right)^t
\left(\begin{array}{cc}
\frac{1}{(1-\xi_0)^2}+\eta_0^2 & \eta_0 (1-\xi_0) \\
\eta_0 (1-\xi_0) & (1-\xi_0)^2
\end{array}
\right)
\left(
\begin{array}{c}
\xi_3 \\
\eta_3
\end{array}
\right).
\end{split}
\]

In order to think of the matrix as a metric, we need to verify that the eigenvalues are positive.
For this we compute the determinant and trace of the associated quadratic form in $\xi_3,\eta_3$.
The trace is $\frac{1}{(1-\xi_0)^2}+\eta_0^2+(1-\xi_0)^2>0$ if $\xi_0\neq 1$.
The determinant is $1$, so near $\xi_0=1$ (but not at $\xi_0=1$) this matrix represents a metric.
In particular, we can bound $p_3^2+q_3^2$ from below by $c_\epsilon(\xi_3^2+\eta_3^2)$ for some constant $c_\epsilon$. We can then set $C_\epsilon=c_\epsilon C_\epsilon^\prime \epsilon^2$ to get the claim.

\medskip

{\bf Claim 2: } \emph{We can find $\delta>0$ and $A_\delta >0$ such that, along $\Sigma$, we have $\Omega_g(X_Q) \geq A_\delta (\xi_3^2+ \eta_3^2)$ for $\xi_0\geq 1-\delta$. In particular, $\Omega_g(X_Q)\geq 0$ for $\xi_0\geq 1-\delta$, with equality only along $B \cap \{\xi_0\geq 1-\delta\}$.}

\medskip

We will verify this claim with a computation. From Equation~\eqref{OmegaXQwithb} we have
\[
\begin{split}
\Omega_g(X_Q)=&f^2\eta_3^2+\Vert \eta \Vert^2f(1+\xi_0 b+M+M_\eta \cdot \eta-M_\xi\cdot \xi+(1-\xi_0)(b+b_\eta \cdot \eta - b_\xi\cdot \xi))\xi_3^2\\
&+\Vert \eta \Vert^2f((1-\xi_0)(\eta_3b_{\eta_3}+\xi_3b_{\xi_3}-2\xi_3\eta_3b_\eta\cdot \xi)+\eta_3 M_{\eta_3}+\xi_3M_{\xi_3}-2\eta_3\xi_3 M_\eta\cdot \xi).
\end{split}
\]
The first term, which we will abbreviate by $T_1$, is obviously non-negative, and the coefficient $f^2$ in front of $\eta_3^2$ is strictly positive on $\Sigma$ (since $g>0$, recalling that $Q=\frac{1}{2}\Vert \eta\Vert ^2f^2=\frac{1}{2}g^2$ along $\Sigma$). Let us now deal with the second term $T_2$. We first consider the limit case $\xi_0=1$, which means that $\xi_1=\xi_2=\xi_3=0$. Hence several terms drop out, and we will use Assumption (A4) to further simplify this expression.

Recall that $M(\xi,\eta)=-\langle \vec \xi, \vec A(\xi,\eta) \rangle$. Using the Moser regularization formula for $\vec \xi$ and Assumption (A2), we see that
\begin{equation}\label{eq:M}
\begin{split}
M&=-\xi_1 A_1(q_1(\eta_0,\xi_1,\eta_1),q_2(\eta_0,\xi_2,\eta_2))-\xi_2A_2(q_1(\eta_0,\xi_1,\eta_1),q_2(\eta_0,\xi_2,\eta_2))\\
&= \frac{2(p_1A_1(q_1,q_2)+p_2A_2(q_1,q_2))}{\Vert \vec p \Vert^2+1}.
\end{split}
\end{equation}
We then immediately see that
\begin{equation}\label{eq:M3}
M_{\xi_3}=M_{\eta_3}=0.
\end{equation}
For $j>0$, with the formulas for Moser regularization, we have
\begin{equation}\label{eq:Meta}
\begin{split}
M_{\eta_j} &=\sum_{i=1}^3\frac{\partial M}{\partial q_i} \frac{\partial q_i}{\partial \eta_j}+\frac{\partial M}{\partial p_i} \frac{\partial p_i}{\partial \eta_j} \\
&= (1-\xi_0)\frac{\partial M}{\partial q_j}.
\end{split}
\end{equation}
This vanishes on the collision locus where $\xi_0=1$.
For $M_{\eta_0}$ we get
$$
M_{\eta_0} = \sum_{i=1}^2 \frac{\partial M}{\partial q_i} \xi_i,
$$
which also vanishes on the collision locus (and so $M_\eta \cdot \xi$ also does). 
This second term hence reduces at $\xi_0=1$ to
\[
T_2=\Vert \eta \Vert^2f(1+b+M-M_\xi\cdot \xi)\xi_3^2,
\]
which is non-negative by Assumption (A4). It follows by compactness of $\Sigma$ that the coefficient of $T_2$ in front of $\xi_3^2$ is strictly positive for $\xi_0\geq 1-\delta,$ for $\delta$ sufficiently small.

We now deal with the third term $T_3$, for which we dissect some of its terms. 
Recall that $b(\eta,\xi)=-(c+1/2)+\frac{1}{2}\Vert \vec A(\eta,\xi)\Vert^2+V_1(\eta,\xi)$. Using Assumption (A1), we compute
\begin{equation*}
   \begin{split}
        b_{\xi_3}=&\sum_{i,j=1}^3A_i\frac{\partial A_i}{\partial q_j} \frac{\partial q_j}{\partial \xi_3}+\sum_{j=1}^3\frac{\partial V_1}{\partial q_j} \frac{\partial q_j}{\partial \xi_3}\\ 
        =& \eta_0\frac{\partial V_1}{\partial q_3}.
    \end{split}
\end{equation*}
Similarly,
\begin{equation*}
    \begin{split}
        b_{\eta_3}=&\sum_{i,j=1}^3A_i\frac{\partial A_i}{\partial q_j} \frac{\partial q_j}{\partial \eta_3}+\sum_{j=1}^3\frac{\partial V_1}{\partial q_j} \frac{\partial q_j}{\partial \eta_3}\\ 
        =& (1-\xi_0)\frac{\partial V_1}{\partial q_3}.
    \end{split}
\end{equation*}
Assumption (A2), implies that near $q_3=0$ we may write
$$
\frac{\partial V_1}{\partial q_3}=h(\vec q)q_3+\mathcal{O}(q_3^2)=h(\vec q)(\eta_0\xi_3+(1-\xi_0)\eta_3)+\mathcal{O}(q_3^2),
$$
where $h(\vec q)=\frac{\partial^2V_1}{\partial^2 q_3}\vert_{q_3=0}(\vec q)$. This implies that $\xi_3 b_{\xi_3}$ and $\eta_3b_{\eta_3}$ can be viewed as quadratic forms in $\eta_3, \xi_3$ (with non-constant coefficients). Also, note that if $\xi_0$ is uniformly close to $1$, then $q_3=\eta_0 \xi_3 + (1-\xi_0)\eta_3$ is uniformly close to zero along the compact manifold $\Sigma$, so that the above Taylor expansion holds for $\xi_0\geq 1-\delta$ for sufficiently small $\delta$.   

From the above discussion, we see that $\frac{T_3}{f \Vert \eta \Vert^2}$ can be written as a quadratic form in $\xi_3$ and $\eta_3$ near the collision locus, such that all its coefficients vanish at the collision locus. We conclude that we can write $T_3$ near $\xi_0=1$ as
$$
T_3 =f \Vert \eta\Vert^2 \cdot
\left(
\begin{array}{c}
\xi_3 \\
\eta_3
\end{array}
\right)^t
P(\eta,\xi)
\left(
\begin{array}{c}
\xi_3 \\
\eta_3
\end{array}
\right),
$$
where $P(\eta,\xi)$ is a symmetric $2\times 2$-matrix whose coefficients are smooth functions in $\xi$ and $\eta$, all which vanish at $\xi_0=1$. If $\delta>0$ is sufficiently small, $T_3$ is then dominated along $\xi_0\geq 1-\delta$ as a quadratic form by $T_1+T_2$ (which is positive definite). Therefore, 
$$
\Omega(X_Q)=\left(
\begin{array}{c}
\xi_3 \\
\eta_3
\end{array}
\right)^t
B(\eta,\xi)
\left(
\begin{array}{c}
\xi_3 \\
\eta_3
\end{array}
\right)
$$
for a matrix $B(\eta,\xi)$ that is positive definite for $\xi_0\geq 1-\delta$. We may then find the lower quadratic bound as stated in Claim 2, similarly as we did in the proof of Claim 1.

\medskip

To complete the proof, we need to fix the cutoff function. We choose $\epsilon$ and $\delta$ as in the above two claims and decrease $\epsilon$ such that $\epsilon \leq \delta/2$.
Choose a cutoff function $\rho_0$ depending only on $\xi_0$ such that
\begin{itemize}
    \item $\rho_0$ is non-decreasing,
    \item $\rho_0(\xi_0)$ vanishes for $\xi_0\leq 1-\delta$,
    \item $\rho_0(\xi_0)$ equals $1$ for $\xi_0 \geq 1-\epsilon$.
\end{itemize}
Define $\rho(\xi_0):= \frac{C_\epsilon}{ K } \rho_0(\xi_0)$, where $K=\max_\Sigma |d\rho_0(X_Q)|$, and note that $|d\rho(X_Q)|\leq C_\epsilon$, $|\rho|\leq C_\epsilon/K$. We choose $A_\delta>0$ in Claim 2 small enough so that $A_\delta\leq K/4$. We now evaluate $\Omega(X_Q)$. 

If $\xi_0< 1-\delta,$ we have $\Omega(X_Q)=\Omega_p(X_Q)\geq 0$ with equality only along $B$ (by Claim 1). If $\xi_0>1-\epsilon,$ we see that 
$$
\Omega(X_Q)=\Omega_p(X_Q)+\Omega_g(X_Q)\geq \Omega_g(X_Q)\geq A_\delta (\xi_3^2+\eta_3^2),
$$
by Claim 1 and Claim 2, and $\Omega(X_Q)$ vanishes if and only if $\Omega_p(X_Q)=\Omega_q(X_Q)=0$, which happens only along $B \cap \{\xi_0> 1-\epsilon\}$.
For the intermediate region,  $1-\delta\leq \xi_0 \leq 1-\epsilon$, we use that $2\xi_3\eta_3\leq \xi_3^2+\eta_3^2$, and we find
\begin{equation}
\label{eq:final_bound_angular_form}
\begin{split}
\Omega(X_Q)
&=
\Omega_p(X_Q) + \rho(\xi_0) \Omega_g(X_Q)+\xi_3 \eta_3 d\rho(X_Q) \\
&\geq 
\left(C_\epsilon-\frac{C_\epsilon}{K}A_\delta - \frac{C_\epsilon}{2} \right)(\xi_3^2+\eta_3^2)\geq \frac{C_\epsilon}{4} (\xi_3^2+\eta_3^2). 
\end{split}
\end{equation}
This verifies the assumptions of Lemma~\ref{lemma:Reeb}, and finishes the proof.
\end{proof}

\textbf{Upper bounds on return time.}
From the bound \eqref{eq:final_bound_angular_form} we will now deduce an upper bound on the return time, needed in order to extend the return map to the boundary. 

\begin{lemma}[Bounded return time]\label{lemma:boundedreturntime} Fix a page for the open book of Lemma \ref{lemma:GSS_Stark_Zeeman}. Then, the return time for the associated return map is uniformly bounded from above.
\end{lemma}
\begin{proof}
Let $\pi: \Sigma \setminus B \to S^1,\; x \mapsto \Theta(x)/\Vert\Theta(x)\Vert$ be the open book of Lemma \ref{lemma:GSS_Stark_Zeeman}. If we take standard angular coordinates $\phi$ on $S^1$, i.e.\
$$
\phi \circ\pi(x)=\atan \frac{\im \Theta(x)}{\re \Theta(x)},
$$
then we can compute for a flow line $x(t)$ of $X_Q$ the rate at which the angle progresses. This is
$$
\frac{d}{dt}\phi\circ \pi(x(t) \, )=\omega_{x(t)}(X_Q)= \frac{\Omega_{x(t)}(X_Q)}{\Vert \Theta(x(t))\Vert^2} \geq \frac{C_\epsilon(\xi_3^2 (t)+\eta_3^2 (t))}{4\left( \re \Theta(x(t))^2 +\im \Theta(x(t))^2 \right)}.
$$
The denominator can be bounded from above by a computation:
\[
\begin{split}
\re \Theta(x)^2 +\im \Theta(x)^2
&= 
\xi_3^2 + \left( 
(1-\xi_0)\eta_0\xi_3 + \left((1-\xi_0)^2+\rho\right)\eta_3
\right)^2 \\
&=\left( 1+(1-\xi_0)^2\eta_0^2 \right)\xi_3^2
+\left(
(1-\xi_0)^2+\rho
\right)^2\eta_3^2 \\
&\phantom{=}
+2(1-\xi_0)\eta_0\left( (1-\xi_0)^2+\rho \right) \xi_3 \eta_3.
\end{split}
\]
This can be written as a quadratic form in $\xi_3,\eta_3$, namely
\[
\left(
\begin{array}{c}
\xi_3 \\
\eta_3
\end{array}
\right)^t
\left(\begin{array}{cc}
1+ (1-\xi_0)^2\eta_0^2 & (1-\xi_0)\eta_0\left( (1-\xi_0)^2+\rho \right) \\
(1-\xi_0)\eta_0\left( (1-\xi_0)^2+\rho \right) & \left((1-\xi_0)^2+\rho \right)^2
\end{array}
\right)
\left(
\begin{array}{c}
\xi_3 \\
\eta_3
\end{array}
\right)
.
\]
We will bound the eigenvalues of the matrix from above. The determinant of the matrix is given by $ \left((1-\xi_0)^2+\rho \right)^2>0$ and its trace is 
$$
1+ (1-\xi_0)^2\eta_0^2 + \left((1-\xi_0)^2+\rho \right)^2.
$$
This means that the largest eigenvalue is bounded from above by 
$$
1+ (1-\xi_0)^2\eta_0^2 + \left((1-\xi_0)^2+\rho \right)^2 \leq 
1+4\eta_0^2 +(4+\rho)^2.
$$
The latter is also bounded on the compact hypersurface $\Sigma$, so we get a positive upper bound on the largest eigenvalue, say $\kappa$.
It follows that 
$$
\re \Theta(x(t))^2 +\im \Theta(x(t))^2 \leq \kappa (\xi_3^2(t)+\eta_3^2(t)).
$$
Hence 
$$
\frac{d}{dt}\phi\circ \pi(x(t) \, ) \geq \frac{C_\epsilon}{4 \kappa}.
$$
The return time is hence uniformly bounded from above by $ \frac{8\pi \kappa }{C_\epsilon}$.
\end{proof}

\section{The case of the restricted three-body problem}\label{sec:3bpcase}
The restricted three-body problem is an example of a Stark-Zeeman system, so we will only need to verify the conditions of Lemma~\ref{lemma:GSS_Stark_Zeeman}. With the expressions for $b$ and $M$ given in equations (\ref{eq:b}) and (\ref{eq:M1}), we find 
\[
\begin{split}
1+b+M-M_\xi \cdot \xi 
&=
1-(c+1/2)-\frac{1-\mu}{\Vert \vec \eta(1-\xi_0) +\vec \xi \eta_0+\vec m - \vec e \Vert}\\
&+(\xi_2\eta_1-\xi_1 \eta_2)(1-\xi_0) -\xi_2(1-\mu)-M_\xi\cdot \xi.
\end{split}
\]

Evaluated at $\xi=(\xi_0,\xi_1,\xi_2,\xi_3)=(1,0,0,0)$, using that $M_{\xi_0}=-(\xi_2\eta_1-\xi_1\eta_2)$, this reduces to
$$
1/2 -c - 1+\mu = \mu-c - 1/2.
$$

We conclude:
\begin{corollary}
\label{cor:RTBP_GHSS_low}
The spatial restricted three-body problem admits an $S^1$-family of global hypersurfaces of section for all energies $c$ with $c <H(L_1).$ 
\end{corollary}
\begin{proof}
Assumptions (A1) and (A2) were checked in Section \ref{sec:form3BP}, and Assumption (A3), in Remark \ref{rem:A3}. To check Assumption (A4), by the above, we only need to verify that  $\mu-c - 1/2>0$, and to see this, we use the fact that $H(L_1)\leq -3/2$; this can be deduced by combining Theorem~5.4.7 and Corollary~5.4.4 from \cite{FvK18}.
\end{proof}
\begin{remark}\label{rk:higherenergy}
Obviously, Assumption (A4) will also hold for higher energy.
However, in Lemma~\ref{lemma:GSS_Stark_Zeeman} we also need  the component $\Sigma$ of the once regularized hypersurface to be closed.
This condition does \emph{not} hold for large energies, i.e.~$c>H(L_3)$. 
\end{remark}

\textbf{Symmetries.} We now prove Proposition \ref{prop:symmetries} from the Introduction, which is a simple observation. 

\begin{proof}[Proof of Prop.\ \ref{prop:symmetries}]
The symplectic involution $r: T^*S^3 \rightarrow T^*S^3$ induced by the smooth reflection along the equator $S^2$ is, in regularized coordinates $(\xi,\eta)$, simply given by the restriction to $\mathbb S^*S^3$ of the map
$$
r\colon (\xi_0,\xi_1,\xi_2,\xi_3,\eta_0,\eta_1,\eta_2,\eta_3)\mapsto (\xi_0,\xi_1,\xi_2,-\xi_3,\eta_0,\eta_1,\eta_2,-\eta_3), 
$$
which flips the sign of the coordinates $\xi_3,\eta_3$. Moreover, it follows from Lemma \ref{thetap} and Equation \ref{thetamap} that $\Theta(r(\xi,\eta))=-\Theta(\xi,\eta)$ away from $B=\{\xi_3=\eta_3=0\}$, which is clearly the fixed point set of $r$. Similarly, the anti-symplectic involutions take the form
$$
\rho_1\colon (\xi_0,\xi_1,\xi_2,\xi_3,\eta_0,\eta_1,\eta_2,\eta_3)\mapsto (\xi_0,-\xi_1,\xi_2,\xi_3,-\eta_0,\eta_1,-\eta_2,-\eta_3), 
$$
$$
\rho_2\colon (\xi_0,\xi_1,\xi_2,\xi_3,\eta_0,\eta_1,\eta_2,\eta_3)\mapsto (\xi_0,-\xi_1,\xi_2,-\xi_3,-\eta_0,\eta_1,-\eta_2,\eta_3),
$$
and so $\Theta(\rho_1(\xi.\eta))=\overline{\Theta}(\xi.\eta)$, $\Theta(\rho_2(\xi.\eta))=-\overline{\Theta}(\xi.\eta)$, and the claim follows.
\end{proof}

\section{Second order estimates}\label{sec:secondorderestimates}

In this section, we carry out the estimates necessary to extend the return map to the boundary. These put the estimates of Lemma~\ref{lemma:boundedreturntime} in a more general setting.
%We consider the Hamiltonian $Q=\frac{1}{2}f^2\Vert \eta \Vert^2$ and 
We consider a Hamiltonian $Q$ on a symplectic manifold $(M,\omega)$, and we assume that $(\Sigma,\alpha)$ is a contact-type compact component of a level set of $Q$. 
Furthermore, we assume the following:
\begin{itemize}
\item $(B,\pi)$ is an open book for $\Sigma$ with adapted contact form $\alpha$. The trivialization of the normal bundle given in (2) of Definition~\ref{OBDS} will be denoted by $\varepsilon: B\times \mathbb{D}^2 \to \nu_\Sigma(B)$; we use $(x,y)$ to denote the coordinates on the $\mathbb{D}^2$-factor.
\item $B$ is invariant under the Reeb flow of $\alpha$ (which is a reparametrization of the flow of $X_Q$).
\item $(u_1=\frac{\partial}{\partial x},u_2=\frac{\partial}{\partial y})$ is a symplectic frame that trivializes the normal bundle $\nu_\Sigma(B)$. We will call this frame an \emph{adapted frame} since it is ``adapted'' to the binding of the open book.
\end{itemize}
We denote the coframe dual to $(u_1,u_2)$ by $(u^1,u^2)$.
Note that $u_1(Q)|_{B}=u_2(Q)|_{B}=0$, since $B$ is invariant under the flow of $X_Q$.
Consider the metric $g=\alpha\otimes \alpha+d\alpha|_{B}(\cdot,J\cdot)+u^1 \otimes u^1 +u^2 \otimes u^2$ for an extension of $u_1, u_2$ to a neighbourhood of $B$, where $J$ is a compatible almost complex structure on $\xi=\ker \alpha$ preserving $\xi\vert_B$.
Let $\nabla$ denote the Levi-Civita connection for this metric.
Choose any symplectic trivialization $\varepsilon_B$ of $\xi|_{B}$ along a Reeb trajectory $x$ in $B$; we denote the associated frame of $\xi|_{B}$ by $\{ e_i \}$ and the coframe by $\{ e^i \}$.
This gives us the trivialization $\varepsilon_\xi =\varepsilon_B\oplus(u_1,u_2)$ of $\xi$ along $x$.

\begin{lemma}
\label{lemma:normal_hessian}
With respect to the trivialization $\varepsilon_\xi$, the Hessian $Hess(Q):=\nabla dQ|_\xi$ has the block form
\[
Hess(Q)=\left(
\begin{array}{cc}
S_\xi & 0 \\
0 & S_N
\end{array}
\right)
,\text{ where }S_\xi \in \xi^* \otimes \xi^* \text{ and }S_N \in N^* \otimes N^*.
\]
\end{lemma}
\begin{proof}
Along $B$ we have $du^i=0$, so from the first structure equation we see that $\nabla u^i=0$ and $\nabla u_i =0$.

With the Einstein summation convention, we compute the Hessian:
\[
\begin{split}
\nabla dQ|_{\xi}
&=\nabla (e_i(Q)e^i+u_k(Q)u^k)\\
&=e_j e_i(Q) e^j \otimes e^i + u_k e_i(Q) u^k \otimes e^i
+e_i(Q) \nabla e^i
+e_j u_k(Q) e^j \otimes u^k\\
&\phantom{=}
+u_\ell u_k(Q) u^\ell \otimes u^k
+u_k(Q) \nabla u^k.
\end{split}
\] 
To prove the lemma we need to show that there are no $u^k\otimes e^i$-terms (in any order).
Since $u_k(Q)|_{B}=0$ by our above observation, the mixed derivatives $e_i u_k(Q)$ and $u_k e_i(Q)$ vanish.
Furthermore, $\nabla u^k$ vanishes by the above.

That leaves the term $\nabla e^i$.
We have 
$$
0=\nabla( e^i(u_k))= (\nabla e^i) (u_k) +e^i(\nabla u_k)=(\nabla e^i) (u_k),
$$
which excludes the term $e^j \otimes u^k$, and we can use torsion freeness to show there is no term $u^k \otimes e^j$, either.
This establishes the claim.
\end{proof}
We will call the bilinear form $S_N$ in the above decomposition the \emph{normal Hessian}.

\begin{proposition}\label{proposition:convexity}
As in the above setup assume that $(\Sigma,\xi=\ker\alpha)$ admits an open book $(B,\pi)$ with an adapted contact form $\alpha$. Let $P:=\overline{\pi^{-1}(1)}$ be a page of the open book, and denote the return map of the Reeb flow $R_\alpha$ by $f:\mathrm{int}(P) \to \mathrm{int}(P)$.
Assume that we have an adapted frame $(U=u_1, V=u_2)$ and that the normal Hessian is positive definite.
Then $f$ extends smoothly to the boundary.
\end{proposition}

\begin{remark}
The positive-definite assumption has some similarities with the condition of \emph{dynamical convexity} in \cite{HWZ98}.
We use it here to get a strong twist around the binding.
\end{remark}

\begin{proof}
By smooth dependence on initial conditions, we know that the return map $f$ is a smooth map on $\mathrm{int}(P)$.
We now define an extension to the boundary.
Take $x_0\in B$ and a sequence $\{ x_n\}_{n=1}^\infty \subset \mathrm{int}(P)$ converging to $x_0$.

For each $x_n$ we get a flow line $x_n(t)=Fl^R_{t}(x_n)$ with a first return time $t_n=t_+(x_n)$.
As a first step, we need to show that the limit $\lim_{n\to \infty} t_n$ is a well-defined (independent of the sequence $x_n$), positive real number $t_0$.
This will give us a candidate extension $\bar f(x_0):=Fl^X_{t_0}(x_0)$.

Let $\theta$ denote the angular coordinate in the neighbourhood $B\times \mathbb{D}^2$ as in (2) of Definition~\ref{OBDS}.
Since $\alpha$ is an adapted contact form, we have 
$$
\frac{d}{dt}\theta (x_n(t)\,)>0.
$$
For large $n$, we can approximate $x_n(t)$ by the linearized Reeb flow along $x_0(t)=Fl^R_t(x_0)$. We do this with a vector field $X_n(t)$ along $x_0(t)$ via
$$
x_n(t)=\exp_{x_0(t)}(\epsilon_n X_n(t) ),
$$
where we have chosen $\epsilon_n$ such that $\Vert X_n(0) \Vert=1$ and we use the adapted metric $g$ as defined above.
We choose coordinate functions
$$
u:B\times \mathbb D^2 \longrightarrow \R, (b;x,y) \longmapsto x,\quad
v:B\times \mathbb D^2 \longrightarrow \R, (b;x,y) \longmapsto y.
$$
Then we may write $\theta= \atan(u/v)$, so we obtain the angular form
$$
\beta =d \theta=\frac{udv-vdu}{u^2+v^2}.
$$
We have
$$
\frac{d}{dt} u\circ x_n(t)= \epsilon_n du \frac{d}{dt}X_n(t) +o(\epsilon_n)=
\epsilon_n du \nabla_{X_n} R_\alpha +o(\epsilon_n), 
$$
so we can see the growth rate of the angular coordinate near $B$ just from the linearized Reeb flow equation $\frac{d X}{dt} =\nabla_X R_\alpha$.
We expand $\beta(\dot x_n)$ near the binding where we write the $u,v$-component of $X_n$ as $X_n^u$ and $X_n^v$.
Since we are considering projections of Reeb orbits to $\mathbb D^2$, we can identify $T\mathbb D^2$ with $\R^2$, and expand the orbits in $\epsilon_n$.
This gives us the approximation $\pi_{\mathbb D^2}(x_n(t)\,)=\epsilon_n(X_n^u(t),X_n^v(t))+o(\epsilon_n)$, which we use to get
$$
\beta_{x_n}(\dot x_n)
=\frac{X_n^u dv(\nabla_{X_n}R_\alpha) - X_n^v du(\nabla_{X_n}R_\alpha)}{(X_n^u)^2+(X_n^u)^2}
+o(1).
$$
Since the Reeb flow is just a reparametrization of the Hamiltonian flow of $X_Q$, we now switch to the linearization of $X_Q$, which has the same qualitative behaviour, so we consider, with a little abuse of notation since we continue to use $X$, the linear ODE
$$
\frac{d X}{dt} =\nabla_X X_Q.
$$
Take the adapted frame $\{u_1=\frac{\partial}{\partial u},u_2= \frac{\partial}{\partial v} \}$ along $B$ and choose a symplectic trivialization $\epsilon_B$ as above along $x_0(t)$.
With respect to this trivialization we can write $X=(X_\xi,X_N)$.

By Lemma~\ref{lemma:normal_hessian} the linearized Hamiltonian flow splits into a $\xi|_{B}$-part and a normal part:
$$
\dot X_\xi =J_\xi S_\xi X_\xi,
\quad
\dot X_N =J_N S_N X_N.
$$
Note that we can write $\beta_X(Y)=\frac{X_N^t (-J_N) Y_N}{X_N^t X_N}$.
From the above expansion (using Hamiltonian flow rather than Reeb flow) we find
\[
\begin{split}
\lim_{n\to \infty} \beta_{x_n}(\dot x_n) & =\beta_X(\dot X) \\
&= \frac{X_N^t (-J_N)J_NS_N X_N}{X_N^t X_N} = \frac{X_N^t S_N X_N}{X_N^t X_N}.
\end{split}
\]
Since $B$ is compact, and $S_N$ is positive definite by assumption, we can bound the smallest eigenvalue of $S_N$ from below by $K>0$.
This gives a lower bound on the turning rate, i.e.~$d\theta(\dot X)\geq K$.
Since the linearized flow is smooth, it follows that there is a unique $t_0>0$ such that $\tilde \theta(X(t_0))=2\pi$ (where we write $\theta=e^{i\tilde \theta}$), and we have the bound $t_0 \leq 2\pi/K$.
By the linear approximation procedure we see that $\lim_{n \to \infty}t_n=t_0$, and because of the block form of the linearization  of Lemma~\ref{lemma:normal_hessian}, we see that this limit is independent of the sequence $\{ x_n \}$.

We now consider the first return time $t_+: P \to \R_{>0}$, and claim that this is a smooth function.
To see this, we take another point of view and adopt the arguments of Section~3.1 of \cite{H20}, which blows up the binding to verify smoothness.
To do so, we first identify a neighbourhood of the binding in $\Sigma$ with $\nu_\Sigma(B) =B\times D^2 \subset B \times \C$, and define the gluing map
$$
\Phi: B \times (0,1) \times S^1\subset B \times (-\infty,1) \times S^1
\longrightarrow
\nu_\Sigma(B) \cong B\times D^2,
(b,r,\phi) \longmapsto (b;re^{i\phi}).
$$
Then we blow up the binding by putting
$$
\Sigma_B :=
\Sigma \setminus B \coprod B \times (-\infty,1) \times S^1 / \sim,
$$
where we identify $(b,r,\phi)\in B \times (0,1) \times S^1\subset B \times (-\infty,1) \times S^1$ with $\Phi(b,r,\phi) \in \Sigma \setminus B$.\footnote{Note that only the part $B \times [0,1) \times S^1$ is directly relevant for the construction.}  
The map $\Phi^{-1}$ is well-defined and smooth on $B \times (D^2 \setminus \{ 0 \}) \subset \Sigma \setminus B$, so we can compute the pullback vector field $\Phi^* X_Q$.
The same computation as in Section~3.1 of \cite{H20} shows that $\Phi^* X_Q$ extends to a smooth vector field $\hat X_Q$ with the following properties:
\begin{itemize}
\item $\hat X_Q|_{\Sigma \setminus B} = X_Q|_{\Sigma \setminus B}$;
\item $\hat X_Q$ is tangent to $B\times \{0 \} \times S^1$;
\item the flow of $\hat X_Q$ is complete.
\end{itemize}
We also note that the open book $\pi: \Sigma \setminus B \to S^1$ extends to a fiber bundle $\hat \pi:\Sigma_B \to S^1$.

Define the function
$$
F: \Sigma_B \times \R \longrightarrow \R,\quad
(x,t) \longmapsto
\int_0^t d\hat \pi\left(\frac{d}{ds}(Fl^{\hat X_Q}_s(x))\right)ds.
$$
By the above argument and the computations we started with, the function $F$ is a smooth function of $x$ and $t$.
Furthermore it is an increasing function of $t$ on the subset 
$$
\{ x \in \Sigma_B~|~x\in \Sigma \setminus B \text{ or } x\in B\times \{0 \} \times S^1 \} \times \R\subset \Sigma_B\times \mathbb{R}
$$
with derivative bounded from below by a positive quantity.
The page $P$ embeds into $\Sigma_B$, so we can apply the implicit function theorem to the function $F|_{P\times \R}-2\pi$ to find a smooth function $t_+: P \to \R$ with the following properties:
\begin{itemize}
\item $F(p,t_+(p))=2\pi$ for all $p\in P$;
\item $t_+$ is positive (since $F(x,0)=0$).
\end{itemize}
This smooth function $t_+: P \to \R_{>0}$ equals the first return time on $\mathrm{int}(P)$, and equals the above first return time corresponding to the linearization of the Reeb flow on $B=\partial P$.
Since this function is smooth, the candidate extension $\bar f(x) =Fl^R_{t_+(x)}(x)$ is also smooth.
\end{proof}

In what follows, we will check the hypothesis of Proposition~\ref{proposition:convexity} for the restricted three-body problem. From Lemma \ref{lemma:boundedreturntime} for Stark-Zeeman systems, we see that we need only check the positive definite assumption. Rather than computing the Hessian along the normal direction to the binding, we will work directly with the linearized flow equation.  

\subsection{Round sphere} We first work out how the case for the round sphere $S^n$, using the notation from Example \ref{sec:geodesic_OB}. The symplectic normal bundle of $\mathbb S^*S^{n-1}\subset \mathbb S^*S^n$ has the following symplectic frame:
$$
U=\frac{\partial}{\partial \eta_n},\quad
V=\frac{\partial}{\partial \xi_n}.
$$
We will directly work out the equations for the linearized equation $\dot X=\nabla_X X_Q$ in terms of this frame and insert $\dot X$ into the angular form $\beta$, defined above.
We only need the $U,V$ part, since other components drop out (using that $\Vert \eta \Vert^2=1$). This is
\[
\begin{split}
\nabla_{uU+vV}X_Q
&=
u\nabla_U X_Q+v\nabla_V X_Q \\
&
=u \partial_{\xi_n}-2 u \eta_n (\xi \cdot \partial_\eta) -v \Vert \eta \Vert^2 \partial_{\eta_n}.
\end{split}
\]
On $B\cap \{\Vert \eta \Vert^2=1\}$, we find
$$
\beta(\dot X)=\frac{u^2 +v^2\Vert \eta \Vert^2}{u^2+v^2}
= 1.
$$

\subsection{The spatial three-body problem}
The same steps can be done for the spatial three-body problem, both in regularized and unregularized coordinates, as follows. 

\begin{comment}
\begin{remark} For future reference, we write the derivatives of $D$ and $f$:
\begin{equation}\label{Dderivatives}
    \begin{split}
       D_{\xi_0}&=-2(1-\xi_0)|\vec \eta|^2-2\eta_0 \vec \eta \cdot \vec \xi +2 \eta_1\\
       D_{\xi_1}&= 2\eta_0(\eta_1(1-\xi_0)+\xi_1\eta_0-1)\\
       D_{\xi_2}&= 2\eta_0(\eta_2(1-\xi_0)+\xi_2\eta_0)\\
       D_{\xi_3}&= 2\eta_0(\eta_3(1-\xi_0)+\xi_3\eta_0)\\
       D_{\eta_0}&= 2(1-\xi_0)\vec \eta \cdot \vec \xi +2\eta_0 |\vec \xi |^2 -2\xi_1\\
       D_{\eta_1}&= 2(1-\xi_0)(\eta_1(1-\xi_0)+\xi_1\eta_0-1)\\
       D_{\eta_2}&= 2(1-\xi_0)(\eta_2(1-\xi_0)+\xi_2\eta_0)\\
       D_{\eta_3}&= 2(1-\xi_0)(\eta_3(1-\xi_0)+\xi_3\eta_0)\\
    \end{split}
\end{equation}
\begin{equation}\label{fderivatives}
    \begin{split}
    f_{\xi_0}&=(\xi_1\eta_2-\xi_2\eta_1+c+1/2)+(1-\mu)\left(\frac{2D+(1-\xi_0)D_{\xi_0}}{2D^{3/2}}\right)\\
    f_{\xi_1}&=-(1-\xi_0)\eta_2+\frac{(1-\mu)(1-\xi_0)D_{\xi_1}}{2D^{3/2}}\\
    f_{\xi_2}&=(1-\xi_0)\eta_1-(1-\mu)+\frac{(1-\mu)(1-\xi_0)D_{\xi_2}}{2D^{3/2}}\\
    f_{\xi_3}&=\frac{(1-\mu)(1-\xi_0)D_{\xi_3}}{2D^{3/2}}\\
    f_{\eta_0}&=\frac{(1-\mu)(1-\xi_0)D_{\eta_0}}{2D^{3/2}}\\
    f_{\eta_1}&=(1-\xi_0)\xi_2+\frac{(1-\mu)(1-\xi_0)D_{\eta_1}}{2D^{3/2}}\\
    f_{\eta_2}&=-(1-\xi_0)\xi_1+\frac{(1-\mu)(1-\xi_0)D_{\eta_2}}{2D^{3/2}}\\
    f_{\eta_3}&=\frac{(1-\mu)(1-\xi_0)D_{\eta_3}}{2D^{3/2}}\\
    \end{split}
\end{equation}
\end{remark}
\end{comment}
\medskip

\textbf{Estimates in unregularized coordinates.} Recall that the unregularized Hamiltonian is
$$
H=\frac{1}{2}\Vert \vec p \Vert^2 -\frac{\mu}{\Vert \vec q-\vec m\Vert}-\frac{1-\mu}{\Vert \vec q-\vec e \Vert}+p_1q_2-p_2q_1,
$$
and so the Hamilton equations give
$$
\dot q_3 = p_3,
\quad
\dot p_3 = -q_3 \cdot
\left(\frac{\mu}{\Vert \vec q-\vec m \Vert^3 }+
\frac{1-\mu}{\Vert \vec q-\vec e \Vert^3}\right),
$$
with respect to the symplectic form $\omega=dp\wedge dq.$ We write the Hamiltonian vector field as
$$
X_H=X_H^B+\dot q_3 \frac{\partial}{\partial q_3}+\dot p_3 \frac{\partial}{\partial p_3},
$$
where $X_H^B$ is tangent to the planar problem $B=\{p_3=q_3=0\}$. The linearized flow satisfies the linear ODE
$$
\dot X=\nabla_X X_H,
$$
and we have the normal symplectic frame along $B$ given by
$$
U=\frac{\partial}{\partial p_3},
\quad
V=\frac{\partial}{\partial q_3}.
$$
We write $X=X_B+xU+yV$, where $X_B \in TB$. Note that the derivatives with respect to the $B$-coordinates $q_1,q_2,p_1,p_2$ of the $U$ and $V$ components $\dot q_3,\dot p_3$ of $X_H$ all vanish along $B$. Therefore, the same holds for the $U$, $V$ components of $\nabla_{X_B}X_H$.

We then compute that, along $B$, we have
$$
\nabla_XX_H\vert_B=\nabla_{X_B}X_H\vert_B+\nabla_{xU+yV}X_H\vert_B,
$$
with $\nabla_{X_B}X_H \vert_B\in TB$. Using that $\nabla_U\dot q_3=1$, $\nabla_V\dot p_3=-\left(\frac{\mu}{\Vert \vec q-\vec m \Vert^3 }+
\frac{1-\mu}{\Vert \vec q-\vec e \Vert^3}\right)$, and $\nabla_{xU+yV}X_H^B=0$, we see that
$$
\nabla_{xU+yV}X_H\vert_B= x\partial_{q_3}-y\left(\frac{\mu}{\Vert \vec q-\vec m \Vert^3 }+
\frac{1-\mu}{\Vert \vec q-\vec e \Vert^3}\right)\partial_{p_3}.
$$
For $\beta=\frac{p_3dq_3-q_3dp_3}{p_3^2+q_3^2},$ we obtain
$$
\beta(\dot X)=\frac{x^2+y^2\left(\frac{\mu}{\Vert \vec q-\vec m \Vert^3 }+
\frac{1-\mu}{\Vert \vec q-\vec e \Vert^3}\right)}{x^2+y^2}
$$

The relevant eigenvalues are $\lambda_1=1$ and $\lambda_2=
\frac{\mu}{\Vert \vec q-\vec m \Vert^3 }+
\frac{1-\mu}{\Vert \vec q-\vec e \Vert^3}$, both positive and bounded away from zero along $B$ (and non-singular away from collisions).

\medskip

\textbf{Estimates in regularized coordinates.} We finish the second order estimate by computing in regularized coordinates and checking the condition of Proposition~\ref{proposition:convexity} over the collision locus. We do our computations for any Hamiltonian of the form $Q=f^2\Vert \eta \Vert^2$, and check for the restricted three-body problem. The linearized flow equation is now
$$
\dot Y = \nabla_Y X_Q,
$$
where $X_Q$ is given by expression (\ref{Hamvf}). The following normal frame along the binding is
$$
U=\frac{\partial}{\partial \eta_3},
\quad
V=\frac{\partial}{\partial \xi_3}.
$$

We choose a metric for which $U,V$ are orthogonal to $B$, and its metric connection $\nabla$, satisfying $\nabla U=\nabla V=0$ along $B$. A straightforward computation gives:

\begin{lemma}\label{secondorderestimate}
Along $B$, we have

\begin{equation*}
    \beta(\dot Y)=\frac{(u,v) S (u,v)^t}{u^2+v^2},
\end{equation*}
where 
$$
S=\left(\begin{array}{cc}
    f(f+\Vert\eta\Vert^2f_{\eta_3\eta_3})& f\Vert\eta \Vert^2(f_{\xi_3 \eta_3}-f_\eta\cdot \xi) \\
    f\Vert\eta \Vert^2(f_{\xi_3 \eta_3}-f_\eta\cdot \xi) & f\Vert \eta\Vert^2 (f_{\xi_3\xi_3}+f+f_\eta \cdot \eta-f_\xi \cdot \xi) 
\end{array}\right)
$$
$\hfill \square$
\end{lemma}

For the three-body problem, after some computation using the expressions for $f$, we get the following explicit formulas for the entries of $S$:

$$
S_{11}=f\left(1+(1-\xi_0)(\xi_2\eta_1-\xi_1\eta_2-c-1/2)-\xi_2(1-\mu)+\frac{(1-\mu)(1-\xi_0)(\Vert \eta \Vert^2(1-\xi_0)^2-D)}{D^{3/2}}\right)
$$
$$
S_{22}=f\Vert \eta \Vert^2\left(\frac{(1-\mu)(1-\xi_0)}{D^{3/2}}(\eta_0^2+\eta_0\vec{\eta}\cdot \vec{\xi}+(1-\xi_0)|\vec{\eta}|^2-\eta_1)+\xi_2\eta_1-\xi_1\eta_2-c+1/2-\frac{1-\mu}{\sqrt{D}} \right)
$$
$$
S_{12}=S_{21}=-f\Vert \eta \Vert^2\left( \frac{(1-\mu)(1-\xi_0)}{2D^{3/2}}(2(1-\xi_0)(\vec{\eta}\cdot \vec{\xi}-\eta_0)+2\eta_0|\vec{\xi}|^2-2\xi_1)\right).
$$
Here,
$$
D
=
\left( \eta_{{3}} \left( 1-\xi_{{0}} \right) +\xi_{{3}}\eta_{{0}}
 \right) ^{2}+ \left( \eta_{{2}} \left( 1-\xi_{{0}} \right) +\xi_{{2}}
\eta_{{0}} \right) ^{2}+ \left( \eta_{{1}} \left( 1-\xi_{{0}} \right) 
+\xi_{{1}}\eta_{{0}}-1 \right) ^{2}.
$$

\medskip

We see that on the collision set $\xi_0=1$, we have $S_{12}=S_{21}=0$, and so $S$ is diagonal. Note that $\Vert \xi \Vert^2=1$ implies that $\xi_1=\xi_2=\xi_3=0$, and so $f\equiv 1$ and $D\equiv 1$ on this set. Therefore the eigenvalues of $S$ are
$$
S_{11}=1,\;S_{22}=\Vert \eta \Vert^2\left(-c-1/2+\mu\right).
$$
If we further restrict to the fixed energy level set $\Sigma_{c}=\{f^2\Vert \eta \Vert^2=\mu^2\}$, and we assume that $c<H(L_1)\leq -3/2$ (so that $\Sigma_{c}$ is closed and contact-type), we have $$S_{22}=\mu^2(-c-1/2+\mu)>0.$$

This finishes the second order estimate, and the proof that the return map extends to the boundary (in the case $c<H(L_1)$), claimed in Theorem \ref{thm:returnmap}.

\section{A global hypersurface of section for the connected sum}

\subsection{Two-center Stark-Zeeman systems}
We consider now a two-center spatial Stark-Zeeman system. Consider two distinct vectors $\vec e, \vec m$ in $\R^3$, and a Hamiltonian of the form
$$
H(\vec q,\vec p)
=\frac{1}{2}\Vert \vec p +\vec A(\vec q) \Vert^2
+V(\vec q),$$
where
$$
V(\vec q)=V_{\vec e}(\vec q)
+V_{\vec m}(\vec q)
+V_s(\vec q),
$$
with $V_{\vec e}(\vec q) =-\frac{g_1}{\Vert \vec q -\vec e\Vert}$,  $V_{\vec m}(\vec q) =-\frac{g_2}{\Vert \vec q -\vec m\Vert}$ for $g_1,g_2>0$, and $V_s$ a smooth function.
For very negative $c$, the sublevel set $\{ \vec q~|~V(\vec q) \leq c \}$ has two distinguished bounded components, which we denote by $C_{\vec e}(c)$ and $C_{\vec m}(c)$.
These are characterized by the property
\begin{itemize}
    \item $\vec e \in \overline{C_{\vec e}(c)}$ and $\vec m \in \overline{C_{\vec m}(c)}$.
\end{itemize}

We make the following assumption: 
\begin{enumerate}
    \item[(A5)] The effective potential $V=V_{\vec e}+V_{\vec m}+V_s$ has only finitely many critical points $L_1,\ldots, L_k$. Furthermore, the first critical value $c_1:=V(L_1)$ has only one preimage $L_1$ under $V$, and this is a critical point of index 1. 
    In addition, this critical point induces the boundary connected sum of the components $C_{\vec e}(c_1)$ and $C_{\vec e}(c_1)$. By this we mean, topologically, that $C_{\vec e}(c_1+\epsilon)= C_{\vec m}(c_1 +\epsilon):=C_{\vec e, \vec m}(c_1+\epsilon)\cong C_{\vec e}(c_1)\natural C_{\vec m}(c_1)$ for small $\epsilon>0$.
\end{enumerate}

In particular, we note that for $c<c_1$, the Hamiltonian $H$ restricts to a single center Stark-Zeeman system on each of the bounded components of the unregularized level sets 
$$
\Sigma_{\vec e,u}(c)=\{ (\vec q,\vec p) \in H^{-1}(c)~|~\vec q \in C_{\vec e}(c) \};
$$
$$
\Sigma_{\vec m,u}(c)=\{ (\vec q,\vec p)\in H^{-1}(c)~|~\vec q \in C_{\vec m}(c) \}.
$$
Smoothly, we have $\Sigma_{\vec e,u}(c)\cong S^2\times \mathbb{R}^3\cong \Sigma_{\vec m,u}(c)$. For $c=c_1+\epsilon>c_1$, the bounded component of the level set $H^{-1}(c)$, which projects to $C_{\vec e, \vec m}(c)$ and which we denote by $\Sigma_{\vec e,\vec m,u}(c)$, is topologically the connected sum $\Sigma_{\vec e,\vec m,u}(c)\cong \Sigma_{\vec e,u}(c_1)\# \Sigma_{\vec m, u}(c_1)$. The Hamiltonian flow is non-complete due to collisions at $\vec q=\vec e$ and $\vec q = \vec m$. We can perform Moser regularization to obtain regularized compact hypersurfaces $\Sigma_{\vec e}(c)$, $\Sigma_{\vec m}(c)$, $\Sigma_{\vec e,\vec m}(c)$, obtained by compactifying the corresponding unregularized versions, so that $\Sigma_{\vec e}(c) \cong S^2\times S^3\cong \Sigma_{\vec m}(c)$ for $c<c_1$, and $\Sigma_{\vec e,\vec m}(c)\cong \Sigma_{\vec e}(c_1)\#\Sigma_{\vec m}(c)$ for $c=c_1+\epsilon$. We make the following observations:
\begin{itemize}
    \item Near the collision locus of $\vec e$, Moser regularization leads to a Hamiltonian $Q_{\vec e}$, which has the same form as in Section \ref{sec:Moserregularization}, namely
$$
Q_{\vec e} = \frac{1}{2} f_{\vec e} (\xi,\eta)^2 \Vert \eta\Vert^2;
$$
\item Similarly, we obtain a Hamiltonian
$$
Q_{\vec m} = \frac{1}{2} f_{\vec m} (\xi,\eta)^2 \Vert \eta\Vert^2
$$
by performing  Moser regularization near $\vec m$.
\end{itemize}

For $c < c_1$, we have $\Sigma_{\vec e}(c)=Q_{\vec e}^{-1}(g_1^2/2)$ and $\Sigma_{\vec m}(c)=Q_{\vec m}^{-1}(g_2^2/2)$. There is also an induced Hamiltonian $Q_{\vec e,\vec m}$ on the boundary connected sum $\mathbb{D}^*S^3\natural \mathbb{D}^*S^3$, defined on a neighbourhood of $\Sigma_{\vec e,\vec m}(c)$, for which this hypersurface is as a level set, whenever $c=c_1+\epsilon$. The Hamiltonian $Q_{\vec e,\vec m}$ coincides with $Q_{\vec e}$ and $Q_{\vec m}$ in the corresponding summands.

The level set $\Sigma_{\vec e,\vec m}(c)$ has the following decomposition (\emph{not} a disjoint union):
\begin{itemize}
    \item A regularized neighbourhood $C_E$ of the collision point $\vec e$. Geometrically this is a neighbourhood of a Legendrian $2$-sphere (the collision locus);
    \item Similarly, a regularized neighbourhood $C_M$ of the collision point $\vec m$ and the corresponding collision locus;
    \item The unregularized level set $H^{-1}(c)$.
\end{itemize}

\begin{lemma}
\label{lemma:GSS_Stark_Zeeman_connected}
Assume that Hamiltonian $H$ is a Stark-Zeeman system satisfying (A1), (A2), (A3), and (A5). Write $f_{\vec e}=1+(1-\xi_0) b_{\vec e}+M_{\vec e}$, $f_{\vec m}=1+(1-\xi_0) b_{\vec m}+M_{\vec m}$. Let $c>c_1$, and $\Sigma_{\vec e,\vec m}(c)$ be the associated (connected) regularized energy hypersurface, which we assume to be closed. Further assume that:
\begin{assumption}
\begin{enumerate}
    \item[(A4')] $1+b_{\vec e}(\xi,\eta)+M_{\vec e}(\xi,\eta) -(M_{\vec e})_\xi(\xi,\eta) \cdot \xi>0$ for all $(\xi,\eta)\in \Sigma_{\vec e,\vec m}(c)$ with $\xi=(1,0,0,0)$, in regularized coordinates near $\vec e$; and similarly for $\vec m$.
\end{enumerate}
\end{assumption}

Then there exists an open book decomposition on $\Sigma_{\vec e,\vec m}(c)$, with binding the planar problem, so that each page is a global hypersurface of section for the Hamiltonian flow. 
\end{lemma}

\begin{proof} 
Fix $c>c_1$. By assumptions (A1) and (A2), we have an invariant set $B\subset \Sigma_{\vec e, \vec m}(c)$, corresponding to the planar problem. The key point now is that the construction of adapted open book of Lemma \ref{lemma:GSS_Stark_Zeeman} is independent on the actual energy value $c$ away from the collision locus (always coinciding with the physical open book, globally defined in unregularized coordinates), and this collision locus lies away from the region where the connected sum takes place. So we have a well-defined open book extending near the index $1$ critical point. Interpolating with the geodesic open book near the two connected components of the collision locus as in the proof of Lemma \ref{lemma:GSS_Stark_Zeeman}, we obtain the desired open book.  

Explicitly, we define the open book projection
\[
\begin{split}
\theta: \Sigma_{\vec e,\vec m}(c) \setminus B &\longrightarrow S^1\\
x &\longmapsto
\begin{cases}
\frac{\Theta_p}{|\Theta_p|}(x) & x\in \Sigma_{\vec e,\vec m}(c) \setminus (C_E \cup C_M) \\
\pi_{\vec e}(x) & x \in C_E \\
\pi_{\vec m}(x) & x \in C_M.
\end{cases}
\end{split}
\]
Here $\Theta_p$ is the map giving physical open book in regularized coordinates, as defined in Lemma~\ref{thetap}, and the maps $\pi_{\vec e}$ and $\pi_{\vec m}$ are the interpolated circle maps as defined in Lemma~\ref{lemma:GSS_Stark_Zeeman}. By the proof of Lemma~\ref{lemma:GSS_Stark_Zeeman}  (using assumptions (A3) and (A4')), the Hamiltonian vector field is transverse to the pages of the open book, so the claim follows.
\end{proof}

\subsection{Case of the three-body problem.} 
\label{sec:RTBP_connected}
We observe that  Assumption (A4') is satisfied provided $c<-1/2$.
Secondly, we note that $H(L_2)\leq -3/2 <-1/2$, see Section \ref{sec:3bpcase} and Remark \ref{rk:higherenergy}.
This means that the bounded component of the regularized $\Sigma_c$, which is diffeomorphic to the connected sum, admits an open book for all $c\in (H(L_1),H(L_2))$.
This finishes the proof of Theorem~\ref{thm:gss_openbooks}.

\begin{remark}
Topologically, each page of the open book for $\Sigma_{\vec e,\vec m}(c_1+\epsilon)$ is the boundary connected sum of the pages of the open books for $\Sigma_{\vec e}(c_1-\epsilon)$ and $\Sigma_{\vec m}(c_1-\epsilon)$ provided by Lemma \ref{lemma:GSS_Stark_Zeeman}. 
The monodromy is the composition of the corresponding monodromies (which commute); see below for more details.
\end{remark}

\section{Summary and details of the proofs of the main results}
\label{sec:summary_proofs}

\subsection{Proof of Theorem~\ref{thm:gss_openbooks}}
The existence of a global hypersurface of section in the case $c<H(L_1)$ follows from Corollary~\ref{cor:RTBP_GHSS_low}.
The ingredient for this corollary is Lemma~\ref{lemma:GSS_Stark_Zeeman}.
This lemma also contains the explicit form of the projection map $\Theta/|\Theta|$, in Equation~\eqref{thetamap}.
This formula tells us that $\xi_3=0$ is a global hypersurface of section.

The case $c\in (H(L_1),H(L_2))$ is handled by combining Lemma~\ref{lemma:GSS_Stark_Zeeman_connected} with the observation in the beginning of Section~\ref{sec:RTBP_connected}.

\subsection{Proof of Theorem~\ref{thm:openbooks}}
The assertion that hypersurface $\Sigma_c$ is a contact manifold follows from \cite{ChoKim}.
To obtain the statement about open books we apply Lemma~\ref{lemma:Reeb} to conclude that the open book $(B,\theta)$ on $\Sigma_c$ carries the underlying contact structure.
The statement about the monodromies can then be obtained as follows:

\subsubsection{Monodromy for $c< H(L_1)$}
We can homotope the Hamiltonian to the Kepler problem via
$$
H_s=\frac{1}{2}\Vert \vec p\Vert^2 - \frac{s+(1-s)\mu}{\Vert \vec q-\vec m\Vert } - (1-s)\frac{1-\mu}{\Vert \vec q-\vec e\Vert } +(1-s)(p_1q_2-p_2q_1). 
$$
These problems can all be regularized, leading to regularized Hamiltonians $Q_s$ on $T^*S^3$.
The level set of the Hamiltonian $Q_1$ corresponds to standard contact form on $S T^*S^3$ together with its natural open book.

We will compute the monodromy using a Lefschetz fibration on the natural filling of $ST^*S^3$, namely $D^*S^3$.
Put
\[
\begin{split}
\Theta: D^*S^3=\{ (\xi,\eta) \;|\; \Vert \xi\Vert^2=1, \langle \xi,\eta \rangle=0,~\Vert \eta \Vert^2 \leq 1 \} &\longrightarrow D^2 \subset \C,\\
(\xi,\eta) &\longmapsto \xi_3+i\eta_3.
\end{split}
\]
We claim that $\Theta$ defines a Lefschetz fibration.
With a computation, we see that the only critical points are at $\xi =(0,0,0,\pm 1)$ and $\eta=0$; expanding $\Theta$ in local coordinates shows that these critical points are of Lefschetz type, i.e.~they are non-degenerate quadratic.
The generic fiber is symplectically deformation equivalent to $D^*S^2$.
Because there are only two critical points, the holonomy of this Lefschetz fibration is the product of two positive Dehn twists, say $\tau_1$ along $L_1$ and $\tau_2$ along $L_2$. 
The corresponding Lagrangian spheres $L_1$ and $L_2$ of each Dehn twist are the zero-section of $D^*S^2$.
Finally, we recall that an exact Lefschetz fibration on $W$ induces an open book on the boundary of $W$ and the monodromy of this open book equals the isotopy class of the holonomy. The claim follows.

\subsubsection{Monodromy for the case $c\in ( H(L_1), H(L_1)+\epsilon)$ and $\mu \in(0,1)$}
We will also construct a Lefschetz fibration.
The essential idea is to glue two copies of the standard Lefschetz fibration on $D^*S^3$ together using several cutoff functions. 
Let us first define the symplectic manifold of interest using a cover.
Write the first Lagrange point $L_1=(\ell_1,p_{\ell_1})$.
Choose $\delta_1>0$ such that
\begin{itemize}
\item the first coordinate satisfies $(\ell_1-\vec m)_1>2\delta_1$.
\item the first coordinate satisfies $(\vec e-\ell_1)_1>2\delta_1$.
\end{itemize}
Put the physical region
$$
W_p:=\{
(\vec q,\vec p)\in T^*\R^3~|~q_1\in( (\ell_1)_1-\delta_1,(\ell_1)_1+\delta_1),~H(\vec q,\vec p)\leq c
\}
.
$$
Intuitively, this region contains points in the unregularized phase-space close to the Lagrange point, which stay away from both the Earth and Moon (by a distance at least $\delta_1$).
Next define
$$
W_{\vec m}:=\left\{ (\xi,\eta) \in T^*S^3\cap \Sigma_c~|~Q^{\vec m}_{\mu,c}(\xi,\eta) \leq \frac{\mu^2}{2},~q_1(\xi,\eta)\leq (\ell_1)_1-\delta_1 \right\}
.
$$
This set contains points in the bounded (due to the intersection with $\Sigma_c$) regularized component that are close to the Moon and that stay away from the Lagrange point $L_1$.
The Hamiltonian $Q^{\vec m}_{\mu,c}$ is the Moser regularization near $\vec m$ of the Jacobi Hamiltonian for parameters $\mu,c$.
Finally, we put 
$$
W_{\vec e}:=\left\{ (\xi,\eta) \in T^*S^3\cap \Sigma_c~|~Q^{\vec e}_{\mu,c}(\xi,\eta) \leq \frac{(1-\mu)^2}{2},~q_1(\xi,\eta)\geq (\ell_1)_1+\delta_1 \right\}
.
$$
The Hamiltonian $Q^{\vec e}_{\mu,c}$ is the Moser regularization near $\vec e$ of the Jacobi Hamiltonian for parameters $\mu,c$, so $W_{\vec e}$ consists of point in the bounded regularized component near $\vec e$.

Put $W:=W_{\vec m} \cup W_p \cup W_{\vec e}$. This is a symplectic manifold with boundary, and it is clear that $W$ is diffeomorphic to the boundary connected sum of two copies $D^*S^3$.
We define a Lefschetz fibration of $W$ over disk in $\C$ by the map
\[
\Theta(x)
=
\begin{cases}
\Theta_g(\xi,\eta)+i\delta_{\vec m}(\xi,\eta) &\text{if } x=(\xi,\eta)\in W_{\vec m}, \\ 
(-p_3+iq_3) \cdot (1-\rho_{\vec m} \xi_0^{\vec m}-\rho_{\vec e} \xi_0^{\vec e}) & \text{if } x=(p_3,q_3) \in W_p, \\
\Theta_g(\xi,\eta)+i\delta_{\vec e}(\xi,\eta) &\text{if } x=(\xi,\eta)\in W_{\vec e}.
\end{cases}
\]
Before we explain the terms in more details, let us briefly give the intuition. On $W_{\vec m}$ we are taking the geodesic Lefschetz fibration that we have also used for the case $c<H(L_1)$, that we modify to interpolate the ``physical'' Lefschetz fibration $-p_3+iq_3$, which we use on $W_p$.
The latter map is also modified by a cutoff function to guarantee smoothness.
On $W_{\vec e}$ we use again the geodesic Lefschetz fibration.
Let us explain the above terms:
\begin{itemize}
\item the term $\delta_{\vec m}$ is
$$
\delta_{\vec m}=( (1-\xi_0)\eta_0 \xi_3+(-2\xi_0+\xi_0^2)\eta_3 )\cdot \rho_1(\xi_3),
$$
where $\rho_1$ is a cutoff function that vanishes near $\xi_3 =\pm 1$, and that equals $1$ for $|\xi_3|<1-\delta'$.
\item the term $\delta_{\vec e}$ has the same expression, and is defined on $W_{\vec e}$.
\item $\xi_0^{\vec m}$ is the $\xi_0$ coordinate (as a function of $\vec p$) as defined via Moser regularization at $\vec m$, and
$\xi_0^{\vec e}$ is the $\xi_0$ coordinate (as a function of $\vec p$) as defined via Moser regularization at $\vec e$.
\item $\rho_{\vec m}$ and $\rho_{\vec e}$ are cutoff functions depending only on $q_1$ with the following property. The function $\rho_{\vec m}$ vanishes for $q_1>(\ell_1)_1-\delta_1/2$ and equals $1$ near $q_1=(\ell_1)_1-\delta_1$. Similarly, the function $\rho_{\vec e}$ vanishes for $q_1<(\ell_1)_1+\delta_1/2$ and equals $1$ near $q_1=(\ell_1)_1+\delta_1$.  
\end{itemize}
These choices guarantee that $\Theta$ is a smooth function on $W$.
We leave out the details, and refer to the expressions in Section~\ref{thetap} to perform the necessary computations.

We now look at the critical points of $\Theta$.
We find that $\Theta$ has two critical points in $W_{\vec m}$, again corresponding to $\xi=(0,0,0,\pm 1)$ and $\eta=0$, two critical points in $W_{\vec e}$ (of the same form), and no critical points in $W_p$. Using the arguments for $c<H(L_1)$, and the fact that the Lagrangian zero sections in $W_{\vec m}$ and $W_{\vec m}$ are disjoint, we find that the monodromy is a composition double Dehn twists along each of these two Lagrangian zero-sections.
This establishes the claim.

\begin{remark}
The map $\Theta$ we defined in the above, is not equal to the formula we used in \eqref{thetamap}.
However, the corresponding bindings $\Theta^{-1}(0)$ are equal, and we may see Formula~\eqref{thetamap} as a reparametrization of $\Theta$.
\end{remark}

\subsection{Argument for Theorem~\ref{thm:returnmap}}
The extension of the return map to the boundary of $P$ follows from Proposition~\ref{proposition:convexity} and the estimates in Section~\ref{sec:secondorderestimates}.
To show that the return map is Hamiltonian, we use Lemma~\ref{lemma:monodromyvsreturnmap} from Appendix~\ref{monodromyvsreturnmap}. 
 to show that the return map can be Hamiltonianly isotoped to a representative of the monodromy, via an isotopy which extends smoothly to the boundary, which is either a double Dehn twist or the composition of two of those along disjoint Lagrangian. 
Since a double Dehn twist along a Lagrangian sphere is actually Hamiltonian, and the isotopy extends smoothly to the boundary, the first claim in Theorem~\ref{thm:returnmap} follows.

For the second claim, we construct the diffeomorphism $G: \mbox{int}(P) \rightarrow \mbox{int}(\mathbb{D}^*S^2)$, extending continuously to the boundary but whose inverse has a smooth extension to it, as in the the proof of Lemma \ref{lemma:monodromyvsreturnmap} in Appendix \ref{monodromyvsreturnmap}. Namely, we first find a collar neighbourhood $B\times \mathbb{D}^2$ in which the contact form looks like $\alpha=A(\alpha_P+r^2d\theta)$; we then use Moser's trick to construct a symplectomorphism $\psi_1:(P,d\alpha\vert_P)\rightarrow (P,\omega_1)$ which is the identity at $B$ and supported near it, and $\omega_1=d((1-r^2)\alpha_P)$ near $B$; we then construct a square root map $Q:(\mathbb{D}^*S^2,\omega_Q)\rightarrow (P,\omega_1)$, which is smooth away from the boundary and only continuous along it, and supported near it, and where $(\mathbb{D}^*S^2,\omega_Q)$ is an honest Liouville filling of $(\mathbb{R}P^3,\alpha_P)$. We may then take $G:=Q^{-1}\circ \psi_1:(P,d\alpha\vert_P)\rightarrow (\mathbb{D}^*S^2,\omega_Q)$. In the case of the Kepler problem, this construction yields $\omega_Q=\omega_{std}\vert_{\mathbb{D}^*S^2}$; by deforming to this problem we see that $\omega_Q$ is always deformation equivalent to $\omega_{std}\vert_{\mathbb{D}^*S^2}$. The rest of the second claim is then immediate from the first claim.

\appendix

\section{Return map for the rotating Kepler problem}\label{app:rotKepler}
In this appendix, we illustrate how to understand the qualitative dynamics for the rotating Kepler problem via a global hypersurface of section. This is a completely integrable system, obtained as the limit of the restricted circular three-body problem by setting $\mu=1$, for which the return map can be written down explicitly. 

In unregularized coordinates the rotating Kepler problem is described by the Hamiltonian $H=K+L$,\footnote{We are using a different convention here than Equation~\eqref{eq:Jacobi_Hamiltonian}, because we will be using the physical interpretation in terms of angular momentum.} where
$$
K= \frac{1}{2} \Vert \vec p \Vert^2  -\frac{1}{\Vert\vec q\Vert},
\quad
L= q_1 p_2 - q_2 p_1.
$$
After Moser regularization, the Hamiltonians $K$ and $L$ both generate circle actions and are in involution, which implies that they are preserved quantities of the motion. The remaining integral is the last component of the Laplace-Runge-Lenz vector. The regularized Hamiltonian is $Q(\xi,\eta)=\frac{1}{2}f^2(\xi,\eta)\Vert \eta \Vert^2$ where $f(\xi,\eta)=1+(1-\xi_0)(-c-1/2+\xi_2\eta_1-\xi_1\eta_2)$, obtained from Equation (\ref{eq:f3bp}) by setting $\mu=1$. 

Instead of using the general open book for Stark-Zeeman systems, we consider the geodesic open book 
$$
(
\xi, \eta) \longmapsto \frac{\xi_3 +i\eta_3}{\Vert \xi_3+i\eta_3 \Vert}.
$$
\begin{lemma}\label{app:lemmarotKepler}
The geodesic open book is a supporting open book for the rotating Kepler problem for $c<-3/2$.
\end{lemma}

\begin{proof}
A geometric way of seeing this is observing that the pages of the geodesic open book, which is adapted to the Kepler problem $K$, are also invariant under the Hamiltonian flow of $L$ (which acts by rotation along the $(\xi_0,\eta_0)$-axis inside a given page; recall Figure \ref{fig:Birkhoff_annulus}). Since $\{K,L\}=0$, we have $\phi^t_H=\phi^t_K \circ \phi^t_L$. While $\phi^t_L$ leaves the pages invariant, the flow $\phi^t_K$ is transverse to them. This implies the claim. This also implies that the return maps associated to different pages have the same exact dynamics, only that the return map on the $t$-page is rotated by an angle $t$ along the $\xi_0$-axis. \end{proof}

From the proof of the above lemma, in order to study the return map, we see that it suffices to consider the page $$
P= \left\{ (\xi;\eta) \in T^*S^3:\;Q(\xi,\eta)=\frac{1}{2},\; \xi_3=0,\; \eta_3 \geq 0 \right\},
$$ which is the easiest to visualize (see Figure \ref{fig:Birkhoff_annulus} in the Introduction). 

Since orbits in the rotating Kepler problem are precessing ellipses, the return time in unregularized coordinates is simply the minimal Kepler period for the corresponding Kepler energy $K$ (which is preserved under the flow of $H$). By Kepler's third law, this return time depends only $K$, and is given by 
$$
T = T(K)=\frac{\pi}{2(-K)^{3/2}}.
$$
From $\phi_K^T=id$, we obtain $\phi^T_H=\phi^T_K \circ \phi^T_L=\phi^T_L$, and therefore we derive the following return map
\[
R: (q_1,q_2,q_3;p_1,p_2,0) \longmapsto
( Rot_{\frac{\pi}{2(-K)^{3/2}}}(q_1,q_2),q_3; Rot_{\frac{\pi}{2(-K)^{3/2}}}(p_1,p_2),0),
\]
where $Rot_\phi$ is the rotation by angle $\phi$. 
This is generated by the Hamiltonian $L$ restricted to the global hypersurface of section (as the time $T$-map). To obtain an explicit formula for a Hamiltonian generating $R$ in time-$1$, we manipulate the above expression for $R$. With the relation $K+L=c$, we can also write 
$$
R=\phi^{X_L}_{T(K)}=\phi^{T(K)X_L}_{1}=\phi^{T(c-L)X_L}_1.
$$
We see that there is a function $g(L)$ such that $R=\phi^{X_{g(L)}}$, by noting that $X_{g(L)}=g'(L)X_{L}$. With $g'(L)=\frac{\pi}{2(L-c)^{3/2}}$, we can compute $g(L)$ as
$$
g(L)=-\pi\left( 2(L-c)\right)^{-1/2}.
$$
We may now describe the return map in the Moser regularized coordinates, which is given by
$$
R_r:(\xi_0,\xi_1,\xi_2,0;\eta_0,\eta_1,\eta_2,\eta_3) \longmapsto
(\xi_0,Rot_{T(c-L)}(\xi_1,\xi_2),0;\eta_0,Rot_{T(c-L)}(\eta_1,\eta_2),\eta_3).
$$
The Hamiltonian $L$ is given in these coordinates by $L=\xi_2 \eta_1 -\xi_1 \eta_2$.
\subsection{Dynamical consequences} We now explore the consequences of the above explicit description.

\medskip

\textbf{Polar orbits.} The two (distinct) periodic polar orbits are clearly visible as two fixed points of the return map: these are $x_\pm=(\pm 1,0,0,0;0,0,0,\eta_3)$, where $\eta_3>0$ is such that they lie on the level set $Q^{-1}(\frac{1}{2})$. The fixed point $x_-$ is the starting point of the vertical periodic collision orbit that lies in the upper half-space $q_3\geq 0$; in unregularized coordinates it is the point in the $q_3$-axis that is maximally far from the origin.
The fixed point $x_+$ is the starting point of the vertical periodic collision orbit that is contained in the lower half-space, and this fixed point corresponds to the periodic collision point.
This orbit and nearby periodic orbits for $\mu\neq 1$ were already studied by Belbruno in \cite{B}.

\begin{figure}
    \centering
    \includegraphics[width=0.7 \linewidth]{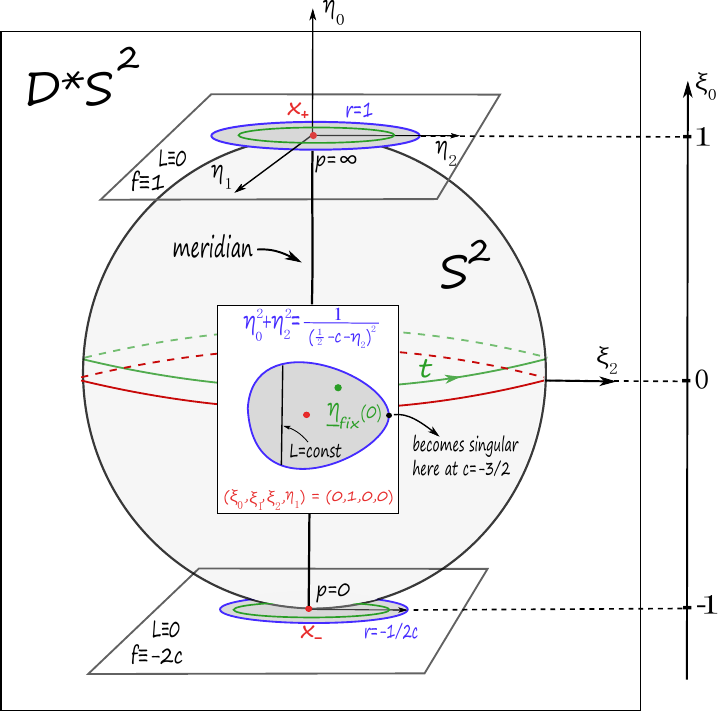}
    \caption{The invariant circles inside the page $P$ (in green).}
    \label{fig:rotatingkepler}
\end{figure}

\medskip

\textbf{Invariant circles.} Other than that, there are lots of invariant circles. For instance, writing $\eta=(\underline{\eta},\eta_3)$ where $\underline{\eta}=(\eta_0,\eta_1,\eta_2)$, these are given by
\begin{equation}
    \begin{split}
C_{x,\underline{\eta}_{fix}}=\Bigg\{&  (\xi;\eta) \in Q^{-1}(1/2):\;\xi_0\equiv x,\;\xi_1^2+\xi_2^2=1-x^2,\;\xi_3=0,\\& \underline{\eta}=R_t(\underline{\eta}_{fix}(x)),\; t \in [0,2\pi],\;\eta_3=\sqrt{\left(f^{-2}(\xi,\eta)\vert_{\xi_0=x}-\Vert \underline{\eta}_{fix}(x)\Vert^2\right)}\geq 0\Bigg\},
    \end{split}
\end{equation}
for given $x\in [-1,1]$, where $\underline{\eta}_{fix}:[-1,1]\rightarrow T^*S^2$ is a section of the cotangent bundle along the meridian passing through $(\xi_0,\xi_1,\xi_2)=(0,1,0):=\underline{\xi}_{fix}$ satisfying $\Vert\underline{\eta}_{fix}(x)\Vert^2 f^2\vert_{\xi_0=x}\leq 1$ (i.e. taking values in $P\subset T^*S^2$), and where $R_t$ is the rotation of angle $t$ along the $\eta_0$-axis. We parametrize this meridian by $x \mapsto \left(x,\cos(2\pi x), \sqrt{\sin(2\pi x)^2-x^2}\right)$. Note that the fibers of $P$ are all invariant under rotation in the $\xi_0$-axis (i.e.\ under the Hamiltonian action of $L$), since $f$ depends only on the angular term $L$ and on $\xi_0$. 
In particular, the fiber over the north pole is a disk of radius $1$, and the fiber on the south pole is a disk of radius $0<-\frac{1}{2c}<1/3$. The singular cases, for which the circles collapse to a point, correspond to $(x=\pm 1,\underline{\eta}_{fix}(\pm 1)=0)$, in which case $C_{\pm 1,\underline{\eta}_{fix}(\pm 1)}=\{x_\pm\}$ are the polar fixed points. See Figure \ref{fig:rotatingkepler}. One may also vary $x$ to obtain invariant annuli
$$
A_{\underline{\eta}_{fix}}=\bigcup_{x \in [-1,1]} C_{x,\underline{\eta}_{fix}},
$$
which in singular cases collapses to a sphere or a disk. 

The action of $R_r$ on such a circle $C_{x,\underline{\eta}_{fix}}$ is by rotation by angle $T(c-L)$. Note that $L$ vanishes at $\xi_0=\pm 1$, and so $R_r$ acts by rotation by angle $T(c)$ on the fibers of the north and south poles. This immediately implies that $x_\pm$ (and the associated orbits) are elliptic, and they are degenerate whenever $T(c)$ is an integer multiple of $2\pi$. Also note that $L\vert_{\underline{\xi}_{fix}}=-\eta_2$, and so $L$ is constant in the vertical lines of the fiber at $\underline{\xi}_{fix}$. Therefore, having fixed $c$, and noting that the function $T(c-L)=\frac{\pi}{2(L-c)^{3/2}}$ is injective in the region $K=c-L<0$ as a function of $L$, we may vary $\underline{\eta}_{fix}(0)$ transversely to these level sets in order to achieve the resonance condition $T(c-L)=2\pi\frac{p}{q}$ for some co-prime integers $p,q$. For such a choice, every point in $C_{0,\underline{\eta}_{fix}}$ is periodic of period $q$. Thus we obtain infinitely many periodic points of arbitrary large period lying over the equator $\{\xi_0=0\}$. 
We may play the same game for different values of $\xi_0$, to obtain infinitely many periodic points of arbitrary large period lying over all parallels.
Whenever $\underline{\eta}_{fix}(x)$ lies in the boundary of the corresponding fiber, fixed points for $R_r$ in the associated circle gives planar orbits; the spatial orbits are detected precisely when $\underline{\eta}_{fix}(x)$ lies in the interior.

\begin{figure}
    \centering
    \includegraphics[width=1.0 \linewidth]{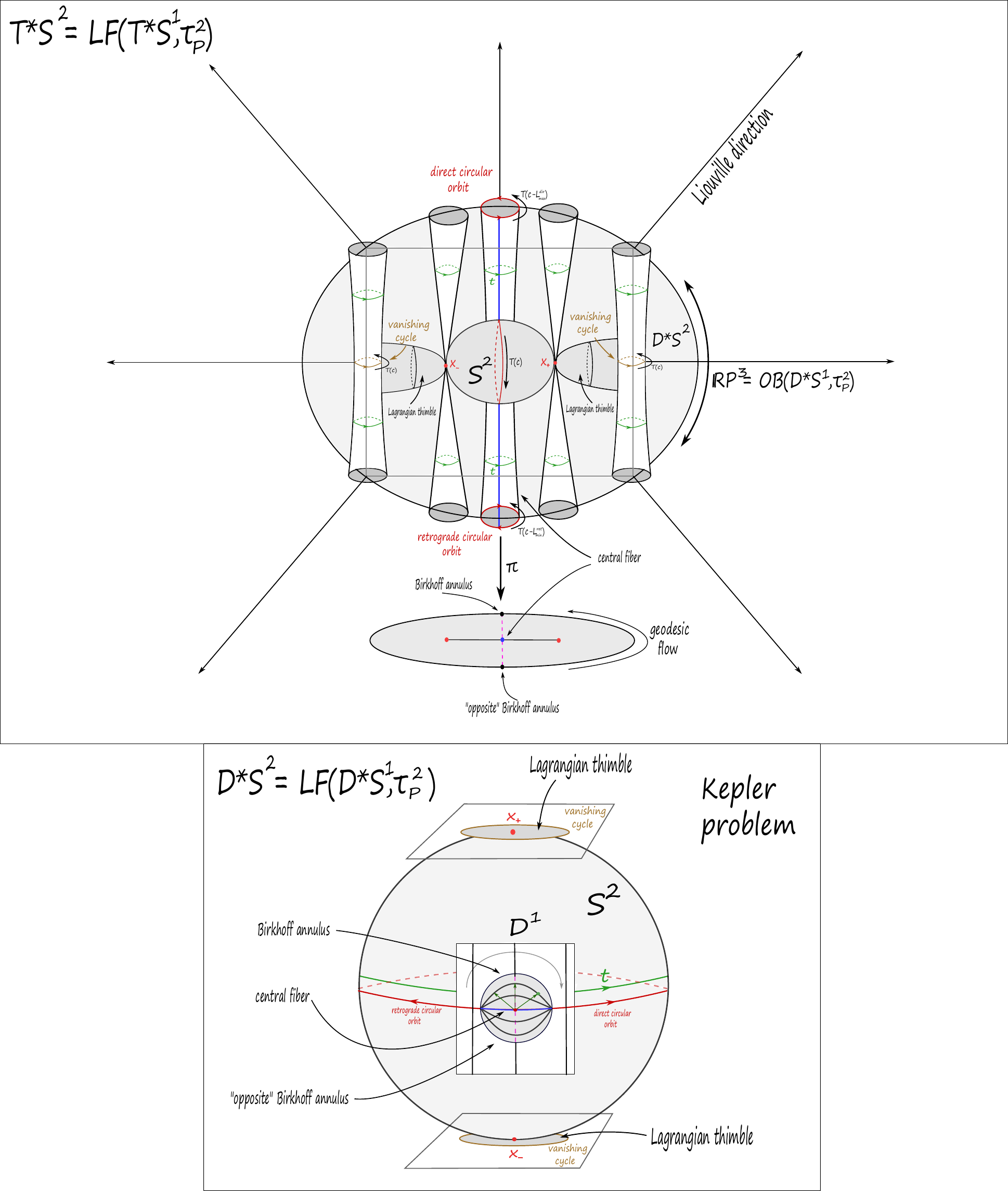}
    \caption{The ideal (i.e.\ non-compact) version of the modified standard symplectic Lefschetz fibration $\pi: T^*S^2\rightarrow \mathbb{C}$ is shown above, where we have isotoped it so that the direct/retrograde circular planar orbits $\gamma_{ret}, \gamma_{dir}$ lie over the origin, and hence looks standard. The return map is a classical twist map on the fibers, which are non-compact and asymptotic to $\gamma_{ret}\cup \gamma_{dir}$ at infinity. To obtain the Lefschetz fibration on $\mathbb{D}^*S^2$, having compact fibers, one has to project in the Liouville direction. The figure below shows the standard (unmodified) Lefschetz fibration $\pi_0$ after projection in the Liouville direction, which is the one suitable for the non-rotating Kepler problem, and the role of direct/retrograde orbits is played by the equator traversed in both directions. By construction of the modified version below, which is adapted to the rotating Kepler problem, this compact standard version can also be deformed to the compact modified version.}
    \label{fig:LF}
\end{figure}

%Moreover, if $-1\leq x_1<x_2 \leq 1$, and we consider $\underline{\eta}_{fix}\vert_{[x_1,x_2]}$, and we further choose the endpoints $\underline{\eta}_{fix}(x_i)$ so that the angles of rotation are different (only possible if either $x_1> -1$ or $x_2<1$), the annulus $A_{\underline{\eta}_{fix}}^{x_1,x_2}:=\bigcup_{x \in [-x_1,x_2]} C_{x,\underline{\eta}_{fix}}\subseteq A_{\underline{\eta}_{fix}}$ is invariant under $R_r$, which acts as a twist map. So the situation is indeed compatible with the classical Poincar\'e-Birkhoff theorem. 

\medskip

\textbf{Liouville tori.} We may also understand the (3-dimensional) Lagrangian Liouville tori in the ambient phase-space. We consider the geodesic flow on the page $P$, which is given by the Hamiltonian flow $\phi_K^t: P \rightarrow P$ of $K\vert_P$, and denote
$$
T_{x,\underline{\eta}_{fix}}=\bigcup_{t\in \mathbb{R}}\phi_K^t(C_{x,\underline{\eta}_{fix}}).
$$
Generically, this is a two-torus lying in $P$, and it is a Liouville torus for the planar problem whenever $\underline{\eta}_{fix}(x)$ lies in the boundary. 
There are singular cases, too.
For instance, one case correspond to $\underline{\eta}_{fix}(x)=0$, which gives a circle over a parallel, point-wise fixed by the action of $X_K$; in particular when $C_{x,\underline{\eta}_{fix}}=\{x_\pm\}$, for which $T_{x,\underline{\eta}_{fix}}=\{x_\pm\}$, a point fixed by the flow of $X_K$ and $X_L$;
The Liouville tori $L_{x,\underline{\eta}_{fix}}$ for the spatial problem are obtained by spinning $T_{x,\underline{\eta}_{fix}}$ around the open book $S^1$-direction, so that the ones corresponding to $x_\pm$ are precisely the associated collision orbits. The planar tori are therefore singular cases of the $T_{x,\underline{\eta}_{fix}}$, since the open book direction is not defined along the binding.

Note that the frequency corresponding to the $S^1$-direction of the Liouville tori $L_{x,\underline{\eta}_{fix}}$ given by the $X_K$-action is zero. This implies that, in action-angle coordinates, the Hessian with respect to action variables of the Hamiltonian is degenerate, and so the original version of the KAM theorem does not apply. On the other hand, the circles $C_{x,\underline{\eta}_{fix}}$ for which the resonance condition is satisfied, with integers $p,q$, correspond, under rotation with the open book, to resonant $2$-tori obtained from $L_{x,\underline{\eta}_{fix}}$ by forgetting the $X_K$-direction, with frequency vector $(p,q) \in \mathbb{Z}^2$. Note that the relative angular rotation from page to page (as explained in the last sentences in the proof of Lemma \ref{app:lemmarotKepler}) induces the slope of orbits on the $2$-tori. The points in all other non-resonant invariant circles, i.e.\ when $T(c-L)$ is an irrational multiple of $2\pi$, have dense orbits. These circles correspond to non-resonant $2$-tori under the open book rotation, some of which have Diophantine behavior (when the frequencies are viewed as vectors in $\mathbb{R}^2$). 

%Therefore, choosing $\underline{\eta}_{fix}$ appropriately and varying $x$, we may achieve annuli $A_{\underline{\eta}_{fix}}$ for which the dynamics alternates between rational/irrational/Diophantine behavior along each of its circles.

\medskip

\textbf{Lefschetz fibration.} We may further find a very natural foliation of the page $P$ by invariant annuli where $R_r$ acts as a twist map in the classical integrable sense (so the situation is indeed compatible with the classical Poincar\'e-Birkhoff theorem). Indeed, the return map preserves the fibers of (an isotoped version of) the standard Lefschetz fibration $\pi$ on $P=\mathbb{D}^*S^2$, whose regular fibers are symplectic annuli, with precisely two singular fibers whose singularities are $x_\pm$, and such that the boundary of all these annuli coincides with the direct/retrograde planar circular orbits $\gamma_{dir},\gamma_{ret}$. This, which we have stated as Theorem \ref{prop:integrablecase} in the Introduction, can be seen as follows.

\begin{figure}
    \centering
    \includegraphics[width=1.12 \linewidth]{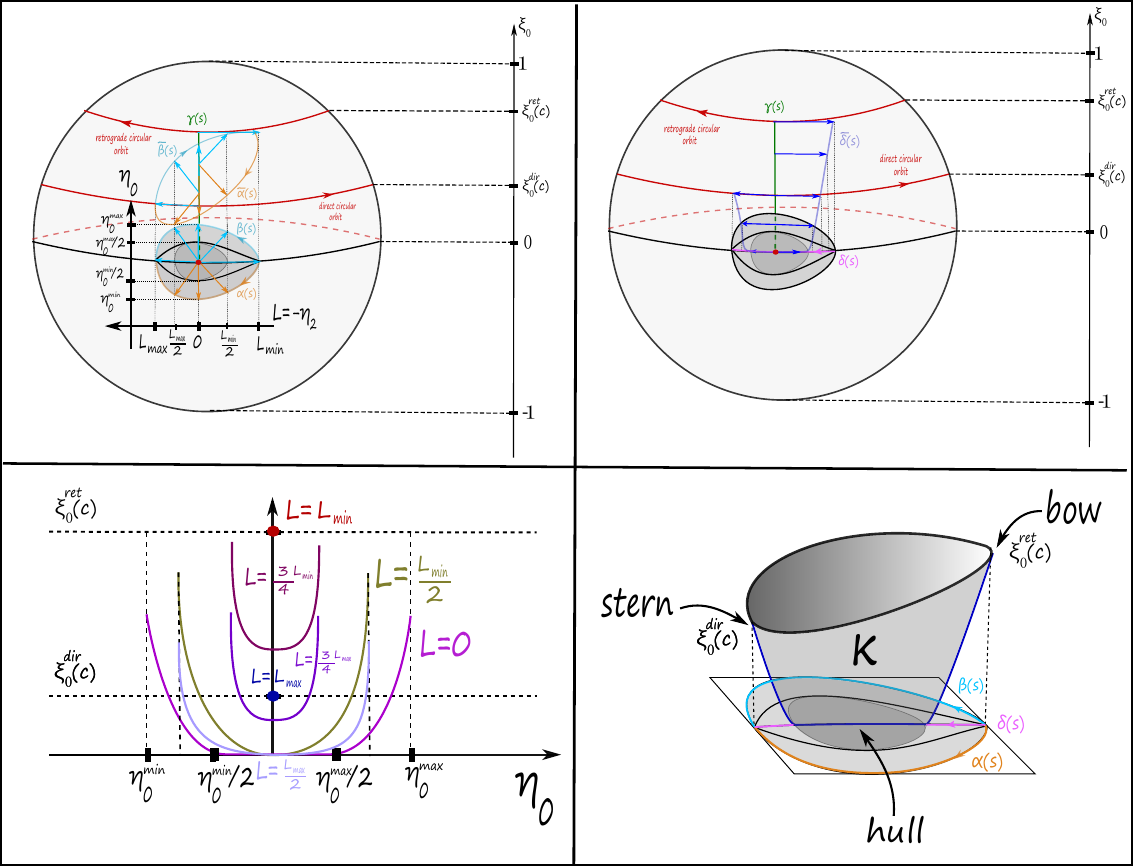}
    \caption{The modified Lefschetz fibration. In the lower right we describe the bump function $\kappa$, shown in the lower left as a $1$-parameter family of functions $\kappa_L$. In the upper left we describe how the modified Birkhoff annulus (in skyblue, parametrized by $\widetilde \beta$) and its opposite (in orange, parametrized by $\widetilde{\alpha}$) look like; in the upper left, how the modified central fiber looks like (in blue, parametrized by $\widetilde \delta$). In green, we have the piece of meridian $\gamma$.}
    \label{fig:modLF}
\end{figure}

First of all, we recall a couple of facts.
\begin{proposition}
The smooth affine variety 
$V=\{ (z_0,\ldots, z_n) \in \C^{n+1} ~|~\sum_{j=0}^n z_j^2 =1 \}$
with its natural K\"ahler form is symplectomorphic to $(T^*S^n,d\lambda_{can} )\subset T^*\R^{n+1}$ via the symplectomorphism
\[
\psi: V \longrightarrow T^*S^n,
z=x+iy \longmapsto (\Vert x\Vert^{-1} x;\Vert x\Vert y)=:(\xi;\eta).
\]
\end{proposition}
Indeed, we can simply pullback the Liouville form to see this. The smooth variety $V$ has a natural Lefschetz fibration obtained by projecting to one of the coordinates,
\[
\Theta_0: V \longrightarrow \C,\quad
(z_0,\ldots,z_n) \longmapsto z_0.
\]
Indeed, we have precisely two critical points of holomorphic Morse type, i.e.\ $p_\pm=(\pm 1,0,\ldots,0)$.
We can pullback this Lefschetz fibration to $T^*S^n$ via $\Theta_0':=\Theta_0\circ\psi^{-1}$, which has the expression
$$
\Theta_0'(\xi;\eta)=\sqrt{\frac{1+\sqrt{1+4\Vert\eta\Vert^2}}{2}} \xi_0+i \sqrt{ \frac{2}{1+\sqrt{1+4\Vert\eta\Vert^2}}} \eta_0.
$$

We now consider the case $n=2$. In order to pullback this Lefschetz fibration to $P$, we first make a couple of observations.
The page $(P,\omega\vert_P=d\alpha|_P)$ has an exact symplectic form that degenerates on the boundary and finite volume. In particular, it is not symplectomorphic to $T^*S^2$.
However, we can modify this form to obtain an ideal Liouville domain (in the sense of \cite{Gir2}):
\begin{proposition}
\label{prop:modified_liouville}
There is  modification $\tilde \alpha$ of $\alpha$ making $(\mbox{int}(P),d\tilde \alpha)$ into an ideal Liouville domain. This modification can be chosen to have the following properties.
\begin{itemize}
\item $\tilde \alpha$ agrees with $\alpha$ in the complement of a collar neighbourhood of the boundary $\nu_P(\partial P)$;
\item $\tilde \alpha = f(\eta_3) \alpha$ for some smooth function $f$ on $\nu_P(\partial P) \setminus \partial P$;
\item $(\mbox{int}(P),d\tilde \alpha)$ is symplectomorphic to $(T^*S^2,d\lambda_{can})$.
\end{itemize}
\end{proposition}
This modification can be constructed with the proof of Proposition~2.19 in \cite{Otto2}. The argument there proves this proposition, after observing that $\eta_3$ is the collar parameter.

Proposition~\ref{prop:modified_liouville} gives us a symplectomorphism $\psi_P:(\mbox{int}(P),d\tilde \alpha) \to (T^*S^2,d\lambda_{can})$.
As a result, we obtain a Lefschetz fibration
$$
\pi_0:(\mbox{int}(P),d\tilde \alpha) \longrightarrow \C,
\quad
p \longmapsto (\Theta_0\circ \psi^{-1}\circ \psi_P)(p).
$$
See Figure \ref{fig:LF}. Note that $R_r$ preserves the fibers of $\Theta_0'$, i.e.\ $\Theta_0'\circ R_r=\Theta_0'$. Moreover, since the modification of Proposition \ref{prop:modified_liouville} only involves the $\eta_3$ coordinate, the same holds for the fibers of $\pi_0$, i.e.\ $\pi_0 \circ R_r=\pi_0$. Topologically, the effect of passing from the Lefschetz fibration on $(\mbox{int}(P),d\tilde\alpha)\cong (T^*S^2,d\lambda_{can})$ to $(P,d\alpha\vert_P)$, $P\cong \mathbb D^*S^2$, may be thought of as projecting along the Liouville direction in the symplectization (see Figure \ref{fig:LF}). With this in mind, the fibers of $\pi_0$ have (ideal) boundary the equator traversed in both directions, which are invariant circles for the planar problem, but these circles are not necessarily closed orbits (this happens e.g.\ in the limit case $c\rightarrow -\infty$ in which we recover the Kepler problem, and $R_r$ is the identity). But we will modify $\pi_0$, relative a neighbourhood of the zero section, in such a way that the boundary of the modified symplectic fibers becomes the disjoint union of the retrograde/direct circular planar orbits $\gamma_{dir}, \gamma_{ret}$, so that the fibers are still symplectic, and invariant under $R_r$. In particular, the modification will happen away from the nodal singularities, so that they are still quadratic.

Both the orbits $\gamma_{dir}, \gamma_{ret}$ are circles with constant $\xi_0$ value, but the values differ for each of them, and furthermore, these values also depend on the energy $c$. In other words, they lie on separate parallels of $P$. Call these values $\xi_0^{ret}(c),\xi_0^{dir}(c)$. We have the inequality
$$
\xi_0^{ret}(c)\geq \xi_0^{dir}(c)\geq 0.
$$
Here is a computation to see this.

\medskip

\noindent
\textbf{Circular orbits in the rotating Kepler problem.}
Write the Hamiltonian in polar coordinates, see for example the appendix of \cite{AFFvK}. This is
$$
H=\frac{1}{2}\left(p_r^2+\frac{L^2}{r^2}\right)^2-\frac{1}{r}+L,
$$
where we have used coordinates $(p_r,r,L,\phi)$ with Liouville form $p_r dr +Ld\phi$.
By writing out the equations of motion, we find a circular orbit must satisfy $r=L^2$.
Substituting this condition into the energy $H=c$, we get two equations. Namely, $-\frac{1}{2r}-\sqrt{r}=c$
for the direct orbit and the orbit in the unbounded component.
For the retrograde orbit, we have the equation $-\frac{1}{2r}+\sqrt{r}=c$.
Rewriting this equation leads to the cubic equation
$$
r^3-c^2r^2-cr-1/4=0,
$$
which we solve with Cardano's formula. 
The $r$-component of orbit in the unbounded component is given by
$$
r_{unbded}=
\left( 1/6\,\sqrt [3]{-54-8\,{c}^{3}+6\,\sqrt {24\,{c}^{3}+81}}+2/3\,
{\frac {{c}^{2}}{\sqrt [3]{-54-8\,{c}^{3}+6\,\sqrt {24\,{c}^{3}+81}}}}
-c/3 \right) ^{2}
$$
The $r$-component of the retrograde orbit is given by
$$
r_{ret}=\left( 1/6\,\sqrt [3]{54+8\,{c}^{3}+6\,\sqrt {24\,{c}^{3}+81}}+2/3\,{
\frac {{c}^{2}}{\sqrt [3]{54+8\,{c}^{3}+6\,\sqrt {24\,{c}^{3}+81}}}}+c
/3 \right) ^{2},
$$
and the $r$-component for the direct orbit is given by 
\[
\begin{split}
r_{dir}&=
{c}^{2}- \left( 1/6\,\sqrt [3]{-54-8\,{c}^{3}+6\,\sqrt {24\,{c}^{3}+81
}}+2/3\,{\frac {{c}^{2}}{\sqrt [3]{-54-8\,{c}^{3}+6\,\sqrt {24\,{c}^{3
}+81}}}}-c/3 \right) ^{2}\\
&\phantom{=}
- \left( 1/6\,\sqrt [3]{54+8\,{c}^{3}+6\,
\sqrt {24\,{c}^{3}+81}}+2/3\,{\frac {{c}^{2}}{\sqrt [3]{54+8\,{c}^{3}+
6\,\sqrt {24\,{c}^{3}+81}}}}+c/3 \right) ^{2}
\end{split}
\]
At the critical value $c=-3/2$, we have
$$
r_{unbded}=r_{dir}=1,
\quad
r_{ret}=1/4.
$$
The corresponding values for the norms of momenta are
$$
p_{dir}(c=-3/2)=1,
\quad
p_{ret}(c=-3/2)=2,
$$
and using the derivative, and recalling that $\xi_0=\frac{\Vert p \Vert^2-1}{\Vert p \Vert^2+1}$, we can verify the above claim.

%The direct circular orbit $\gamma_{dir}$ actually never attains the value $\vert p\vert=1$ (i.e.\ $\xi_0=0$) on the regular energy level sets $c<-3/2$, although $\vert p\vert=1$ is attained for the singular energy level set $c=-3/2$, i.e.\ $\xi_0^{dir}(-3/2)=0$. For the retrograde orbit $\gamma_{ret}$, the value $\vert p\vert=1$ is not attained until the much higher energy $c=-1/2$ (which we will not consider), i.e.\ $\xi_0^{ret}(-1/2)=0$. \textcolor{red}{Otto, please add some details here if possible, i.e.\ of the relevant computation for the circular orbits} 

\medskip

\noindent
\textbf{Deforming cylinders.}
Parametrize a piece of the meridian joining $\underline{\xi}_{fix}=(0,1,0)$ to $\underline{\xi}_{ret}:=(\xi_0^{ret}(c),\sqrt{1-(\xi_0^{ret}(c))^2},0)$ via the path
$$
\gamma: [0,\xi_0^{ret}(c)]\rightarrow S^2,
$$
$$
\gamma(s)=(s,\sqrt{1-s^2},0),
$$
and consider the time-$s$ parallel transport $PT_s: T_{\underline{\xi}_{fix}}S^2\rightarrow T_{\gamma(s)}S^2$ along $\gamma$ (with respect to the round metric). Recall that $L\vert_{T_{\underline{\xi}_{fix}}S^2}=-\eta_2$, and hence is constant in the vertical $\eta_0$-lines. Moreover, $L\vert_{T_{\gamma(s)}S^2}=-(\sqrt{1-s^2})\eta_2$, so the analogous statement holds over $\gamma(s)$. Denote by $P_{\underline{\xi}}$ the fiber of $P$ over $\underline{\xi}\in S^2$. We view $P_{\underline{\xi}_{fix}}\subset T_{\underline{\xi}_{fix}}S^2$ as a closed star-shaped domain in the $(\eta_0,L)$-plane (see Figure \ref{fig:modLF}).  Let 

$$
L_{min}=\min( L\vert_{P_{\underline{\xi}_{fix}}})=\max(\eta_2\vert_{P_{\underline{\xi}_{fix}}}),\; L_{max}=\max( L\vert{P_{\underline{\xi}_{fix}}})=\min(\eta_2\vert_{P_{\underline{\xi}_{fix}}}),
$$
$$
\eta_0^{min}=\min(\eta_0\vert_{P_{\underline{\xi}_{fix}}}),\;\eta_0^{max}=\max(\eta_0\vert_{P_{\underline{\xi}_{fix}}})=-\eta_0^{min},
$$
$$
L_{min}^{dir}=\min( L\vert_{P_{\gamma(\xi_0^{dir}(c))}})=\left(\sqrt{1-(\xi_0^{dir}(c))^2}\right)L_{min},$$
$$
L_{max}^{ret}=\max( L\vert_{P_{\gamma(\xi_0^{ret}(c))}})=\left(\sqrt{1-(\xi_0^{ret}(c))^2}\right)L_{max}.$$
Then $L_{min}^{dir}, L_{max}^{ret}$ are respectively the angular momenta of the $\gamma_{dir}$ and $\gamma_{ret}$. We then define a smooth bump function
$$
\kappa: P_{\underline{\xi}_{fix}} \rightarrow [0,\xi_0^{ret}(c)] 
$$
$$
\kappa(\eta_0,L)=\kappa_L(\eta_0),
$$
shown qualitatively in Figure \ref{fig:modLF}, where we also describe it as a $1$-parameter family $\kappa_L$ for clarity; intuitively, $\kappa$ has the shape of a ``boat''. In particular, we impose that $\kappa$ vanishes on $\frac{1}{2}P_{\underline{\xi}_{fix}}$ (the ``hull'' of the boat), that it satisfies the symmetry $\kappa(\eta_0,L)=\kappa(-\eta_0,L)$, and that $\kappa(0,L_{min})=\xi_0^{ret}(c)$ (the ``bow''), $\kappa(0,L_{max})=\xi_0^{dir}(c)$ (the ``stern''). We then let
$$
PT: P_{\underline{\xi}_{fix}} \rightarrow P,
$$
$$
PT(\eta_0,L)=PT_{\kappa_L(\eta_0)}(\eta_0,L),
$$
and denote $\widetilde P_{\underline{\xi}_{fix}}=PT(P_{\underline{\xi}_{fix}})$ the (diffeomorphic) image of this map. We now foliate $P_{\underline{\xi}_{fix}}$ by line segments joining $(0,L_{min})$ to $(0,L_{max})$, so that the resulting foliation looks like an ``eye'' (see Figure \ref{fig:modLF}). If $l:[L_{min},L_{max}]\rightarrow P_{\underline{\xi}_{fix}}$ is a non-singular parametrization of such a segment, we denote $\widetilde l=PT(l)\subset \widetilde P_{\underline{\xi}_{fix}}$, which taken together provide a foliation of $\widetilde P_{\underline{\xi}_{fix}}$. Figure \ref{fig:modLF} shows what this path looks like in the case of the upper boundary of $P_{\underline{\xi}_{fix}}$ (denoted $l=\beta$), the lower one ($l=\alpha$), and the central segment ($l=\delta$). If $\phi_L^t: P\rightarrow P$ denotes the time-$t$ Hamiltonian flow of $L$ (i.e.\ rotation along the $(\xi_0,\eta_0)$-axis), for $l$ a segment in the foliation, we define
$$
\widetilde C_l=\bigcup_{t\in \mathbb{R}} \phi_L^t(\widetilde l).
$$
By construction, this is a cylinder with boundary $\gamma_{ret}\cup \gamma_{dir}$. Moreover, note that each $l$ is by construction positively transverse to the level sets of $L\vert_{T_{\underline{\xi}_{fix}}S^2}$, and the parallel transport map $PT_s$ preserves the vertical level sets of $L$. Therefore $d\tilde \alpha(\partial_t\widetilde{l}, X_L)=dL(\partial_t \widetilde{l})>0$, which means that $\widetilde C_l$ is a symplectic cylinder. Note that $\widetilde C_l$ is invariant under $R_r$, since it is a union of the invariant circles of the form $C_{x,\underline{\eta}_{fix}}$ considered above. We then flow the $\widetilde C_l$ with the action of $K$, i.e.\ we consider
$$
\widetilde C_l^t=\phi_K^t(\widetilde C_l), \; t\in \mathbb{R}.
$$
These are symplectic cylinders which glue with the fibers of the standard Lefschetz fibration on $(1/2)P$, a neighbourhood of the zero section. Note that as constructed the $\widetilde C_l^t$ do not intersect the zero section away from the equator, but after gluing them with the standard fibers, they do. After this modification, the $\widetilde C_l^t$ are the fibers of a Lefschetz fibration $\pi$ which coincides with $\pi_0$ along $(1/2)P$, and has the desired properties. This map is given by
$$
\pi:P \rightarrow \mathbb{C},
$$
$$
\pi=\pi_0 \circ \Psi_1^{-1},
$$
where $\Psi_t: P \rightarrow P$, $t\in [0,1]$, is the smooth isotopy induced by doing the same deformation construction with $\kappa$ replaced with $t\kappa$, so that $\Psi_0=id_P$, and $\Psi_t$ is the identity on $(1/2)P$ for all $t$.

By construction, the return map acts as a twist map on each fiber, preserving the horizontal circles where it acts as a rotation whose angle changes from circle to circle, rotating each boundary component ($\gamma_{dir}$ resp.\ $\gamma_{ret}$) with angle $T(c-L^{dir}_{min})$ resp. $T(c-L^{ret}_{max})$. The dynamics then alternates between rational/irrational/Diophantine behavior along each circle in the fiber. The critical points of $\pi$ are precisely $x_\pm$, which are fixed by the return map. This is depicted in Figure \ref{fig:LF}, where we view the fibers as copies of $T^*S^1$ which become asymptotic at infinity to the direct/retrograde orbits (we only draw their intersection with $\mathbb{D}^*S^2)$, after an isotopy of $\pi$ so that it looks standard. The pages of the standard open book in $\mathbb{R}P^3=\partial \mathbb{D}^*S^2$ are obtained as the radial projections along the Liouville direction of the fibers of $\pi$. In the second picture of Figure~\ref{fig:LF} we sketch how some of these fibers look like from an alternative perspective, for the case of the Kepler problem, where we have projected all fibers to $\mathbb{D}^*S^2$ so that now they are copies of $\mathbb{D}^*S^1$ (including the Birkhoff annulus consisting of co-vectors which point towards the upper hemisphere along the equator, corresponding to the segment $\beta$; its opposite version, consisting of co-vectors which point towards the lower hemisphere along the equator and corresponding to $\alpha$; and the central fiber, a regular annulus which intersects the Lagrangian zero section along the equator, and corresponds to $\delta$).

\medskip

\textbf{Outlook and further comments.} Theorem~\ref{thm:gss_openbooks} in particular provides a global hypersurface of section for any nearby perturbed system (where $\mu$ is close to $1$). These can be thought of as the original hypersurface $P$, but where the disks in the fiber are now perturbed to have boundary the level set for the corresponding planar problem. From the above discussion one can ask what part of the above structure survives a perturbation of $\mu$, and/or whether we can detect Arnold diffusion. This is naturally the realm of KAM theory, as well as weaker versions like Aubry-Mather theory, and we shall not pursue this direction here. In work of the first author \cite{M20}, it will be argued that the underlying geometric structure for the page (i.e.\ the symplectic Lefschetz fibration) also holds non-perturbatively; this will be used to define an associated Reeb dynamics on $S^3$ whenever the planar dynamics is dynamically convex. 

\section{Symplectic monodromy and return maps}\label{monodromyvsreturnmap}

In this appendix, we prove a general fact that implies that the return map in the statement of Theorem \ref{thm:returnmap} is Hamiltonian. Namely, we shall establish that return maps arising from an adapted Reeb flow, under a suitable concavity assumption near the boundary, are always symplectically isotopic to a representative of the monodromy, via an isotopy which preserves the boundary (Lemma \ref{lemma:monodromyvsreturnmap}). If the monodromy happens to be Hamiltonian (as is the case for an even power of the Dehn-Seidel twist), and the page $P$ happens to have $H^1(P;\mathbb{R})=0$ (so that all symplectic isotopies are Hamiltonian), it follows that the return map is Hamiltonian.

Fix a concrete open book $(B,\theta)$ on a closed $(2n+1)-$manifold $\Sigma$, supporting a contact structure $\xi$. Assume that $\alpha$ is a Giroux form, and let $f: P \rightarrow P$ be the associated Poincar\'e return map on a fixed page, which is an element of $\mathrm{Symp}(\mbox{int}(P),d\alpha|_P)$. We assume that $f$ admits a (unique) smooth extension to the boundary $B$, and we will denote the symplectic monodromy of the open book by $[\phi]$.

We need some notation to state an additional assumption.
The binding $B$ is a contact submanifold, so by a neighbourhood theorem $\alpha=A(\alpha_B+r^2 d\theta)$. 
Because $B$ is invariant under the Reeb flow of $\alpha$, whose Reeb vector field coincides with that of $\alpha_B$ along $B$, we see that
\begin{itemize}
\item[(i)] $dA|_{B}=0$ (including the normal direction);
\item[(ii)] $A|_B=1$.
\end{itemize}
We will consider the Hessian along the binding below, and make the assumption that the Hessian of $A$ in the normal direction of $B$ is negative definite.
To see that this condition is independent of the chosen trivialization of the normal bundle, simply write out the second derivative of $A$, and apply the chain rule; the term not corresponding to the Hessian vanishes due to the assumption that $dA$ vanishes along $B$.

\begin{lemma}\label{lemma:monodromyvsreturnmap}
Assume that the Hessian of $A$ in the normal direction of $B$ is negative definite.
Then the return map $f$ is symplectically isotopic to a representative $\phi$ of the symplectic monodromy, as elements in $\mathrm{Symp}(\mbox{int}(P),d\alpha|_P)$. The isotopy $\psi_t$ is supported near $B$, extends smoothly to $B$, and preserves $B$. 
\end{lemma}

\begin{proof} For clarity, we split the proof in several steps. 

\medskip

\textbf{Step 1: Isotopy to standard form.} 
We shall consider $\alpha_t$ the $1$-parameter family of Giroux forms following \cite[Prop.\ 3.1]{DGZ}
This family joins the original contact form $\alpha=\alpha_0$ to an adapted contact form $\alpha_1$ in standard form near $B$, and satisfying $\alpha_t=\alpha$ away from a fixed neighbourhood of $B$, for all $t$. 
This means that $\alpha_1=h_1(r)\alpha_B+h_2(r)d\theta$ in a neighbourhood $B\times \mathbb{D}^2$, for suitable functions $h_1,h_2$ and $\alpha_B=\alpha\vert_B$, where $r$ is the radial coordinate so that $B=\{r=0\}$, and $\theta$ is the open book coordinate. For convenience of the reader, we recall some of the arguments in \cite[Prop.\ 3.1]{DGZ}. 

Using the characteristic foliation on the pages, one first finds a neighbourhood $B\times \mathbb{D}^2$ in which $\alpha=A(\alpha_B+r^2 d\theta)$ for a smooth positive function $A$ satisfying $A\equiv 1$ along $r=0$ and $\partial_rA<0$ for $r>0$ (note that this last condition is to ensure that $d\alpha$ is positive on each page). From condition (i), we have in particular that $\partial_rA\vert_{r=0}=0$. 
Since $Hess(A)$ is negative definite along $B$, we find a neighbourhood $r<\delta$ where this property also holds.
On this neighbourhood we also find a constant $a>0$ with $a< \min (-Hess(A))$.
After that we choose a decreasing cutoff function $\lambda$ satisfying $\lambda \equiv 1$ near $r=0$ and $\lambda=0$ for $r=\delta<1$.
Define the deformed form $\beta$ by
$$
\beta:=(\lambda(r)(1-a r^2)+(1-\lambda(r))A)(\alpha_B+r^2d\theta).
$$
We claim that $\beta$ is a contact form. 
This is evident for $r\geq \delta$, as $\beta|_{r\geq \delta}=\alpha$.
For $r<\delta$, we note that $ \tilde A:=(\lambda(r)(1-a r^2)+(1-\lambda(r))A)>0$, and for $0<r<\delta$ it suffices to check that $\partial_r \widetilde A<0$. We have 
\begin{equation}
    \begin{split}
       \partial_r \tilde A 
&=\lambda'(1-ar^2)-2ar\lambda-\lambda' A+(1-\lambda) \partial_r A\\
&=\lambda'(-Hess(A)(r,r)-ar^2 +o(r^2) )-2r\lambda +(1-\lambda)\partial_r A, 
    \end{split}
\end{equation}
The last two terms are negative for $r>0$, and for the first term we observe that the expression
$$
-Hess(A)(r,r)-ar^2 +o(r^2)
$$
is positive by our choice of $a$. As $\lambda'\leq 0$, the claim follows.
Hence the linear interpolation of these forms
$$
\alpha_t=(1-t)\alpha+t\beta=(t\tilde A+(1-t)A)(\alpha_B+r^2d\theta)
$$
brings $\alpha$ to the desired form, through a family of Giroux forms. It follows that $\alpha_t\vert_P=F_t\alpha_B$ near $B$, with $F_t:=t\tilde A+(1-t)A$ satisfying $\partial_r F_t<0$ for $r>0$ and $F_t\vert_B=1$. Note that we may further choose $h_1(r)=1-ar^2$ for some $a>0$, and $r\leq \delta$. We may also slightly isotope $h_2$ so that $h_2(r)=r^2$ for $r\leq \delta$. 
%However, we wish to have $\alpha_t\vert_B=\alpha_B$ for all $t$. To arrange this, we scale $\alpha$ by a sufficiently small constant if necessary before considering $\beta$ (this does not change its Reeb dynamics, and therefore $f$ stays fixed), and so we make $A$ uniformly $C^\infty$-small. Then we may still arrange that the $r$-derivative of $A^\prime$ is negative even if we choose $c=1$, and the proof still goes through; and so this isotopy $\alpha_t$ can be taken fixed at $B$ and we may fix $c=1$.
%we simply scale $\alpha_t$ by $1/c_t$, and so we assume $F_t\vert_B=1$ for all $t$.
% 

\medskip

\textbf{Step 2: Moser trick.} We let $\omega_t:=d\alpha_t\vert_P$. This is then a $1$-parameter family of $2$-forms on $P$ which are symplectic on int$(P)$ (but become degenerate at $B$). We wish to appeal to the standard Moser trick to obtain a symplectic isotopy $\psi_t$ of int$(P)$, satisfying $\psi_t^*\omega_t=\omega_0$. For this, because of degeneracy of $\omega_t$ at $B$, we need to study the behavior of the Moser isotopy at the boundary. This isotopy would be generated by a vector field $X_t$ tangent to $P$ and defined via $i_{X_t}\omega_t=-\dot{\alpha_t}\vert_P$, which would make $\psi_t$ supported near $B$. 
Moreover, $\psi_t$ would need to be tangent to the $r$-direction: indeed, near $B$ we have $\dot{\alpha_t}\vert_P=(\tilde A-A)\alpha_B$ and $\omega_t=dF_t \wedge \alpha_B + F_td\alpha_B$, so this uniquely determines $X_t=\gamma_t \partial_r$ with $\gamma_t=\frac{A-\tilde A}{\partial_rF_t}=\frac{A-\tilde A}{t\partial_r\tilde A +(1-t)\partial_rA}$. 

To make sure its flow is well-defined at $B$, we study the limit of $\gamma_t$ as $r\rightarrow 0$. Expand numerator and denominator in powers of $r$.
Near $r=0$, we have $\lambda \equiv 1$, so this yields
$$
\gamma_t=\frac{1+Hess(A)(r,r)+o(r^2)-1+ar^2}{-2art +(1-t)\partial_r^2 A \cdot r +o(r)}
=\frac{r\partial_r^2 A/2+ar+o(r)}{-2at +(1-t)\partial_r^2 A  +o(1)}.
$$
The denominator is negative and can be bound from above by $-2a$. The numerator vanishes for $r=0$, so we conclude that the Moser flow exists and it extends as the identity to $B$.

\medskip

\textbf{Step 3: Holonomy maps.} 
Let $f_t: (\mbox{int}(P),\omega_t)\rightarrow (\mbox{int}(P),\omega_t)$ be the symplectic holonomy map associated to the symplectic connection $\omega_t$, i.e.\ the return map associated to a vector field spanning the $1$-dimensional horizontal distribution Hor$_t=\{v \in T(\Sigma\backslash B): i_v \omega_t=0\}$. Define $g_t:=\psi_t^{-1} \circ f_t \circ \psi_t$. Then $g_t$ is a symplectic isotopy in Symp$(\mbox{int}(P),\omega_0)$, i.e.\ for the same symplectic form on int$(P)$, joining (rel $B$) the original return map $f=f_0=g_0$ to a symplectomorphism $g_1$ for $\omega_0$. In particular, $g_t$ (and $f_t$) extends smoothly to $B$ for every $t$, agreeing with the extension of $f$. The $\omega_{1}$-horizontal distribution is spanned by $$
X_{1}=\partial_\theta -\frac{h_2^\prime}{h_1^\prime}R_B,
$$
and so $X_{1}=\partial_\theta +\frac{1}{a}R_B$ for $r\leq \delta$. 

\medskip

\textbf{Step 4: Monodromy.} Consider the symplectic form $\omega_{g}=d( \frac{\alpha_1}{g^2})$, where $g$ is a function satisfying $g=r$ near $r=0$, and $g=1$ for $r\geq \delta-\epsilon$ for some small $0<\epsilon\ll \delta$. The $\omega_{g}$-horizontal distribution is spanned by $$
X_{g}=\partial_\theta -2r\left(\frac{g-rg^\prime}{gh_1^\prime-2h_1g^\prime}\right)R_B=:\partial_\theta+\Phi(r)R_B.
$$
Note that $\Phi\equiv 0$ near $r=0$, and $\Phi\equiv 1/a$ on $r\geq \delta-\epsilon$, so that $X_1=X_g$ on $r\geq \delta-\epsilon$ and also on $P\backslash\{r\in [0,1]\}$. By definition, the return map for $\omega_g$ is a representative of the monodromy (see e.g.\ the proof of \cite[Prop.\ 2.19]{Otto2}). Note that $\omega_g$ has infinite volume since $g$ vanishes at $r=0$; it is in fact symplectomorphic to the Liouville completion of $(P_{r_1},\omega_1\vert_{P_{r_1}})$, where $P_{r_1}=\{r\geq r_1\}\subset P$ with $0<r_1< \delta$, which is a Liouville domain filling the contact manifold $(B,\xi_B=\ker \alpha_B)$ (whose completion is independent of such $r_1$). 

\medskip

\textbf{Step 5: Square root map.} We now make the symplectic manifold $(P,\omega_1)$ into an honest Liouville domain by considering a square root map $Q:P\rightarrow P$, defined as follows. Choose a continuous function $q:[0,1]\rightarrow \mathbb{R}$ which satisfies $q(r)=\sqrt{r/a}$ for $r\leq \delta-\epsilon<1$, $q(r)=r$ for $r\geq \delta$, and $q$ is smooth away from $r=0$ where it satisfies $q^\prime >0$. Define $Q(b,r,\theta):=(b,q(r),\theta)$ along $B\times \mathbb{D}^2$, and extend it to $P$ via the identity. Note that this map is not smooth at $r=0$ (only continuous), although its inverse is. Define $\omega_Q:=Q^*\omega_1$ on $P$. Note that along $B\times \mathbb{D}^2$ we have $\omega_Q=d(h_1(q(r))\alpha_B),$ and so $\omega_Q=d((1-r)\alpha_B)$ for $r\leq \delta-\epsilon$. It follows that $\omega_Q$ is a symplectic form which is also non-degenerate at $B$, and so $(P,\omega_Q)$ is a Liouville domain filling the \emph{strict} contact manifold $(B,\alpha_B)$, and having the same symplectic volume as $(\mbox{int}(P),\omega_1)$.

\medskip

\textbf{Step 6: Continuous conjugation to a Liouville domain.} By squashing the $r$-direction, we can construct a symplectomorphism between the (honest) Liouville domains 
$$
F:(P,\omega_Q)\rightarrow (P_{r_1},\omega_g),
$$
where the precise value of $r_1$ is determined by the total volume of $(P,\omega_Q)$. We put $$F_Q:=F\circ Q^{-1}:(P,\omega_1)\rightarrow (P_{r_1},\omega_g),$$ and $G:= F_Q\circ \psi_1$, which satisfies $G^*\omega_g=\omega_0$. Note that $G$ has a smooth extension to $B$, although $G^{-1}$ extends only continuously.

\medskip

\textbf{Step 7: the symplectic isotopy.} The vector field $X_{g}$ differs from $X_1$ by $$X_g-X_1=\left(\Phi+\frac{h_2^\prime}{h_1^\prime}\right)(r)R_B=:\Psi(r)R_B,$$ with $\Psi\equiv 0$ on $r\geq \delta-\epsilon$ and on $P\backslash\{r\in [0,1]\}$, and $\Psi\equiv 1/a$ near $r=0$. This vector field is Hamiltonian for $\omega_1$, generated by $H:P\rightarrow \mathbb{R}$ which on $B\times [0,1]$ is given by $H(b,r)=H(r)=\int_0^r\Psi(s)h_1^\prime(s)ds$, and extends to $P\backslash\{r\in [0,1]\}$ by zero. Let $\phi_s^H$ be the corresponding Hamiltonian flow (computed with respect to $\omega_1$), which we may take supported near $B$, which moreover extends smoothly to $B$, and preserves $B$. Then $\phi_s:=F_Q \circ \phi_s^H \circ f_1 \circ F_Q^{-1}$ joins $\phi_0=F_Q\circ f_1 \circ F_Q^{-1}$ to $\phi^\prime:=\phi_1$. By construction, $\phi^\prime$ preserves $\omega_g$ on $P_{r_1}$, and extends from $P_{r_1}$ to $P$ via the identity on $r\leq r_1$, to a representative of the monodromy, which we also call $\phi^\prime$. We have $\phi_s^{H \circ \psi_1}=\psi_1^{-1}\circ \phi^H_s\circ \psi_1$ is Hamiltonian for $\omega_0$. Finally, the symplectic isotopy $g_s^\prime:=\phi_s^{H \circ \psi_1}\circ g_1= \psi_1^{-1}\circ \phi_s^H \circ f_1 \circ \psi_1$ joins $g_1$ to $\phi:=G^{-1}\circ \phi^\prime \circ G$, which also represents the monodromy. Concatenating $g_s^\prime$ to $g_t$ we join the original return map $f=g_0$ to $\phi$, as elements in Symp$(\mbox{int}(P),\omega_0)$. But note that $g_s^\prime$ does \emph{not} involve the function $G^{-1}$, which was only used to define $\phi$ (which agrees with the identity at $B$), and in fact has a \emph{smooth} extension to $B$. This finishes the proof.
\end{proof}

\end{document}